\documentclass[final,3p,8pt]{elsarticle}

\usepackage{amssymb}
\usepackage{array}
\usepackage{mathrsfs}
\usepackage{amscd}
\usepackage{amsmath}
\usepackage{amsfonts}
\usepackage{amsthm}
\usepackage{psfrag}
\usepackage{textcomp}
\usepackage{bm}
\usepackage{dsfont}
\usepackage{algorithm}
\usepackage{algpseudocode}

\journal{arXiv.org}

\usepackage{amssymb}

\biboptions{square,sort&compress}

\usepackage[figuresright]{rotating}
\usepackage{moreverb}
\usepackage{amssymb}
\usepackage{array}
\usepackage{mathrsfs}
\usepackage{amscd}
\usepackage{amsmath}
\usepackage{amsfonts}
\usepackage{amsthm}
\usepackage{psfrag}
\usepackage{textcomp}
\usepackage{bm}
\usepackage{color}
\usepackage{float}
\usepackage[latin1]{inputenc} 	
\usepackage{afterpage}
\usepackage{tabularx}
\usepackage{booktabs}
\usepackage{romannum}
\usepackage{multirow}
\usepackage[normalem]{ulem}

\newcommand{\ds}	{\displaystyle    }

\newcommand{\p}[1]{\left( #1 \right)}
\newcommand{\cor}[1]{\left[ #1 \right]}

\newcommand{\Om}	{\Omega   }

\newcommand{\om}	{\omega   }

\newcommand{\dt}	{\Delta t   }
\newcommand{\ep}	{\epsilon   }

\newcommand{\mbf} {\mathbf  }

\newcommand{\mrm} {\mathrm  }

\newcommand{\abs}[1]{\left| #1 \right|}

\newcommand{\mbR} {\mathbb{R}}

\newcommand{\mbu}{\mathbf{u}}

\newcommand{\bphi}[1] {\bm{\phi}_{#1}}

\newcommand{\norm}[1]{\left\| #1 \right\|}

\newcommand{\cccolumna}[3]{\left\{\begin{array}{c} #1 \\ #2 \\ #3 \end{array}\right\}}
\newcommand{\mmatriz}[4]{\left[\begin{array}{cc} #1 & #2 \\ #3 & #4 \end{array}\right]}

\newcommand{\der}[2]{\frac{\mrm{d} #1}{\mrm{d} #2}}

\newcommand{\mbPhi}{\mathbf{\Phi}}
\newcommand{\mbM}{\mathbf{M}}
\newcommand{\mbK}{\mathbf{K}}

\begin{document}

\begin{frontmatter}




\title{Damping perturbation based time integration asymptotic method for structural dynamics}


\author[upv1,upv2]{Mario L\'azaro\corref{cor}}
\ead{malana@mes.upv.es}

\cortext[cor]{Corresponding author. Tel +34 963877000 (Ext. 76732). Fax +34 963877189}
\address[upv1]{Department of Continuum Mechanics and Theory of Structures\\
Universitat Polit\`ecnica de Val\`encia 46022 Valencia, Spain}
\address[upv2]{Instituto Universitario de Matemática Pura y Aplicada\\
Universitat Polit\`ecnica de Val\`encia 46022 Valencia, Spain}

\begin{abstract}

	The light damping hypothesis is usually assumed in structural dynamics since dissipative forces are in general weak with respect to inertial and elastic forces. In this paper a novel numerical method of time integration based on the artificial perturbation of damping is proposed. 
The asymptotic expansion of the transient response results in an infinite series which can be summed, leading to a well-defined explicit iterative step-by-step scheme. Conditions for convergence are rigorously analyzed, enabling the determination of the methodology boundaries in form of maximum time step. The numerical properties of the iterative scheme, i.e. stability, accuracy and computational effort are also studied in detail. The approach is validated with two numerical examples, showing a high accuracy and computational efficiency relative to other methods.

\end{abstract}


\end{frontmatter}


\section{Introduction}


In this paper, the dynamic analysis of linear structures with time invariant properties and viscous damping is considered. Equilibrium of inertia, elastic and damping forces leads to well--known system of second-order differential equations
\begin{equation}
\mbf{M} \, \ddot{\mbf{u}} + \mbf{C} \, \dot{\mbf{u}} + \mbf{K} \, \mbf{u} = \mbf{f}(t) \ , \quad
\mbf{u}(0) = \mbf{u}_0 \ , \quad \dot{\mbf{u}}(0) = \dot{\mbf{u}}_0
\label{eq001}
\end{equation}
where $\mbf{M}, \mbf{C}, \mbf{K} \in \mathbb{R}^{N\times N}$ are the mass, damping and stiffness matrices respectively and $\mbf{f}(t) \in \mathbb{R}^{N}$ represents the time-dependent vector of external forces at time $t$.  The mass matrix is assumed to be symmetric and positive definite, while $\mbf{K}$ and $\mbf{C}$ are positive semidefinite symmetric matrices. The $N$ degrees of freedom (dof) are arranged in the column vector $\mbf{u}(t) \in \mathbb{R}^{N}$.  \\


Time-integrator algorithms have become an essential tool in the numerical simulation of structural dynamics. In general, two families of methods can be distinguished: explicit and implicit schemes~\cite{Dokainish-1989a,Dokainish-1989b}. In the former, the response at each time step is evaluated from the previous step. Implicit integration requires the solution of one or more systems of linear equations at each iteration. In linear problems, the coefficient matrices are time invariant and a proper factorization of these allows to optimize the solution of the response at each step. Most implicit algorithms exhibit unconditional stability: Newmark's method~\cite{Newmark-1959}, Wilson's method \cite{Wilson-1968}, Hilber-Hughes-Taylor's method (HHT) \cite{Hilber-1977}, Wood-Bozak-ZienKiewicz's method (WBZ) \cite{Wood-1980}. The HHT and WBZ approaches were proposed with the aim of introducing an artificial dissipation in order to eliminate the high-frequency part of the response. As consequence, damping inconsistencies arose, which were addressed by the $\alpha$-scheme from Chung and Hulbert~\cite{Chung-1993}. Krenk~\cite{Krenk-2006} presented a full energy analysis of Newmark based time integration algorithms. The Bathe's method~\cite{Bathe-2005,Bathe-2007,Bathe-2012b} proposes a composite time-integration approach, splitting up the time-step into two sub-steps, so that the trapezoidal rule is used for the first sub-step and the 3-point Euler backward method for the second one. The Bathe's method has been recently improved with the introduction of new parameters as a controllable spectral radius~\cite{Bathe-2019b,Bathe-2019c} or a prescribing desired numerical dissipation~\cite{Bathe-2019a}.  \\

The computational effort of implicit methods increases as the number of degrees of freedom of the system becomes larger. In such cases, explicit methods are in general more competitive. However, a known drawback of explicit methods is their conditional stability, such as the central difference method (CDM) or the 4th order Runga-Kutta method (RK4)\cite{Butcher-2008,Burden-2001}. This is partially overcome by the family of precise integration methods (PIM) of special interest in linear models of structural dynamics. The main contribution of the PIM, originally proposed by Zhong and Williams~\cite{Zhong-1994,Zhong-2004}, is the proposal of an explicit algorithm based on the precise computation of the exponential matrix by means of the $2^p$ algorithm in combination with the evaluation of the incremental part~\cite{Fung-2004}. The method exhibits conditional stability although within such a wide range of time steps that in practice it can be considered unconditionally stable. In addition, it shows highly precise results in the homogeneous problem, although the accuracy drops when there exist  non-homogeneous terms in the equations. Zhong and Williams~\cite{Zhong-1994} proposed to linearly interpolate  the external applied force in the time step, which needs the computation of an inverse matrix with the consequent loss of accuracy. Since the first attempt of Zhong and Williams, a variety of proposals to improve the efficiency and accuracy of the method have been published. Lin~\cite{Lin-1995} improved the PIM by developing in trigonometric Fourier series the external force term within the time step. Wang et al.~\cite{WangMF-2005} proposed the so-called Modified Precise Integration Method (MPIM), which consisted of evaluating the non-homogeneous term by using Gaussian quadrature integration. The approach allows achieving high accuracy in the result in exchange for increasing the computational cost. Wang and Au~\cite{WangM-2007} presented a PIM avoiding computation of the inverse matrix for the evaluation of the non-homogeneous terms. They proposed two methods for the evaluation of the independent term based on interpolation with Chebyshev~\cite{WangMF-2008} and Lagrange polynomials~\cite{WangMF-2009a}. Fung~\cite{Fung-1997} used the principle of precise integration, combining the $2^p$ algorithm with the evaluation of the incremental part, to determine the step-response and impulse-response matrices of the system, as well as their derivatives. Fung~\cite{Fung-2006} also proposed an efficient algorithm based on the precise integration, using the Krylov subspace and the Pad\'e approximations. In the work of Zhu and Law~\cite{Zhu-2001}, the PIM is applied for a continuous Euler-Bernouilli beam under moving loads. Caprani~\cite{Caprani-2013} developed a PIM based on modal analysis for moving forces in footbridge vibration response.\\

Perturbation techniques have been applied for the numerical solution of linear and nonlinear, algebraic and differential equations, by expanding the solution in terms of certain parameter which is known a prior to be very small respect to other terms of the model~\cite{Hassan-2004}. Real structures are usually assumed to be lightly damped, that is, dissipative forces can be considered much smaller than inertia and elastic forces. Artificial perturbation of dissipative terms in linear structural dynamics has been successfully applied mainly for the computation of complex frequencies and modes of non-classically viscously~\cite{Meirovitch-1979,Meirovitch-1985,Chung-1986,Da-Silva-1995,Cha-2005,Lazaro-2016a} and non-viscously~\cite{Daya-2001,Daya-2003,Lazaro-2013d,Zoghaib-2014,Lazaro-2015d,Lazaro-2016b,Lazaro-2018a} damped structures. In the time domain, techniques based on asymptotic perturbation are very useful to obtain solutions in nonlinear mechanics~\cite{Gallager-1975,Cochelin-1994,Mei-2008}.  In linear structural dynamics, Fafard et al.~\cite{Fafard-1997} and Berhama-Chekroun~\cite{Fafard-2001} have developed asymptotic methods by expanding the transient response in time power series.  However, hardly any references can be found in the literature on asymptotic time-integration methods for the linear model of Eq. \eqref{eq001}. Perhaps because they intrinsically involve two iterative processes: one advancing in the time domain and the other one along the asymptotic dimension. This fact makes them uncompetitive with other explicit methods that only involve step-by-step iterations in time. 
\\

In this paper, an asymptotic time integration algorithm based on artificial perturbation of damping for linear structural dynamics is proposed. A step-by-step explicit algorithm arises after the sum of the infinite series resulting from asymptotic expansion of the response, thus avoiding double iteration. Conditions for convergence of the method have been explored allowing the determination boundary limits of allowable time-step. Algorithms for the computation of the main matrices of the scheme are developed. Stability, accuracy and computational effort are addressed in detail. The proposed method is validated through two numerical examples, covering discrete and continuous structures. Dependency of error with damping level and time step size are evaluated and compared to other implicit and explicit methods. In the usual ranges of damping present in real structures, the proposed approach shows highly accurate results with respect to the other methods. Moreover, the computational time is very competitive with respect to other explicit methods, such as Modified PIM.


\section{Homotopy analysis based on artificial perturbation of damping}

Most problems of structural dynamics are governed by weak dissipation mechanisms so that the viscous damping forces are generally small compared to the inertial and elastic forces. Mathematically, the terms $\mathbf{C} \dot{\mbu}$ of Eqs.~\eqref{eq001}  can be considered as a perturbation of the undamped case, something that can be reproduced by an artificial parameter. We propose to modify Eq.~\eqref{eq001} with the parameter $\epsilon$, with $0 \leq \epsilon \leq 1$, which multiplies the damping matrix in Eq. \eqref{eq001} resulting 
\begin{equation}
\mbf{M} \ddot{\bm{\mathfrak{u}}} + \ep \, \mbf{C} \dot{\bm{\mathfrak{u}}} + \mbf{K} \bm{\mathfrak{u}} = \mathbf{f}(t)
\label{eq002}
\end{equation}
Now, the solution of Eq.~\eqref{eq002}  is a two-variable vector $\mathfrak{u}(t,\epsilon)$. The $\epsilon$--dependent solution can be expanded as
\begin{equation}
\bm{\mathfrak{u}}(t,\epsilon) = \sum_{n=0}^\infty \mbf{x}^{(n)}(t) \, \epsilon^n 
= \mbf{x}^{(0)}(t) + \mbf{x}^{(1)}(t) \, \epsilon 
+ \mbf{x}^{(2)}(t) \, \epsilon^2 + \cdots
\label{eq003}
\end{equation}
where $\{\mbf{x}^{(n)}(t)\}_{n=0}^\infty \in \mathbb{R}^N$ is a sequence of functions to be determined.  If series \eqref{eq003} is convergent for $0 \leq \epsilon \leq 1$ then the solution of our problem can be written as
\begin{equation}
\mbf{u}(t) = \bm{\mathfrak{u}}(t,1) = \sum_{n=0}^\infty \mbf{x}^{(n)}(t)  = \mbf{x}^{(0)}(t) + \mbf{x}^{(1)}(t) + \mbf{x}^{(2)}(t) + \cdots
\label{eq004}
\end{equation}
At the lower limit of the range, $\epsilon=0$, it is clear that we are in front of the undamped problem
\begin{equation}
\bm{\mathfrak{u}}(t,0) = \mbf{x}^{(0)}(t)
\label{eq005}
\end{equation}
The initial conditions are satisfied by the first iteration $\mbf{x}^{(0)}(t)$, so that $\mbf{x}^{(0)}(0)= \mbf{u}_0$ and $\dot{\mbf{x}}^{(0)}(0) = \dot{\mbf{u}}_0$. Thus, the problem for $n=0$ can be written as
\begin{equation}
\begin{cases}
\mbf{M} \ddot{\mbf{x}}^{(0)}+  \mbf{K} \, \mbf{x}^{(0)} = \mbf{f}(t)  & \\
\mbf{x}^{(0)}(0) = \mbf{u}_0 \ , \quad \dot{\mbf{x}}^{(0)}(0) = \dot{\mbf{u}}_0 &
\end{cases}
\label{eq006}
\end{equation}
Substituting Eq.~\eqref{eq003} into Eq.~\eqref{eq002} and identifying the coefficients of $\epsilon^n$, the $n$th function can be determined. After some straight operations it yields
\begin{equation}
\begin{cases}
\mbf{M} \ddot{\mbf{x}}^{(n)}+  \mbf{K} \, \mbf{x}^{(n)} = -  \mbf{C} \, \dot{\mbf{x}}^{(n-1)}& \\
\mbf{x}^{(n)}(0) = \mbf{0} \ , \quad \dot{\mbf{x}}^{(n)}(0) = \mbf{0} &, \  n \geq 1
\end{cases}
\label{eq007}
\end{equation}
where homogeneous initial conditions have been imposed for $n \geq 1$ because the real ones have already been applied for $n=0$ in Eq.~\eqref{eq006}. Let us see in detail the derivation of the resulting double iterative scheme, i.e. in the time domain and in the domain of the asymptotic expansion terms. 

\subsection{Solution of the first iteration, $n=0$}

The time domain will be sampled by the time-step $\Delta t$, i.e. $\{t_k, \  k\geq 0\} $, with $t_0 = 0$ and $t_{k+1} = t_k + \Delta t$. The first iteration ($n=0$) at $t=t_k$ will be denoted by $\mbf{x}_k^{(0)} = \mbf{x}^{(0)}(t_k)$. A recursive scheme is built using the explicit method based on Green's functions~\cite{Mansur-2007} yielding
\begin{eqnarray}
\mbf{x}_{k+1}^{(0)} &=&\mbf{G}(\dt) \, \mbf{x}_{k}^{(0)} + \mbf{H}(\dt) \, \dot{\mbf{x}}_{k}^{(0)} 
+ \int_{t=t_k}^{t_{k+1}} \mbf{H}(t_{k+1}-t) \mbf{M}^{-1} \mbf{f}(t) \, dt  \nonumber \\
\dot{\mbf{x}}_{k+1}^{(0)} &=&\dot{\mbf{G}}(\dt) \, \mbf{x}_{k}^{(0)} + \dot{\mbf{H}}(\dt) \, \dot{\mbf{x}}_{k}^{(0)} 
+ \int_{t=t_k}^{t_{k+1}} \mbf{G}(t_{k+1}-t) \mbf{M}^{-1} \mbf{f}(t) \, dt  
\label{eq008}										
\end{eqnarray}
The Green's functions, $\mbf{G}(t)$ and $\mbf{H}(t)$,  are the solutions of the following system of matrix differential equations
\begin{equation}
\begin{cases}
\dot{\mbf{G}}(t) = - \mbf{A} \, \mbf{H}(t) &  , \ \mbf{H}(0) = \mbf{0} \nonumber \\
\dot{\mbf{H}}(t) = \mbf{G}(t) & , \ \mbf{G}(0) = \mbf{I}_N \nonumber 
\end{cases}
\label{eq010}
\end{equation}
where $\mbf{A} = \mbf{M}^{-1} \mbf{K}$. It is straightforward that both $\mbf{G}(t)$ and $\mathbf{H}(t)$ are solutions of the second order matrix differential equations $\mbf{M}\ddot{Z} + \mbf{K} \mbf{Z} = \mbf{0}$, where $\mbf{Z}(t) \in \mathbb{R}^{N\times N}$ and hence both functions can be written as
\begin{equation}
\mbf{G}(t) = \cos( \sqrt{\mbf{A}} \, t)  = \sum_{j=0}^{\infty} \frac{(-1)^{j} \, \mbf{A}^{j}}{(2j)!} t^{2j} \ , \qquad 
\mbf{H}(t) = \mbf{A}^{-1/2} \, \sin( \sqrt{\mbf{A}} \, t)   = \sum_{j=0}^{\infty} \frac{(-1)^{j} \,\mbf{A}^{j}}{(2j+1)!} t^{2j+1}   
\label{eq009}
\end{equation}
Sometimes $\mathbf{G}(t)$ and $\mathbf{H}(t)$ are named step-response and impulsive-response functions~\cite{Fung-1997}, respectively. In order to obtain a closed form of the recursive scheme, the integrals of the nonhomogeneous part must be evaluated in the interval $\cor{t_k,t_{k+1}}$. For that, we will assume a cubic interpolation of the external applied force using third order Lagrange polynomials, namely
\begin{equation}
\mbf{f}(t) \approx \sum_{i=1}^4 \mathcal{L}_i \p{\frac{t-t_k}{\dt}} \, \mbf{f} \p{t_k + \frac{(i-1)\dt}{3}} \quad , \quad t_k \leq t \leq t_{k+1} 
\label{eq009b}
\end{equation}
where the interpolation polynomials in terms of the parameter $\xi = (t - t_k) / \dt \in \left[0,1\right]$ are
%
\begin{equation}
\mathcal{L}_1(\xi) = \frac{1}{2}(1 - \xi)(2 - 3\xi) (1-3\xi)  \ , \quad
\mathcal{L}_2(\xi) = \frac{9}{2}(1 - \xi)(2 - 3\xi) \xi  		 \ , \quad
\mathcal{L}_3(\xi) = -\frac{9}{2}(1 - \xi)(1 - 3\xi) \xi   \ , \quad
\mathcal{L}_4(\xi) = \frac{1}{2}(1 - 3\xi) (2-3\xi) \xi 		\label{eq009c}
\end{equation}
Plugging Eq.~\eqref{eq009b} into the integrals of Eq.~\eqref{eq008} and after some operations we obtain
\begin{eqnarray}
\int_{t=t_k}^{t_{k+1}} \mbf{H}(t_{k+1}-t) \mbf{M}^{-1} \mbf{f}(t) \, dt  & \approx &
\sum_{i=1}^4 \mbf{L}_{ui}  \, \mbf{M}^{-1} \ \mbf{f}\p{t_k + \frac{(i-1)\dt}{3}} \nonumber \\
\int_{t=t_k}^{t_{k+1}} \mbf{G}(t_{k+1}-t) \mbf{M}^{-1} \mbf{f}(t) \, dt    &\approx &
\sum_{i=1}^4 \mbf{L}_{vi}  \, \mbf{M}^{-1} \ \mbf{f}\p{t_k + \frac{(i-1)\dt}{3}}
\label{eq011}
\end{eqnarray}
where
\begin{equation}
\mbf{L}_{ui} = \int_{t=t_k}^{t_{k+1}} \mbf{H}(t_{k+1}-t) \, \mathcal{L}_i \p{\frac{t-t_k}{\dt}} \, dt  \  , \quad 
\mbf{L}_{vi} = \int_{t=t_k}^{t_{k+1}} \mbf{G}(t_{k+1}-t) \, \mathcal{L}_i\p{\frac{t-t_k}{\dt}} \, dt \ , \quad 1 \leq i \leq 4
\label{eq012}
\end{equation}
Rearranging together displacements and velocities of the degrees of freedom in the same $2N$-size vector, the iterative process can be summarized as the following relationship
\begin{equation}
\mathbf{X}_{k+1}^{(0)} = \mbf{T} \, \mbf{X}_k^{(0)}  + \mbf{L} \, \mbf{g}_k \quad , \quad k \geq 0
\label{eq116}
\end{equation}
where 
\begin{equation}
\mbf{X}_k^{(0)} = 
\left\lbrace 
\begin{array}{c}
\mbf{x}_{k}^{(0)} \\
\dot{\mbf{x}}_{k}^{(0)} 
\end{array}
\right\rbrace  \in \mathbb{R}^{2N} 
\ , \quad
\mbf{X}_0^{(0)} = 
\left\lbrace 
\begin{array}{c}
\mbf{u}_0 \\
\dot{\mbf{u}}_{0}
\end{array}
\right\rbrace 
\ , \quad
\mbf{T} =
\left[ 
\begin{array}{cc}
\mbf{G}(\dt)  & \mbf{H}(\dt) \\
\dot{\mbf{G}}(\dt)  &\dot{\mbf{H}}(\dt) 
\end{array}
\right]  \in \mathbb{R}^{2N \times 2N} 
\end{equation}
and 

\begin{equation}
\mbf{L} = 
\left[ 
\begin{array}{cccccc}
\mbf{L}_{u1} & \mbf{L}_{u2} & \mbf{L}_{u3} & \mbf{L}_{u4}  \\
\mbf{L}_{v1} & \mbf{L}_{v2} & \mbf{L}_{v3} & \mbf{L}_{v4}  
\end{array}
\right]  \in \mathbb{R}^{2N \times 4N} 
\ , \quad
\mbf{g}_k =
\left\lbrace 
\begin{array}{l}
\mbf{M}^{-1} \mbf{f}\p{t_k}  \\
\mbf{M}^{-1} \mbf{f}\p{t_k+\Delta t/3} \\
\mbf{M}^{-1} \mbf{f}\p{t_k+2\Delta t/3} \\
\mbf{M}^{-1}  \mbf{f}\p{t_{k+1}} 
\end{array}
\right\rbrace  \in \mathbb{R}^{4N} 
\label{eq083}
\end{equation}
Integrals of Eq. \eqref{eq012} can analytically be determined  resulting in  power series expressions of matrix $\mbf{A}= \mbf{M}^{-1} \mbf{K}$. After some straight  operations, the matrices  $\mbf{T}$ and $\mbf{L}$ can be written in compact form as
\begin{equation}
\mbf{T} =\sum_{j=0}^{\infty} \mbf{t}_j(\dt) \otimes \mbf{A}^j  \quad , \quad
\mbf{L} = \sum_{j=0}^{\infty} \bm{l}_j(\Delta t) \otimes \mbf{A}^j 
\label{eqMatrixL}
\end{equation}
where $ \otimes $ denotes the Kronecker product and the sequence of $\dt$--dependent matrices $\{\mathbf{t}_j(\dt) \}_{j=0}^\infty \in \mathbb{R}^{2 \times 2}$ and  $\{\bm{l}_j(\Delta t)\}_{j=0}^\infty \in \mathbb{R}^{2 \times 4} $,  are

\begin{eqnarray}
\mbf{t}_j(\dt) &=& 
\frac{(-1)^j \, \Delta t^{2j-1}}{(2j)!}
\left[ 
\begin{array}{cccc}
\dt	&	\dt^2 / (2j + 1) \\
2j	& \dt 
\end{array}
\right] 
\quad , \quad j \geq 0 \label{eq014a} \\
\bm{l}_j(\Delta t) &=& 
\frac{(-1)^j \, \Delta t^{2j+1}}{(2j+4)!}
\displaystyle
\left[ 
\begin{array}{cccc}
\ds \frac{(j+1)(8j^2 + 18j + 13) \Delta t}{2j+5}	&	\ds	\frac{36(j+1)^2 \Delta t}{2j+5}	&\ds	-\frac{9(2j^2 + j-1) \Delta t}{2j+5}  &\ds \frac{2(1+2j^2 ) \Delta t}{2j+5} \\
\ds (2j+1)(4j^2 + 5j + 3)							    	&	9(2j+1)^2	&\ds	-9 (j-1)(2j+1)  & 4j^2 - 4j + 3 
\end{array}
\right]  \nonumber \\
\label{eq014b}
\end{eqnarray}
The scheme of Eq. \eqref{eq116} corresponds to the initial iteration (undamped) of the asymptotic expansion ($n=0$). Now it is the turn to obtain the recursive formula that allows to find the rest of the terms of the expansion ($n\geq1$) for all time samples ($k\geq 0$).

\subsection{Solution of $n$th iteration, $n \geq 1$}

According to Eq.~\eqref{eq007}, the $n$th iteration of the asymptotic expansion, $\mbf{x}^{(n)}(t)$, is in turn solution of the undamped problem with an applied force proportional to the previous iteration. Thus, we can use again the Green's functions-based approach, yielding
\begin{eqnarray}
\mbf{x}_{k+1}^{(n)} &=&\mbf{G}(\dt) \, \mbf{x}_{k}^{(n)} + \mbf{H}(\dt) \, \dot{\mbf{x}}_{k}^{(n)} 
- \int_{t=t_k}^{t_{k+1}} \mbf{H}(t_{k+1}-t) \mbf{M}^{-1} \mbf{C} \ \dot{\mbf{x}}^{(n-1)}(t)  \, dt  \nonumber \\
\dot{\mbf{x}}_{k+1}^{(n)} &=&\dot{\mbf{G}}(\dt) \, \mbf{x}_{k}^{(n)} + \dot{\mbf{H}}(\dt) \, \dot{\mbf{x}}_{k}^{(n)} 
- \int_{t=t_k}^{t_{k+1}} \mbf{G}(t_{k+1}-t) \mbf{M}^{-1} \mbf{C} \ \dot{\mbf{x}}^{(n-1)}(t) \, dt  
\label{eq015}										
\end{eqnarray}
where now $\mbf{x}_k^{(n)} = \mbf{x}^{(n)}(t_k)$ and $\dot{\mbf{x}}_k^{(n)} = \dot{\mbf{x}}^{(n)}(t_k)$, with initial conditions $\mbf{x}_0^{(n)}=\dot{\mbf{x}}_0^{(n)} = \mbf{0}$. Above, $\dot{\mbf{x}}^{(n-1)}(t)$ is not explicitly available, since the solution is being determined from the scheme at the samples $ t_k, \ k=0,1,2,\ldots$. Eqs.~\eqref{eq015} will be transformed into an iterative scheme in two stages. First, let us integrate by parts to express the integrals of Eq.~\eqref{eq015} in terms of the displacements field $\mbf{x}^{(n-1)}(t)$. Thus,  using the definition of $\mbf{G}(t)$ and $\mbf{H}(t)$ we have
\begin{eqnarray}
- \int_{t=t_k}^{t_{k+1}} \mbf{H}(t_{k+1}-t)  \, \mbf{M}^{-1} \mbf{C} \ \dot{\mbf{x}}^{(n-1)}(t)  \, dt &=&
\mbf{H}(\dt) \mbf{M}^{-1} \mbf{C} \, \mbf{x}_{k}^{(n-1)} \nonumber \\				
&-&  \int_{t=t_k}^{t_{k+1}} \mbf{G}[t_{k+1}-t]   \, \mbf{M}^{-1} \mbf{C} \ \mbf{x}^{(n-1)} (t) \, dt \label{eq016}			 \\
- \int_{t=t_k}^{t_{k+1}} \mbf{G}(t_{k+1}-t) \, \mbf{M}^{-1} \mbf{C} \ \dot{\mbf{x}}^{(n-1)}(t) \, dt  &=&
- \mbf{G}(0) \, \mbf{M}^{-1} \mbf{C} \, \mbf{x}_{k+1}^{(n-1)}  + 
\mbf{G}(\dt) \, \mbf{M}^{-1} \mbf{C} \, \mbf{x}_{k}^{(n-1)}  \nonumber \\
&+&  
\mbf{A}  \int_{t=t_k}^{t_{k+1}} \mbf{H}[t_{k+1}-t] \, \mbf{M}^{-1} \mbf{C} \  \mbf{x}^{(n-1)} (t)\, dt 			
\end{eqnarray}
The second stage consist of the approximation of $\mbf{x}^{(n-1)} (t)$ by cubic splines using both displacements and velocities at $t=t_k$ and $t=t_{k+1}$. Within the interval $\cor{t_k,t_{k+1}}$ the dimensionless variable $\xi$ is defined as
\begin{eqnarray}
\xi = \frac{t - t_k}{t_{k+1}- t_k} \quad , \quad 0 \leq \xi \leq 1 \ , \quad t_k \leq t \leq t_{k+1}
\label{eq016b}
\end{eqnarray}
The value of $\mbf{x}^{(n-1)} (t)$ can then be approximated by
\begin{equation}
\mbf{x}^{(n-1)} \p{t} \approx \mathcal{N}_1(\xi) \, \mbf{x}_k^{(n-1)}		+ 		
\mathcal{D}_1(\xi) \, \dt \, \dot{\mbf{x}}_k^{(n-1)} + 
\mathcal{N}_2(\xi) \, \mbf{x}_{k+1}^{(n-1)}		+ 		
\mathcal{D}_2(\xi) \, \dt \, \dot{\mbf{x}}_{k+1}^{(n-1)}
\label{eq018}
\end{equation}
where 
\begin{equation}
\mathcal{N}_1(\xi) 	 =  (1 - \xi)^2 (1 + 2 \xi)   \ , \quad 
\mathcal{N}_2(\xi)  = (3-2 \xi ) \xi ^2 \ , \quad 
\mathcal{D}_1(\xi)  = (1 - \xi)^2 \xi 		\ , \quad 
\mathcal{D}_2(\xi)  = - (1 - \xi) \xi^2 \label{eq017}
\end{equation}
By means of this interpolation, the function and its derivative coincide with those evaluated from the approximation at the two ends of the interval, at $t=\{t_k, t_{k+1}\}$.  After plugging Eq.~\eqref{eq018} into Eq~\eqref{eq016}, the integrals can explicitly be written  in terms of the $(n-1)$th iteration at $t_k$ and $t_{k+1}$, yielding
\begin{eqnarray}
- \int_{t=t_k}^{t_{k+1}} \mbf{H}(t_{k+1}-t)  \, \mbf{M}^{-1} \mbf{C} \ \dot{\mbf{x}}^{(n-1)}(t)  \, dt &=&
\bm{\alpha}_{uu} \, \mbf{x}_k^{(n-1)}	+
\bm{\alpha}_{uv} \, \dot{\mbf{x}}_k^{(n-1)} + 
\bm{\beta}_{uu} \, \mbf{x}_{k+1}^{(n-1)}	 +
\bm{\beta}_{uv} \, \dot{\mbf{x}}_{k+1}^{(n-1)} \nonumber \\
- \int_{t=t_k}^{t_{k+1}} \mbf{G}(t_{k+1}-t) \, \mbf{M}^{-1} \mbf{C} \ \dot{\mbf{x}}^{(n-1)}(t) \, dt  &=&
\bm{\alpha}_{vu} \, \mbf{x}_k^{(n-1)}	+
\bm{\alpha}_{vv} \, \dot{\mbf{x}}_k^{(n-1)} + 
\bm{\beta}_{vu} \, \mbf{x}_{k+1}^{(n-1)}	 +
\bm{\beta}_{vv} \, \dot{\mbf{x}}_{k+1}^{(n-1)} \label{eq019}
\end{eqnarray}
where the $N\times N$ matrices of coefficients are 
\begin{align}
\bm{\alpha}_{uu} &=  \mbf{H}(\dt) \mbf{M}^{-1} \mbf{C}  -
\dt  \int_{\xi=0}^{1} \mbf{G}[\dt(1-\xi)]   \, \mbf{M}^{-1} \mbf{C} \ \mathcal{N}_1(\xi) \, d\xi 
& \bm{\alpha}_{uv} &= - \dt^2  \int_{\xi=0}^{1} \mbf{G}[\dt(1-\xi)]   \, \mbf{M}^{-1} \mbf{C} \ \mathcal{D}_1(\xi) \, d\xi 		\nonumber \\
\bm{\alpha}_{vu} &=  \mbf{G}(\dt) \, \mbf{M}^{-1} \mbf{C}  +
\dt \, \mbf{A}  \int_{\xi=0}^{1} \mbf{H}[\dt(1-\xi)]   \, \mbf{M}^{-1} \mbf{C} \ \mathcal{N}_1(\xi) \, d\xi 
& \bm{\alpha}_{vv} &= \dt^2 \, \mbf{A}  \int_{\xi=0}^{1} \mbf{H}[\dt(1-\xi)]   \, \mbf{M}^{-1} \mbf{C} \,  \mathcal{D}_1(\xi) \, d\xi \nonumber \\					
\bm{\beta}_{uu} &=  - \dt  \int_{\xi=0}^{1} \mbf{G}[\dt(1-\xi)]   \, \mbf{M}^{-1} \mbf{C} \ \mathcal{N}_2(\xi) \, d\xi 
& \bm{\beta}_{uv} &= - \dt^2  \int_{\xi=0}^{1} \mbf{G}[\dt(1-\xi)]   \, \mbf{M}^{-1} \mbf{C} \  \mathcal{D}_2(\xi) \, d\xi 		\nonumber \\
\bm{\beta}_{vu} &=  - \mbf{M}^{-1} \mbf{C}  +
\dt \, \mbf{A}  \int_{\xi=0}^{1} \mbf{H}[\dt(1-\xi)]   \, \mbf{M}^{-1} \mbf{C} \ \mathcal{N}_2(\xi) \, d\xi 
&\bm{\beta}_{vv} &= \dt^2 \, \mbf{A}  \int_{\xi=0}^{1} \mbf{H}[\dt(1-\xi)]   \, \mbf{M}^{-1} \mbf{C} \,  \mathcal{D}_2(\xi) \, d\xi \label{eq020}				
\end{align}
After some rearrangements, Eqs.~\eqref{eq015} lead to the following double iterative scheme, both in $k$ and $n$
\begin{equation}
\mathbf{X}_{k+1}^{(n)} = \mbf{T} \, \mbf{X}_k^{(n)}  + \bm{\alpha} \, \mbf{X}_k^{(n-1)}  +\bm{\beta} \, \mbf{X}_{k+1}^{(n-1)} 
\ , \quad n \geq 1 \ , \quad k \geq 0
\label{eq021}
\end{equation}
where the $n$th state vector  and the initial conditions are 
\begin{equation}
\mbf{X}_k^{(n)} = 
\left\lbrace 
\begin{array}{c}
\mbf{x}_{k}^{(n)} \\
\dot{\mbf{x}}_{k}^{(n)} 
\end{array}
\right\rbrace 
\ , \quad
\mbf{X}_0^{(n)} = 
\left\lbrace 
\begin{array}{c}
\mbf{0} \\
\mbf{0}
\end{array}
\right\rbrace \ , \quad n\geq 1
\label{eq022}
\end{equation}
and the block matrices  $\bm{\alpha},\bm{\beta} \in \mathbb{R}^{2N \times 2N}$
\begin{equation}
\bm{\alpha} =
\left[ 
\begin{array}{cc}
\bm{\alpha}_{uu}  & \bm{\alpha}_{uv} \\
\bm{\alpha}_{vu}  & \bm{\alpha}_{vv}
\end{array}
\right] 
\ , \quad
\bm{\beta} =
\left[ 
\begin{array}{cc}
\bm{\beta}_{uu}  & \bm{\beta}_{uv} \\
\bm{\beta}_{vu}  & \bm{\beta}_{vv}
\end{array}
\right] 
\label{eq023}
\end{equation}
Eqs.~\eqref{eq020} can analytically be solved using the series expansion of $\mbf{G}(t)$ and $\mbf{H}(t)$ given by Eq.~\eqref{eq009}. After some algebra, the matrices $\bm{\alpha} $ and $\bm{\beta} $ can be written as 
\begin{eqnarray}
\bm{\alpha} &=& \sum_{j=0}^\infty \bm{\alpha}_j(\dt) \otimes \p{\mbf{A}^{j} \, \mbf{M}^{-1} \mbf{C}  } 
= \cor{\sum_{j=0}^\infty \bm{\alpha}_j(\dt) \otimes \mbf{A}^{j} } \cdot  \p{\mbf{I}_2 \otimes \, \mbf{M}^{-1} \mbf{C}  }  \nonumber \\
\bm{\beta} &=& \sum_{j=0}^\infty \bm{\beta}_j(\dt) \otimes \p{\mbf{A}^{j} \, \mbf{M}^{-1} \mbf{C}  } 
= \cor{\sum_{j=0}^\infty \bm{\beta}_j(\dt) \otimes \mbf{A}^{j} } \cdot  \p{\mbf{I}_2 \otimes \, \mbf{M}^{-1} \mbf{C}  }  \label{eq037}
\end{eqnarray}
where the sequences $\{\bm{\alpha}_{j} (\dt)\}$ and $\{\bm{\beta}_{j} (\dt)\}$, $j \geq 0$, of $2\times2$ matrices are defined as
\begin{eqnarray}
\bm{\alpha}_{j} (\dt) &=& 
\frac{(-1)^j \dt^{2j}}{(2j+4)!}
\left[ 
\begin{array}{rr}
12 (j+1) \dt 		&		- 2 (2j + 1) (j+1) \dt^2		\\
12 (2j + 1)(j+2) &	 -4 j (2 j+1) (j+2) \dt
\end{array}
\right]  \ , \ j = 0, 1, 2, \ldots
\label{eq038} \\
\bm{\beta}_{j} (\dt) &=&
\frac{(-1)^j\dt^{2j}}{(2j+4)!}
\left[ 
\begin{array}{rr}
-12 (j+1) \dt 		&		2 (2j + 1) \dt^2		\\
-12 (2j + 1)(j+2) &		 8j(j+2) \dt
\end{array}
\right]  \ , \ j = 0, 1, 2, \ldots
\label{eq039}
\end{eqnarray}
Eqs.~\eqref{eq116} and \eqref{eq021} provide a well-defined iterative framework that enables the recursive computation in the $k$- and $n$-domains  of each one of the asymptotic expansion terms, i.e. $\mathbf{X}_k^{(n)}$, starting from $\mathbf{X}_0^{(0)}= \mathbf{U}_0$ and $\mathbf{X}_0^{(n)}= \mathbf{0}, \ n \geq 1$. In the next section the series which gives rise to the asymptotic solution will be summed.

\section{Proposed step-by-step explicit scheme}
\label{NumericalScheme}
Reaching $n$ iterations, then the transient solution at $t = t_k$ can be approximated by the sequence of partial sums, denoted by $\mathbf{U}_k^{(n)}$.  Namely,
\begin{equation}
\mbf{U}_k^{(n)} = 
\left\lbrace 
\begin{array}{c}
\mbf{u}^{(n)}(t_{k}) \\
\dot{\mbf{u}}^{(n)}(t_{k}) 
\end{array}
\right\rbrace 
\equiv
\left\lbrace 
\begin{array}{c}
\mbf{u}_{k}^{(n)} \\
\dot{\mbf{u}}_{k}^{(n)} 
\end{array}
\right\rbrace 
=  \mbf{X}_k^{(0)} + \mbf{X}_k^{(1)}  + \cdots + \mbf{X}_k^{(n)} = \sum_{\nu=0}^n \mbf{X}_k^{(\nu)} 
\label{eq028}
\end{equation}
The limit of the above sequence at $t=t_k$, provided that it exists, represents the proposed approximated solution to the transient problem. It will be denoted by 
\begin{equation}
\mbf{U}_k = \left\lbrace 
\begin{array}{c}
\mbf{u}_k \\
\dot{\mbf{u}}_k 
\end{array}
\right\rbrace 
\equiv \lim_{n\to \infty} \mbf{U}_k^{(n)} = \sum_{n=0}^\infty \mbf{X}_k^{(n)}  \ , \quad k \geq 0
\label{eq028b}
\end{equation}
It turns out that the above limit can be evaluated leading to an explicit time integration algorithm. Indeed, writing the $n+1$ first iterations, $0 \leq \nu \leq n$
\begin{align}
\mathbf{X}_{k+1}^{(0)} &= \mbf{T} \, \mbf{X}_k^{(0)}  + \mbf{L} \, \mbf{g}_k 	& \nu=0  \nonumber \\
- \bm{\beta} \, \mbf{X}_{k+1}^{(0)} + \mathbf{X}_{k+1}^{(1)} &= \mbf{T} \, \mbf{X}_k^{(1)}  + \bm{\alpha} \, \mbf{X}_k^{(0)} 
& \nu = 1\nonumber \\
\vdots &  \qquad \vdots				&   \nonumber \\
- \bm{\beta} \, \mbf{X}_{k+1}^{(n-1)} + \mathbf{X}_{k+1}^{(n)} &= \mbf{T} \, \mbf{X}_k^{(n)}  + \bm{\alpha} \, \mbf{X}_k^{(n-1)} 
& \nu = n\label{eq029}
\end{align}
Summing the left and the right side terms of the above equations and taking into account that $\mbf{U}_k^{(n)} = \sum_{\nu=0}^n \mbf{X}_k^{(\nu)} $, yields
\begin{equation}
- \bm{\beta} \, \mbf{U}_{k+1}^{(n-1)}  + \mathbf{U}_{k+1}^{(n)} = \mbf{T} \, \mbf{U}_k^{(n)}  + \bm{\alpha} \, \mbf{U}_k^{(n-1)}  + \mbf{L} \mbf{g}_k
\ , \quad n \geq 1
\label{eq030}
\end{equation}
Taking limits now on both sides as $n \to \infty $ and assuming that $\mbf{U}_k = \lim_{n\to \infty} \mbf{U}_k^{(n)} $ exists, the following implicit algorithm arises
\begin{equation}
\p{\mbf{I}_{2N} - \bm{\beta} } \, \mbf{U}_{k+1}   = \p{\mbf{T} \,  + \bm{\alpha} }\, \mbf{U}_k+ \mbf{L} \mbf{g}_k
\label{eq031}
\end{equation}
which, however, can be written in explicit form just multiplying by the corresponding inverse matrix, namely
\begin{equation}
\mbf{U}_{k+1}   = \mbf{a}  \, \mbf{U}_k+ \mbf{b}_k \ , \quad k = 0,1,2,\ldots
\label{eq032}
\end{equation}
where the main matrices of the algorithm, $\mbf{a} \in \mathbb{R}^{2N \times 2N}$ and $\mbf{b}_k \in \mathbb{R}^{2N}$, are 
\begin{eqnarray}
\mbf{a} &=& \p{\mbf{I}_{2N} - \bm{\beta} }^{-1} \p{\mbf{T} \,  + \bm{\alpha} } \ \in \mathbb{R}^{2N \times 2N}\label{eq033a} \\
\mbf{b}_k &=& \p{\mbf{I}_{2N} - \bm{\beta} }^{-1} \mbf{L} \, \mbf{g}_k   \ \in \mathbb{R}^{2N}\quad , \quad k = 0,1,2,\ldots \label{eq033b}
\end{eqnarray}
The calculation of the inverse matrix does not add any drawback with respect to the loss of accuracy because, as we will see later, the corresponding Neumann series can be invoked. Essentially, the iterative procedure shown in Eq. \eqref{eq032} summarizes the main contribution of this article. In the following sections, the properties of the numerical approach will be discussed, indeed
\begin{description}
	\item [Sec. \ref{Convergence}:]  Conditions for convergence of asymptotic expansion $\mathbf{U}_k = \sum_{n=0}^\infty \mbf{X}_k^{(n)} $
	\item [Sec. \ref{Computation}:]  Algorithm for the computation of main matrices, $\mathbf{a}$ and $\mathbf{b}_k$.
	\item [Sec. \ref{Stability}:] Analysis of stability 
	\item [Sec. \ref{Accuracy}:] Order of accuracy
	\item [Sec. \ref{ComputationalEffort}:] Computation effort
\end{description}

\section{Analysis of convergence}
\label{Convergence}

The iterative scheme outlined in Eq.~\eqref{eq032}  is the result of assuming as true that the limit of Eq.~\eqref{eq028b} exists, i.e.  the series of the asymptotic method is convergent. In this section,  the conditions under which such limit actually exists will be addressed. For this purpose, we will deduce a recursive relation in the space of the asymptotic expansion ($n$-space) that allows us to explicitly derive the sequence of partial sums. Evaluations of Eq.~\eqref{eq021} for $n \geq 2$ and for $k=\{0,1,\ldots,k-1\}$ yields
\begin{eqnarray}
\mathbf{X}_{1}^{(n)} &=& \mbf{T} \, \mbf{X}_0^{(n)}  + \bm{\alpha} \, \mbf{X}_0^{(n-1)}  +\bm{\beta} \, \mbf{X}_{1}^{(n-1)}  \nonumber \\
\mathbf{X}_{2}^{(n)} &=& \mbf{T} \, \mbf{X}_1^{(n)}  + \bm{\alpha} \, \mbf{X}_1^{(n-1)}  +\bm{\beta} \, \mbf{X}_{2}^{(n-1)}  \nonumber \\
\vdots					 &=&  \vdots \nonumber \\
\mathbf{X}_{k}^{(n)} &=& \mbf{T} \, \mbf{X}_{k-1}^{(n)}  + \bm{\alpha} \, \mbf{X}_{k-1}^{(n-1)}  +\bm{\beta} \, \mbf{X}_{k}^{(n-1)} 
\  , \quad  n \geq 2 \label{eq056}
\end{eqnarray}
Since  $\mbf{X}_0^{(n)} = \mbf{X}_0^{(n-1)}= \mbf{0}$ for $n \geq 2$, then  the first step $(k=0)$ is reduced to 
$\mathbf{X}_{1}^{(n)} = \bm{\beta} \, \mbf{X}_{1}^{(n-1)} $. The terms $\mbf{X}_1^{(n)},  \ldots ,\mbf{X}_{k-1}^{(n)}$ of the right side of each step equation can be replaced by their corresponding expression of the previous step. Repeating this process systematically for each step and after some operations, Eqs.~\eqref{eq056} can be written separating the terms associated to the $n$th and $(n-1)$th iterations in both sides of equations, i.e.
\begin{eqnarray}
\mathbf{X}_{1}^{(n)} &=& \bm{\beta} \, \mbf{X}_{1}^{(n-1)}  \nonumber \\
\mathbf{X}_{2}^{(n)} &=& \bm{\gamma} \mbf{X}_1^{(n-1)}  +\bm{\beta} \, \mbf{X}_{2}^{(n-1)}  \nonumber \\
\mathbf{X}_{3}^{(n)} &=& \mbf{T}\bm{\gamma} \mbf{X}_1^{(n-1)}  +\bm{\gamma} \, \mbf{X}_{2}^{(n-1)}  +  \bm{\beta} \, \mbf{X}_{3}^{(n-1)}  \nonumber \\
\vdots					 &=&  \vdots \nonumber \\
\mathbf{X}_{k}^{(n)} &=& \sum_{r=1}^{k-1} \mbf{T}^{k-1-r} \, \bm{\gamma} \, \mbf{X}_{r}^{(n-1)}  +\bm{\beta} \, \mbf{X}_{k}^{(n-1)} 
\  , \quad  n \geq 2 \label{eq056b}
\end{eqnarray}
where $\bm{\gamma} = \mbf{T}  \bm{\beta} +  \bm{\alpha}$. Rearranging $\mbf{X}_1^{(n)}, \ldots, \mbf{X}_k^{(n)}$  in one single $2N k$--size column vector denoted by $\mbf{X}^{(n)} \in \mathbb{R}^{2Nk}$, and after some algebra, the above relations can be expressed as  
\begin{equation}
\mbf{X}^{(n)} = \mbf{R} \, \mbf{X}^{(n-1)} \ , \qquad n \geq 2
\label{eq025}
\end{equation}
where the $n$th iteration array $\mbf{X}^{(n)} $ and matrix $\mbf{R}$ are formed by the blocks 
\begin{equation}
\mbf{X}^{(n)} =
\cccolumna{\mbf{X}_1^{(n)}}{\vdots}{\mbf{X}_k^{(n)}}  \in \mathbb{R}^{2Nk}
\ , \quad
\mbf{R} = 
\left[ 
\begin{array}{cccccc}
\bm{\beta}  &  \mbf{0}_{2N} &  \mbf{0}_{2N}  & \cdots & \mbf{0}_{2N} & \mbf{0}_{2N} \\
\bm{\gamma}  & \bm{\beta} &  \mbf{0}_{2N}  & \cdots & \mbf{0}_{2N} & \mbf{0}_{2N} \\
\mbf{T} \bm{\gamma}  & \bm{\gamma} &\bm{\beta} & \cdots & \mbf{0}_{2N} & \mbf{0}_{2N} \\
\vdots			& \vdots           & \vdots	& \ddots 	& \vdots 		  &  \vdots          \\
\mbf{T}^{k-3} \bm{\gamma}   &  \mbf{T}^{k-4} \bm{\gamma} & \mbf{T}^{k-5} \bm{\gamma}  & \cdots & \bm{\beta}    &  \mbf{0}_{2N}   \\
\mbf{T}^{k-2} \bm{\gamma}   &  \mbf{T}^{k-3} \bm{\gamma} & \mbf{T}^{k-4} \bm{\gamma}  & \cdots & \bm{\gamma}    &  \bm{\beta}  \\
\end{array}
\right]  \in \mathbb{R}^{2Nk \times 2Nk}
\end{equation}
%

%
%
The first iteration for $n=1$ leads to the relationship between $\mbf{X}^{(0)}$ and $\mbf{X}^{(1)}$. It can be derived straightforward when Eq.~\eqref{eq056} is evaluated and taking into account that $\mbf{X}_0^{(1)} = \mathbf{0}$ and $\mbf{X}_0^{(0)}= \mbf{U}_0$. As before, using some matrix algebra, the following expression can be obtained
\begin{equation}
\mbf{X}^{(1)} = \mbf{R} \, \mbf{X}^{(0)} + \mbf{S} \, \mbf{U}_0 
\label{eq025b}
\end{equation}
where 
\begin{equation}
\mathbf{S} = 
\left[ 
\begin{array}{c}
\bm{\alpha} \\
\mathbf{T} \, \bm{\alpha} \\
\vdots \\
\mathbf{T}^{k-1} \, \bm{\alpha}  
\end{array}
\right]  \in \mathbb{R}^{2Nk \times 2N}
\label{eq025c}
\end{equation}
We can build the $2Nk$--size column vector $\mbf{U}$ with the response of all time steps, resulting the series expansion
\begin{equation}
\mbf{U} = \cccolumna{\mbf{U}_1}{\vdots}{\mbf{U}_k}  
= \cccolumna{\sum_{n=0}^\infty \mbf{X}_1^{(n)}  }{\vdots}{\sum_{n=0}^\infty \mbf{X}_k^{(n)}  } 
= \sum_{n=0}^\infty \cccolumna{\mbf{X}_1^{(n)}}{\vdots}{\mbf{X}_k^{(n)}}  = \sum_{n=0}^\infty \mbf{X}^{(n)}
\label{eq058}
\end{equation} 
Using the recursive relations of Eq.~\eqref{eq025} for $n\geq 2$ together with the first iteration  $(n=1)$ from Eq.~\eqref{eq025b}, the response $\mbf{U}$ results
\begin{equation}
\mbf{U} = \mbf{X}^{(0)} + \mbf{X}^{(1)} + \mbf{X}^{(2)} + \cdots + \mbf{X}^{(n)} + \cdots 
= \mbf{X}^{(0)} + \p{ \sum_{n=0}^\infty \mbf{R}^{n} }  \mbf{X}^{(1)}  
=  \p{ \sum_{n=0}^\infty \mbf{R}^{n} }  \p{\mbf{X}^{(0)} +  \mathbf{S} \, \mbf{U}_0}
\label{eq059}
\end{equation}
This closed-form relationship highlights the intrinsic nature of the proposed approach as an asymptotic expansion-based method. Note that both matrices $\mathbf{R}$ and $\mathbf{S}$ are directly proportional to the damping matrix, so that for undamped problems we have $\mathbf{U} \equiv \mathbf{X}^{(0)}$. The terms arranged within $\mathbf{X}^{(n)}$ represent the $n$th order contribution to the response and their effect is expected to decrease significantly for low damping problems. The proposed method leads to a solution as long as the Neumann series converges. The necessary and sufficient condition for this series to be convergent is that all the eigenvalues of matrix $\mbf{R}$ are strictly smaller than the unit in absolute value~\cite{Wilkinson-1988}. Therefore, the series is convergent provided that $\rho(\mbf{R}) <1$, where $\rho(\bullet)$ denotes the spectral radius of a matrix. Since $\mbf{R}$ is a triangular matrix by blocks and $\bm{\beta}$ is the matrix repeated in the block diagonal, then eigenvalues of $\mbf{R}$ are those of $\bm{\beta}$, and the convergence of the method holds provided that 
\begin{equation}
\rho(\bm{\beta}) < 1
\label{eq057}
\end{equation}
The  matrix $\bm{\beta}$ has already been presented in Eqs. \eqref{eq037} and \eqref{eq039} as a power series expansion of the system matrix $\mbf{A} = \mbf{M}^{-1} \mbf{K}$. Moreover, in view of the form of the matrix all terms are proportional to the damping $ \mbf{M}^{-1} \mbf{C}$. Therefore the spectral radius of the matrix will depend on the one hand on the order of truncation of the series and on the other hand it is expected to be increased proportionally to the intensity of the dissipative forces. Analytically, the expression of $\bm{\beta}$ given above can be expressed as a product of two matrices
\begin{equation}
\bm{\beta} = \dt \, \cor{\lim_{m \to \infty}\bm{\sigma}_m   \p{\dt \sqrt{\mbf{M}^{-1} \mbf{K}}} }  \cdot \p{\mbf{I}_2 \otimes \, \mbf{M}^{-1} \mbf{C}  } 
\label{eq060}
\end{equation}
where 
\begin{equation}
\bm{\sigma}_m(\bm{\tau}) = \sum_{j=0}^{m/2}
\frac{(-1)^j}{(2j+4)!}
\left[ 
\begin{array}{rr}
-12 (j+1) \, \bm{\tau}^{2j}   		&		2 (2j + 1) \dt	\, \bm{\tau}^{2j}	\\
-12 (2j + 1)(j+2) \, \bm{\tau}^{2j} /\dt &		 8j(j+2) \, \bm{\tau}^{2j} 
\end{array}
\right] 
\label{eq059b}
\end{equation}
The matrix $\bm{\sigma}_m(\bm{\tau})$ has $2 \times 2$ elements as long as the argument $\bm{\tau} \in \mathbb{R}$, but it will be $\bm{\sigma}(\bm{\tau}) \in \mathbb{R}^{2N \times 2N}$ if $\bm{\tau} \in \mathbb{R}^{N \times N}$, as shown in Eq.~\eqref{eq060} with $\bm{\tau} =\dt \sqrt{\mbf{M}^{-1} \mbf{K}}$ . The index $m$ denotes the maximum order of $\bm{\tau}$ in the power series expansion. In practice, a finite value for $m$ means the truncation of the series. \\

The nature of the matrix $\bm{\beta}$ is revealed in the decomposition shown in Eq.~\eqref{eq060}. Both time step $\dt$ and damping matrix $\mbf{C}$ play a relevant role in its spectral radius. Let us see that, under certain conditions, some practical outcomes can be drawn to estimate a priori an upper bound of the spectral radius $\rho(\bm{\beta})$ as a function of $\dt$ and of the system matrices. Indeed, it is known that the spectral radius of a product of matrices is less than or equal to the product of the respective spectral radii whenever the matrices involved commute. In the particular case of Eq.~\eqref{eq060}, the matrices $\bm{\sigma}_m\p{\dt \sqrt{\mbf{M}^{-1} \mbf{K}}} $ and $\p{\mbf{I}_2 \otimes \, \mbf{M}^{-1} \mbf{C}  } $  commute provided that $\mbf{M}^{-1} \mbf{K}$ and $ \mbf{M}^{-1} \mbf{C} $ commute, something that holds for classically (or proportionally) damped systems~\cite{Caughey-1965}. In such case, the spectral radius is bounded by
\begin{equation}
\rho(\bm{\beta}) \leq  \rho \cor{\bm{\sigma}_m   \p{\dt \sqrt{\mbf{M}^{-1} \mbf{K}}} } \, \dt \, \rho \p{ \mbf{M}^{-1} \mbf{C}  } 
\label{eq061}
\end{equation}
The four blocks of the matrix $\bm{\sigma}_m  \p{\dt \sqrt{\mbf{M}^{-1} \mbf{K}}}$ can be diagonalized in the modal space of the undamped problem, allowing the analytical determination of the eigenvalues. By modal analysis,  the undamped natural frequencies and eigenvectors $\om_j \in \mathbb{R}$ and $\bphi{j} \in \mathbb{R}^{N} $, are related throughout
\begin{equation}
\mbf{K} \bphi{j}  =  \om_j^2 \mbf{M} \bphi{j}   \ , \ 1 \leq j \leq N
\label{eq063}
\end{equation}
Assuming that there are not repeated frequencies, the eigenvectors can be mass-normalized using the biorthogonality relations $\bphi{j}^T \mbf{M} \bphi{k} = \delta_{jk}$ and  $\bphi{j}^T \mbf{K} \bphi{k} = \om_j^2 \, \delta_{jk} $
%
%
where $\delta_{jk} $ is the Kronecker delta. The eigenvectors $\bphi{j}$ can be arranged as the columns of the square matrix (modal matrix) $\mbPhi = \cor{\bphi{1},\ldots,\bphi{N}} \in \mbR^{N \times N}$ so that the biorthogonality relations can be written in matrix form as 
$\mbPhi^T \, \mbM \,  \mbPhi = \mbf{I}_N, \  \mbPhi^T \,  \mbK \,  \mbPhi = \mathbf{\Omega}^2 $
%
%
%
%
where $\mbf{I}_N \in \mathbb{R}^{N\times N} $ denotes the identity matrix and $ \mbf{\Om} = \textnormal{diag} \cor{\om_1,\ldots,\om_N}$ is a diagonal matrix with the natural frequencies. From these relations it is verified that 
\begin{equation}
\mbf{M}^{-1} \mbf{K} = \bm{\Phi} \, \bm{\Omega}^2 \, \bm{\Phi}^{-1} \quad , \quad 
\sqrt{\mbf{M}^{-1} \mbf{K} }= \bm{\Phi} \, \bm{\Omega} \, \bm{\Phi}^{-1}
\end{equation}
and therefore
\begin{equation}
\bm{\sigma}_m   \p{\dt \sqrt{\mbf{M}^{-1} \mbf{K}}} = 
\p{\mbf{I}_2 \otimes \, \bm{\Phi}  } \cdot 
\bm{\sigma}_m   \p{\dt \bm{\Omega} } \cdot 
\p{\mbf{I}_2 \otimes \, \bm{\Phi}^{-1}  } 
\label{eq067}
\end{equation}
From this matrix relation, both $\bm{\sigma}_m   \p{\dt \sqrt{\mbf{M}^{-1} \mbf{K}}} $ and $\bm{\sigma}_m   \p{\dt \bm{\Omega} } $ are similar matrices and consequently they share the same eigenvalues. In addition, the computation of eigenvalues from matrix $\bm{\sigma}_m   \p{\dt \bm{\Omega} } $ is straightforward because it is formed by four square diagonal blocks. The set of $2N$ eigenvalues yields 
\begin{equation}
\bigcup_{i=1}^N \{\mu_{1}(m,\om_i \dt), \mu_{2}(m,\om_i \dt)\}
\label{eq068}
\end{equation}
where $\{\mu_1(m,\tau),\mu_2(m,\tau)\}$ denotes the two eigenvalues of the $2 \times 2$ matrix $\bm{\sigma}_m(\tau)$ when $\tau \in \mathbb{R}$, which can analytically be computed. Thus, for instance for $m=0$ and $m=2$ the results are
\begin{align}
\bm{\sigma}_0(\tau) &= \mmatriz{-1/2}{\dt/12}{-1/\dt}{0}  & 
\mu_{1,2}(0,\tau) &=  \frac{1}{12} \left(-3 \pm i \sqrt{3}  \right) \nonumber \\
\bm{\sigma}_2(\tau) &= \mmatriz{-1/2 + \tau^2/30}{(1/12 - \tau^2/120)\dt}{(-1 + 3\tau^2/20)/\dt}{-\tau^2/30}  &
\mu_{1,2}(2,\tau) &=  \frac{1}{120} \left(-30 \pm i \sqrt{300 + 2\tau ^4 - 60 \tau ^2}\right) 
\label{eq070}
\end{align}
\begin{figure}[ht]%
	\begin{center}
		\includegraphics[width=10cm]{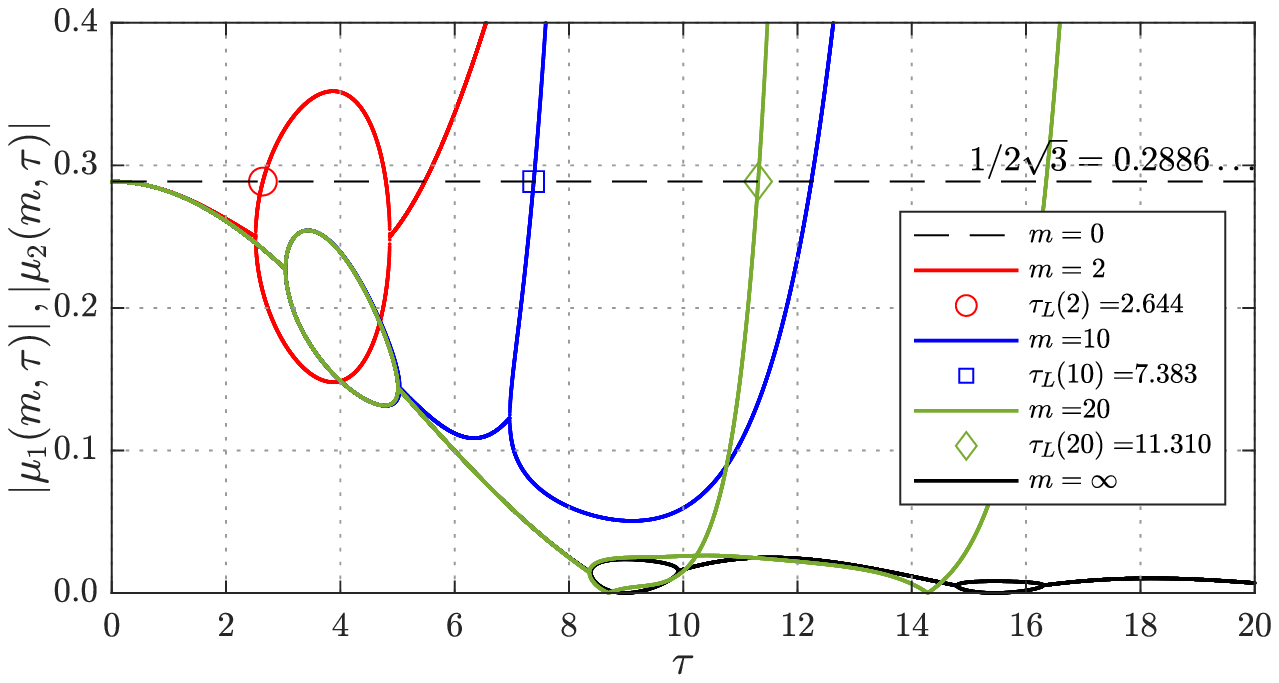} \\	
		\caption{Eigenvalues of $\bm{\sigma}_m(\tau)$ in absolute value: $\left|\mu_{1}(m,\tau) \right|$, $\left| \mu_{2}(m,\tau) \right|$   }%
		\label{fig01}%
	\end{center}
\end{figure}
As the order $m$ increases, the order of $\tau$ within $\mu_1(m,\tau)$ and $\mu_2(m,\tau)$ will increase as well. In general, the spectral radius of matrix $ \bm{\sigma}_m(\dt \sqrt{\mbf{M}^{-1} \mbf{K}}) $ can be written as
\begin{equation}
\rho \cor{\bm{\sigma}_m   \p{\dt \sqrt{\mbf{M}^{-1} \mbf{K}}}} = \max_{1 \leq i \leq N} \rho \cor{\bm{\sigma}_m(  \om_i \, \dt   ) } 
\label{eq062}
\end{equation}
Substituting the above expression into Eq.~\eqref{eq061}, then the general condition for the proposed numerical scheme in classically damped systems to be convergent is
\begin{equation}
\rho(\bm{\beta}) \leq \dt \, \rho \p{ \mbf{M}^{-1} \mbf{C}  } \,  \max_{1 \leq i \leq N} \rho \cor{\bm{\sigma}_m(  \om_i \, \dt   ) }   < 1
\label{eq073}
\end{equation}
The distribution of eigenvalues (in absolute value) of $\bm{\sigma}_m(\tau)  $ is shown in fig.~\ref{fig01} for different values of the truncation order $m$. For $m=0$, the eigenvalues of $\bm{\sigma}_0(\tau)$ are independent of $\tau$ and the spectral radius is $\rho \cor{\bm{\sigma}_0} = 1 / 2\sqrt{3}$. On the other side,  the eigenvalues of $\bm{\sigma}_\infty(\tau) $, as function of $\tau$, draw a decreasing curve as $\tau$ increases, as shown in fig.~\ref{fig01} (solid black line). Assuming infinite terms, i.e. $m = \infty$, it can be ensured that $\rho \cor{\bm{\sigma}_\infty(\tau) } \leq \rho (\bm{\sigma}_0) =  1 / 2\sqrt{3}$ and therefore $\rho(\bm{\beta}) < 1$ provided that
\begin{equation}
\dt \, \rho \p{ \mbf{M}^{-1} \mbf{C}  } < 2 \sqrt{3} \ , \quad (m=\infty)
\label{eq073b}
\end{equation}
In practice, the order of truncation $m$ must be finite. Hence, for each $m < \infty$, convergence does not hold within the whole range of $\om_i \dt$, existing a limit value of $\tau$ from which Eq.~\eqref{eq073} is unbounded. Let us denote by $\tau_L(m)$  to the solution of the equation $\rho \cor{\bm{\sigma}_m(\tau)} = 1 / 2 \sqrt{3}$. In fig.~\ref{fig01} the values of $\tau_L(m)$ have been highlighted for $m=2, \ 10$ and $ 20$ using a square, circle and diamond marker, respectively. Then, for each value of $m$, it is verified that 
\begin{equation}
\rho \cor{\bm{\sigma}_m   \p{\dt \sqrt{\mbf{M}^{-1} \mbf{K}}}} = \max_{1 \leq i \leq N} \rho \cor{\bm{\sigma}_m(  \om_i \, \dt   ) } 
\leq \frac{1}{2 \sqrt{3}} \quad , \quad 0 \leq \om_i \, \dt < \tau_L(m)
\label{eq062b}
\end{equation}
The inequalities $\om_i \dt < \tau_L(m), \ 1 \leq i \leq N$ hold provided that it holds for the last frequency, i.e.  $\om_{\max{}} \, \dt < \tau_L(m)$, where $\om_{\max{}} = \rho (\sqrt{\mbf{M}^{-1} \mbf{K}} ) = 2\pi / T$ stands for the maximum natural frequency in the system (rad/s) and $T$ for the minimum period (in seconds). Substituting Eq.~\eqref{eq062b} into Eq.~\eqref{eq073}, the condition $\rho(\bm{\beta})<1$ can be written as
\begin{equation}
\rho(\bm{\beta}) \leq \dt \, \rho \p{ \mbf{M}^{-1} \mbf{C}  } \,  \max_{1 \leq i \leq N} \rho \cor{\bm{\sigma}_m(  \om_i \, \dt   ) }  
\leq  \dt \, \rho \p{ \mbf{M}^{-1} \mbf{C}  }  \frac{1}{2 \sqrt{3}}  < 1 \ , \quad 0 \leq \om_{\max{}}  \, \dt < \tau_L(m)
\label{eq073c}
\end{equation}
Writing together both restrictions of Eq.~\eqref{eq073c} and fixed a truncation order $m$, it can be concluded that for classically damped systems, $\rho(\bm{\beta}) < 1$ provided that both $\dt \, \rho  (\mbf{M}^{-1} \mbf{C} ) < 2 \sqrt{3}$ and $\om_{\max{}} \dt < \tau_L(m)$  hold simultaneously. The following condition encapsulates both restrictions in a single expression as a limitation of the time step, 
\begin{equation}
\dt \leq \min \left\lbrace    \frac{2 \sqrt{3} }{\rho  (\mbf{M}^{-1} \mbf{C} ) }, 
\frac{\tau_L(m)}{\rho  (\sqrt{\mbf{M}^{-1} \mbf{K}}) } \right\rbrace =
\min \left\lbrace    \frac{2 \sqrt{3} }{\rho  (\mbf{M}^{-1} \mbf{C} ) }, 
\frac{\tau_L(m)}{2\pi  } T \right\rbrace 
\label{eq117}
\end{equation}
where the limit values $\tau_L(m)$ are shown in Table \ref{tab04} for some truncation orders. 
\begin{table}
	\begin{center}
		\begin{tabular}{lrr}
Truncation order $m$						&					$\tau_L(m)$					&					$ \tau_L(m)/2 \pi$				\\
\hline
$m = 2$				 								&					 2.64303						& 					0.42065																\\
$m = 4$				 								&					 9.59335						& 					1.52683								\\
$m = 6$												&					 5.48854							&					0.87353										\\
$m = 8$												&					 6.68027							&					1.06320										\\
$m=10$											  &  					7.38332								&					1.17509													\\
$m = 20$										 &						11.3105								&				  1.80012										\\
$m = 30$									     &					   15.1700								&					2.41438										\\
$m = 40$										&						19.0203								&					3.02718										\\
\hline
\end{tabular}
	\end{center}
	\caption{Values of $\tau_L(m)$ for different truncation orders, $m$. $\tau_L(m)$ is the solution of equation $\rho \cor{\bm{\sigma}_m(\tau)} = 1 / 2 \sqrt{3}$. Once the truncation order has been set, the parameter $\tau_L(m)$ is used to calculate the maximum time step admissible by the numerical method, for convergence reasons, being $\dt \leq \{2 \sqrt{3}/ \rho  (\mbf{M}^{-1} \mbf{C} ), \tau_L(m) T/2\pi \}$ }
	\label{tab04}
\end{table}
Although the condition derived in Eq. \eqref{eq117} is strictly valid for  classically damped systems, it can be considered acceptable for systems with light damping, which can be assumed to be slightly non-proportional, allowing a priori an estimation of the time step $\dt$. Furthermore, Eq. \eqref{eq117} provides information about the constraints on the application of the numerical model. Thus, in the presence of weak damping, the condition $\dt/T<\tau_L(m)/2\pi$ will be limiting. On the other hand, highly damped systems can also be analyzed with the proposed method, provided that the time step is sufficiently small. \\

\section{Algorithms for computation of main matrices}
\label{Computation}

To complete the proposed methodology, specific algorithms are needed to calculate the main matrices of the numerical approach, i.e. from Sec.~\ref{NumericalScheme}, Eq.~\eqref{eq032}, the matrix $\mbf{a}(\dt)$ and the vector $\mbf{b}_k(\dt)$. Throughout this section it is necessary to highlight the dependence of the parameter $\dt$ in both matrices. Rewriting their expressions
\begin{eqnarray}
\mbf{a}(\dt) &=& \p{\mbf{I}_{2N} - \bm{\beta} }^{-1} \p{\mbf{T} \,  + \bm{\alpha} } \label{eq081a} \\
\mbf{b}_k(\dt) &=& \p{\mbf{I}_{2N} - \bm{\beta} }^{-1} \mbf{L} \, \mbf{g}_k  \quad , \quad k = 0,1,2,\ldots \label{eq081b}
\end{eqnarray}
Both $\mbf{a}(\dt)$ and $\mbf{b}_k(\dt)$ are computed from matrices $\mbf{T}$, $\bm{\alpha}$, $\bm{\beta}$ and $\mbf{L}$, as well as the vector of external forces given by $\mbf{g}_k$. Apart from the latter, the rest of matrices are series expansions  in terms of the time step $\dt$ and of the system matrices, $\mbf{M}$, $\mbf{C}$ and $\mbf{K}$, in particular in terms of $\mbf{A} = \mbf{M}^{-1} \mbf{K}$ and $\mbf{M}^{-1} \mbf{C}$. Additionally,  the computation of the  inverse matrix $\p{\mbf{I}_{2N} - \bm{\beta} }^{-1} $ is also required. The condition for the convergence of the method $\rho(\bm{\beta})<1$, enables now to use the Neumann series for its computation by means of the formula
\begin{equation}
\p{\mbf{I}_{2N} - \bm{\beta} }^{-1}  = \sum_{r=0}^\infty  \bm{\beta}^r \approx \mbf{I}_{2N}  +  \bm{\beta} +  \bm{\beta}^2 + \cdots + \bm{\beta}^r 
\label{eq074}
\end{equation}
In the following points the algorithms for the efficient computation of both $\mathbf{a}(\dt)$ and $\mathbf{b}_k(\dt)$ will be described in detail. Although the same matrix $\bm{\beta}$ arises in both Eqs.~\eqref{eq081a} and \eqref{eq081b}, it  will be computed differently in both cases, using different time steps and truncation orders in the series involved.

\subsection{Precise computation of matrix $\mbf{a}(\dt)$}
\label{PreciseComp_a}

Accuracy in the evaluation of the main matrix of the algorithm $\mbf{a}(\dt)$ is a key issue to achieve satisfactory numerical results. In the absence of external forces it yields
\begin{equation}
\mbf{U}_{k+1}   = \mbf{a}(\dt)  \, \mbf{U}_k 
\label{eq034}
\end{equation}
Therefore, it is verified then that $\mbf{a}(\dt)  = \mbf{a}(\dt/2)   \mbf{a}(\dt/2) =  \cor{\mbf{a}(\dt/4)}^4 = \cdots $. Consider a positive integer $p$ and the reduced time step defined as $\dt_0 = \dt / 2^p$ then
\begin{equation}
\mbf{a}(\dt)  = \cor{ \mbf{a} \p{\frac{\dt}{2^p} } }^{2^p} = \cor{ \mbf{a}(\dt_0) }^{2^p}
\label{eq085}
\end{equation}
Moreover, $\mathbf{a}(\dt_0)$ can be expressed in incremental form as
\begin{equation}
\mbf{a}(\dt_0)  =   \mbf{I}_{2N}  + \delta \mbf{a}(\dt_0) 
\label{eq076}
\end{equation}
where $\delta \mbf{a}(\dt_0) $ is a very small matrix. Strictly it can be stated that $\lim_{\dt \to 0}  \delta \mbf{a}(\dt)  = \mbf{0} $. 
The relationship of Eq.~\eqref{eq085} allows the use of the so-called $2^p$ algorithm, commonly used for the evaluation of exponential matrices \cite{Moler-2003}. Moreover, Eq.~\eqref{eq076} enables to use this technique on the incremental part $\delta \mbf{a}(\dt)$ instead of on the total matrix, improving significantly the accuracy. Thus, the product $\mbf{a}(\dt/2)   \mbf{a}(\dt/2) = \mbf{a}(\dt) =  \mbf{I}_{2N}  + \delta \mbf{a}(\dt) $ produces 
\begin{equation}
\delta \mbf{a}(\dt) =2 \, \delta \mbf{a}(\dt/2) +  \delta \mbf{a}(\dt/2)  \cdot \delta \mbf{a}(\dt/2)
\label{eq035} 
\end{equation}
Consider then the reduced time-step introduced above, i.e. $\dt_0 = \dt/2^p$,  the first step is to compute $\delta \mbf{a}(\dt_0)$. Note that this latter is very small with $\norm{\delta \mbf{a}(\dt_0)} \ll 1$, so that using the iterative process of  Eq.~\eqref{eq035} results much more effective from a round-off error perspective than using the total matrices $ \mbf{a}(\dt)$. The iterative process consists of $p$ steps, namely
\begin{eqnarray}
\delta \mbf{a}(\dt/2^{p-1}) &=& 2 \, \delta \mbf{a}(\dt/2^p) +  \delta \mbf{a}(\dt/2^p)  \cdot \delta \mbf{a}(\dt/2^p) \nonumber \\
\delta \mbf{a}(\dt/2^{p-2}) &=& 2 \, \delta \mbf{a}(\dt/2^{p-1}) +  \delta \mbf{a}(\dt/2^{p-1})  \cdot \delta \mbf{a}(\dt/2^{p-1}) \nonumber \\
\vdots & & \quad  \text{($p$ times)} \qquad \vdots \nonumber \\
\delta \mbf{a}(\dt) &=& 2 \, \delta \mbf{a}(\dt/2) +  \delta \mbf{a}(\dt/2)  \cdot \delta \mbf{a}(\dt/2) \label{eq036}
\end{eqnarray}
and finally the main matrix of the algorithm is determined as $ \mbf{a}(\dt) =   \mbf{I}_{2N} + \delta\mbf{a}(\dt)$. For practical applications, it is sufficient to take $p=20$ to obtain results with the required accuracy. \\

In order to complete the procedure, it remains to describe how to compute the incremental matrix $\delta\mbf{a}(\dt_0) =  \mbf{a}(\dt_0) -  \mbf{I}_{2N} $. As Eq.~\eqref{eq081a} shows, the computation of $\mbf{a}(\dt_0) $ requires the truncation of two series: (i) those of $\mathbf{T}, \ \bm{\alpha}$ and $\bm{\beta}$ and (ii) the Neumann series corresponding to $(\mbf{I}_{2N} - \bm{\beta})^{-1}$. This lead us to introduce the following two parameters
\begin{itemize}
	\item [(i)] The parameter $m_a$ stands for the truncation order of matrices $\mbf{T}$, $\bm{\alpha}$, $\bm{\beta}$.   For the sake of convenience, we will use  a common criterion,  that is, fixed the parameter $m_a$ the three series will have the same number of terms. In order to avoid confusion between the exact values and those approximated by truncation, the subscript $(\bullet)_a$ will be used for the latter. Thus, given an even integer $m_a$ we introduce the following approximations
	\begin{eqnarray}
	\mbf{T}  &\approx & \mbf{I}_{2N} + \delta \mbf{T}_a(\dt_0) = \mbf{I}_{2N}  + 
	\left[ 
	\begin{array}{rr}
	\delta \mbf{G}_a(\dt_0)  & \mbf{H}_a(\dt_0) \\
	- \mbf{A} \, \mbf{H}_a(\dt_0)   & 	\delta \mbf{G}_a(\dt_0) 
	\end{array}
	\right]  \label{eq077} \\
	\delta  \mbf{G}_a(\dt_0) & = & \sum_{j=1}^{m_a/2} \frac{(-1)^{j} }{(2j)!} \p{\dt_0^2 \mbf{A}}^{j}  \quad , \quad 
	\mbf{H}_a(\dt_0) = \dt_0  \, \sum_{j=0}^{m_a/2} \frac{(-1)^{j} }{(2j+1)!} \p{\dt_0^2 \mbf{A}}^{j}  \nonumber \\
	\bm{\alpha} &\approx & \bm{\alpha}_a =  \sum_{j=0}^{m_a/2} \bm{\alpha}_j(\dt_0) \otimes \p{\mbf{A}^{j} \, \mbf{M}^{-1} \mbf{C}  } 
	= \cor{\sum_{j=0}^{m_a/2} \bm{\alpha}_j(\dt_0) \otimes \mbf{A}^{j} } \cdot  \p{\mbf{I}_2 \otimes \, \mbf{M}^{-1} \mbf{C}  }  \nonumber \\
	\bm{\beta} &\approx & \bm{\beta}_a  = \sum_{j=0}^{m_a/2} \bm{\beta}_j(\dt_0) \otimes \p{\mbf{A}^{j} \, \mbf{M}^{-1} \mbf{C}  } 
	= \cor{\sum_{j=0}^{m_a/2} \bm{\beta}_j(\dt_0) \otimes \mbf{A}^{j} } \cdot  \p{\mbf{I}_2 \otimes \, \mbf{M}^{-1} \mbf{C}  }  
	\label{eq075}
	\end{eqnarray}
	where $\bm{\alpha}_j(\dt_0), \bm{\beta}_j(\dt_0) \in \mathbb{R}^{2\times 2}$ are defined in Eqs.~\eqref{eq038} and~\eqref{eq039} for $j \geq 0$. 
	Since the order of $\dt_0$ grows by two, a truncation up to $j=m_a/2$ means that the order of $\dt_0$ will be at least $m_a$.  
	\item [(ii)] The other parameter needed to compute $\delta \mbf{a}(\dt_0)$ will be designed by $r_a$ and it represents the order of truncation of the Neumann series for $(\mbf{I}_{2N} - \bm{\beta}_a)^{-1}$, yielding 
	\begin{equation}
	\p{\mbf{I}_{2N} - \bm{\beta}_a }^{-1} \approx \mbf{I}_{2N}  + \bm{\beta}_a+  \bm{\beta}_a^2 + \cdots + \bm{\beta}_a^{r_a} 
	\equiv \mbf{I}_{2N} +  \delta \bm{\beta}_a 
	\label{eq078}
	\end{equation}
	where the matrix $\delta \bm{\beta}_a = \bm{\beta}_a+  \bm{\beta}_a^2 + \cdots + \bm{\beta}_a^{r_a} $
\end{itemize}
Substituting Eqs.~\eqref{eq077}, ~\eqref{eq075} and~\eqref{eq078}   into Eq.~\eqref{eq081a}, and after some matrix products we have
\begin{eqnarray}
\mbf{a}(\dt_0) &=& \p{\mbf{I}_{2N} - \bm{\beta} }^{-1} \p{\mbf{T} \,  + \bm{\alpha} }  
\approx  \p{\mbf{I}_{2N} +  \delta \bm{\beta}_a } \p{ \mbf{I}_{2N}  + \delta \mbf{T}_a \,  + \bm{\alpha}_a }  \nonumber \\
&=&  \mbf{I}_{2N}  
+ \p{\delta \mbf{T}_a + \bm{\alpha}_a + \delta \bm{\beta}_a  + \delta \bm{\beta}_a \,  \delta \mbf{T}_a + \delta \bm{\beta}_a \, \bm{\alpha}_a }
\equiv  \mbf{I}_{2N}  + \delta \mbf{a}(\dt_0) 
\label{eq082}
\end{eqnarray}
where consequently the incremental part $\delta \mbf{a}(\dt_0)  $ is
\begin{equation}
\delta \mbf{a}(\dt_0)  
= \delta \mbf{T}_a + \bm{\alpha}_a + \delta \bm{\beta}_a  
+ \delta \bm{\beta}_a \,  \delta \mbf{T}_a + \delta \bm{\beta}_a \, \bm{\alpha}_a 
\label{eq098}
\end{equation}
Taking $p = 20$, then $\dt_0 = \dt / 1 \, 048 \, 576$ is a very small quantity so that it is sufficient to take $m_a=2$ to achieve a highly precise  estimation of $\delta \mbf{a}(\dt_0)$. Furthermore, as will be shown later, the case $m_a=2$ behaves very favorably respect to the stability for any value of the damping forces, even the undamped case. The form of the matrices is transcribed below for this particular case ($m_a=2$). 
\begin{eqnarray}
\delta \mbf{T}_a  & =  & 
\left[ 
\begin{array}{rr}
- \frac{\dt_0^2}{2} \mbf{A} 						   & \dt_0 \mbf{I}_N - \frac{\dt_0^3}{6} \mbf{A} \\
- \dt_0 \mbf{A} + \frac{\dt_0^3}{6} \mbf{A}^2   & 	- \frac{\dt_0^2}{2} \mbf{A} 	
\end{array}
\right] 
\label{eq099} \\
\bm{\alpha}_a   & =  & 
\left[ 
\begin{array}{rr}
\frac{\dt_0}{2} \mbf{M}^{-1}\mbf{C} 	- \frac{\dt_0^3}{30} \mbf{A} \mbf{M}^{-1}\mbf{C} 
& - \frac{\dt_0^2}{12} \mbf{M}^{-1}\mbf{C} 	+ \frac{\dt_0^4}{60} \mbf{A} \mbf{M}^{-1}\mbf{C}  \\
\mbf{M}^{-1}\mbf{C} 	- \frac{3\dt_0^2}{20} \mbf{A} \mbf{M}^{-1}\mbf{C}  
&  \frac{\dt_0^3}{20} \mbf{A}   \mbf{M}^{-1}\mbf{C}  
\end{array}
\right] 
\label{eq100} \\
\bm{\beta}_a   & =  & 
\left[ 
\begin{array}{rr}
- \frac{\dt_0}{2} \mbf{M}^{-1}\mbf{C} 	+ \frac{\dt_0^3}{30} \mbf{A} \mbf{M}^{-1}\mbf{C} 
&  \frac{\dt_0^2}{12} \mbf{M}^{-1}\mbf{C} 	- \frac{\dt_0^4}{120} \mbf{A} \mbf{M}^{-1}\mbf{C}  \\
- \mbf{M}^{-1}\mbf{C} 	+ \frac{3\dt_0^2}{20} \mbf{A} \mbf{M}^{-1}\mbf{C}  
&  -\frac{\dt_0^3}{30} \mbf{A}   \mbf{M}^{-1}\mbf{C}  
\end{array}
\right] 
\label{eq101} \\
\delta \bm{\beta}_a  &=& \bm{\beta}_a+  \bm{\beta}_a^2 + \cdots + \bm{\beta}_a^{r_a} 
\label{eq102}
\end{eqnarray}
%
The Algorithm \ref{alg01} outlines the steps necessary to obtain the matrix $\mathbf{a}(\dt)$. Since in general $\dt_0 \ll T$, the two parameters $m_a$ and $r_a$ will be taken as $m_a = r_a = 2$. As will be shown later, the scheme is conditionally stable for these values for systems with nonzero damping. 
\begin{algorithm}
	\caption{Computation of main matrix $\mathbf{a}(\Delta t)$}
	\begin{algorithmic}[1]
		\footnotesize 
		\State Fix main parameters, $\Delta t$, $r_a = 2$, $m_a=2$
		\State Compute matrices $\mathbf{M}^{-1}\mathbf{C}$, $\mathbf{A}=\mathbf{M}^{-1}\mathbf{K}$
		\State Evaluate $\Delta t_0 = \Delta t / 2^p \ (p=20)$
		\State Compute matrix $\delta \mathbf{T}_a$, \quad Eq. \eqref{eq099}
		\State Compute matrix $\bm{\alpha}_a$, \quad Eq.  \eqref{eq100}
		\State Compute matrix $\bm{\beta}_a$, \quad Eq.  \eqref{eq101}
		\State Compute matrix $\delta \bm{\beta}_a = \bm{\beta}_a + \bm{\beta}_a^2$, \quad Eq.  \eqref{eq102}
		\State Compute matrix $\delta \mathbf{a}(\Delta t_0) = \delta \mbf{T}_a + \bm{\alpha}_a + \delta \bm{\beta}_a  
		+ \delta \bm{\beta}_a \,  \delta \mbf{T}_a + \delta \bm{\beta}_a \, \bm{\alpha}_a $, \quad Eq.  \eqref{eq098}
		\State Initialize   $\delta \mathbf{a} := \delta \mathbf{a}(\Delta t_0)$
		\For {$j=1\ldots q$}
		\State $\delta \mathbf{a} = 2 \delta \mathbf{a} + \delta \mathbf{a} \cdot  \delta \mathbf{a}$, \quad Eqs . \eqref{eq036}
		\EndFor
		\State Update $\delta \mathbf{a}(\Delta t) = \delta \mathbf{a}$
		\State Evaluate $\mathbf{a}(\Delta t) = \mathbf{I}_{2N} +  \delta \mathbf{a}(\Delta t)$
	\end{algorithmic}
	\label{alg01}
\end{algorithm}
\subsection{Computation of vector $\mbf{b}_k(\dt)$}

The vector $\mbf{b}_k(\dt)$ is directly related to the nonhomogeneous terms of the transient problem. As shown in Eq.~\eqref{eq081b}, it depends on the inverse matrix $(\mbf{I}_{2N} - \bm{\beta})^{-1}$ and on the matrix $\mathbf{L}$ defined in turn in terms of a matrix series, as shown in Eqs.~\eqref{eqMatrixL} and \eqref{eq014b}. As for the main matrix $\mathbf{a}(\dt)$ studied above, the two series expansions in $\mathbf{b}_k$ must be truncated for numerical practice. For that, let us introduce again two parameters: 
\begin{itemize}
	\item [(i)] $m_b$: stands for the order of truncation of the matrix series $\bm{\beta}$ and $\mathbf{L}$. In more detail we have, only with respect to the calculation of $\mathbf{b}_k$
	\begin{equation}
	\mbf{b}_k = \p{\mbf{I}_{2N} - \bm{\beta} }^{-1} \mbf{L} \, \mbf{g}_k \approx 
	\p{\mbf{I}_{2N} - \bm{\beta}_b }^{-1} \mbf{L}_b \, \mbf{g}_k \quad , \quad k = 0,1,2,\ldots \label{eq084}
	\end{equation}
	where
	\begin{equation}
	\mbf{L}_b  =  \sum_{j=0}^{m_b/2} \bm{l}_j(\Delta t) \otimes \mbf{A}^j  \ , \qquad
	\bm{\beta}_b  =  \sum_{j=0}^{m_b/2} \bm{\beta}_j(\dt) \otimes \p{\mbf{A}^{j} \, \mbf{M}^{-1} \mbf{C}  } 
	\label{eq089}
	\end{equation}
	with $ \bm{l}_j(\Delta t) $ and $ \bm{\beta}_j(\dt)  $,  $2 \times 2$ matrices defined in Eqs.~\eqref{eq014b} and \eqref{eq039}, respectively. It is important to highlight that both $\mbf{L}_b$ and $	\bm{\beta}_b$ are computed using the time step of the problem $\dt$ and not the reduced one, $\dt_0$. As before, the parameter $m_b$ should be an even number because $m_b/2$ denotes the number of terms to be taken. 
	\item [(ii)]  $r_b$: denotes the order of truncation of the Neumann series associated to the evaluation of $(\mbf{I}_{2N} - \bm{\beta}_b)^{-1} \approx \mbf{I}_{2N}  +  \bm{\beta}_b +  \bm{\beta}_b^2 + \cdots + \bm{\beta}_b^{r_b} $. Hence, the final evaluation of $\mathbf{b}_k$ leads to
	\begin{equation}
	\mbf{b}_k \approx \p{\sum_{n=0}^{r_b} \bm{\beta}_b^n } \mbf{L}_b \, \mbf{g}_k  
	\label{eq088}
	\end{equation}
	In order to optimize the number of matrix products, the above Neumann series can be computed using the algorithm
	\begin{equation}
	\mathbf{I}_{2N} + \bm{\beta} + \bm{\beta}^2 + \cdots + \bm{\beta}^{r} + \cdots = 
	\mathbf{I}_{2N} + \bm{\beta} + \bm{\beta}^2 ( \mathbf{I}_{2N} + \bm{\beta} +  \bm{\beta}^2 ( \mathbf{I}_{2N} + \bm{\beta} +   \cdots))
	\end{equation}
	requiring consequently just $r_b/2$ matrix products. To ensure that the above sum converges, we know that the spectral radius of $\bm{\beta}_b$ must verify $\rho(\bm{\beta}_b) <1$. We have already checked in Fig.~\ref{fig01} that the higher $m_b$, the wider the range of $\dt$ within the eigenvalues of $\bm{\beta}_b$ are bounded. In general, increasing the parameter $m_b$ will increase the accuracy of the results and decrease the spectral radius of $\bm{\beta}_b$, however it will increase the computational cost. A more detailed study of the latter as a function of these parameters will be carried out later. 
\end{itemize}

In the Algorithm \ref{alg02} the  necessary steps to obtain the matrix $\mathbf{b}_k(\dt)$ are listed. As shown by the numerical experiments, the parameters $m_b$ and $r_b$ result to be significant for the accuracy in certain cases, hence no particular value are recommended here. Later in the numerical examples, a study of the sensitivity of $m_b$ and $r_b$ will be carried out. 
\begin{algorithm}
	\caption{Computation of main matrix $\mathbf{b}_k(\Delta t)$}
	\begin{algorithmic}[1]
		\footnotesize 
		\State Fix main parameters, $\Delta t$ and $r_b$, $m_b$ (even numbers)
		\State Compute matrices $\mathbf{M}^{-1}\mathbf{C}$, $\mathbf{A}=\mathbf{M}^{-1}\mathbf{K}$
		\State Compute and store matrices $\mathbf{A}^{j}, \ j=2,\ldots,m_b/2$ 
		\State Compute and store matrices $\mathbf{A}^{j} (\mathbf{M}^{-1}\mathbf{C}), \ j=2,\ldots,m_b/2$
		\State Compute sequence of matrices 
		$\bm{l}_j(\dt) \in \mathbb{R}^{2 \times 4} \ , \ 1 \leq j \leq m_b/2$, \quad Eq. \eqref{eq014b}
		\State Compute sequence of matrices 
		$\bm{\beta}_j(\dt) \in \mathbb{R}^{2 \times 2} \ , \ 1 \leq j \leq m_b/2$, \quad Eq. \eqref{eq039}
		\State Compute matrix $\mathbf{L}_b$, \quad Eq.  \eqref{eq089}
		\State Compute matrix $\bm{\beta}_b$, \quad Eq.  \eqref{eq089}
		\State Compute matrix $\mathbf{B} = \mathbf{I}_{2N} + \bm{\beta}_b + \bm{\beta}_b^2$
		\If {$r_b \geq 4$}
		\For  {$n=2\ldots r_b/2$}
		\State   $\mathbf{B} = \mathbf{I}_{2N} + \bm{\beta}_b + \bm{\beta}_b^2 \, \mathbf{B}$
		\EndFor
		\EndIf
		\State Compute vector $\mathbf{g}_k$, \quad Eq. \eqref{eq083} 
		\State Compute $\mathbf{b}_k(\Delta t) = \mathbf{B} \, \mbf{L}_b \, \mbf{g}_k $, \quad Eq. \eqref{eq088} 
	\end{algorithmic}
	\label{alg02}
\end{algorithm}

\section{Algorithm stability}
\label{Stability}

\begin{figure}[ht]%
	\begin{center}
		\begin{tabular}{ccc}
			\includegraphics[width=7.6cm]{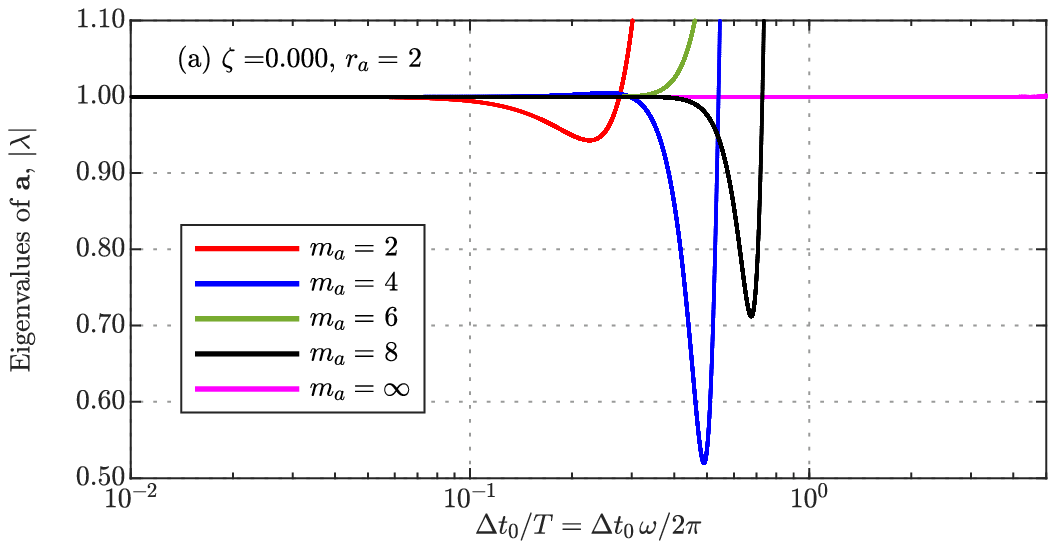} & 
			\includegraphics[width=7.6cm]{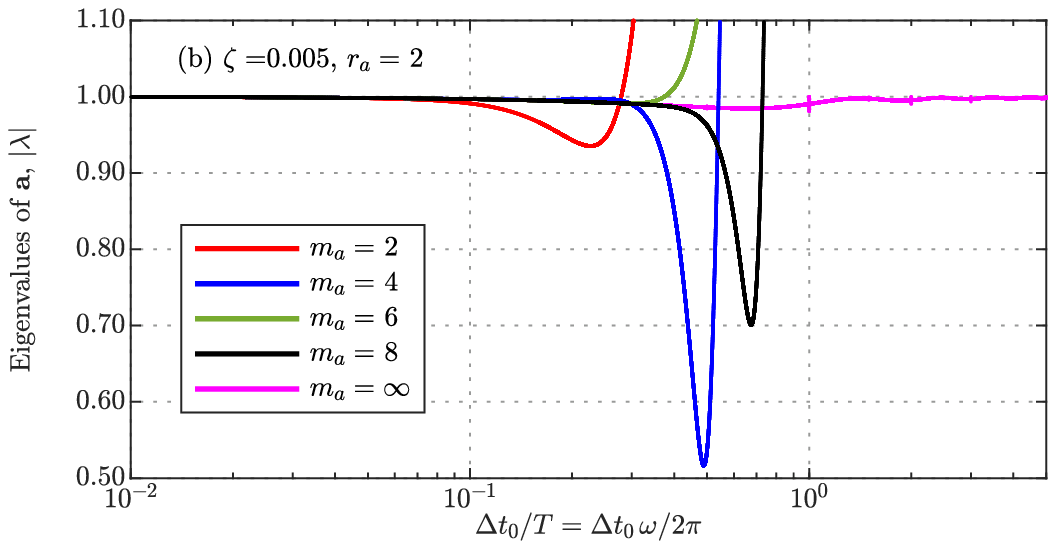} \\			\\
			\includegraphics[width=7.6cm]{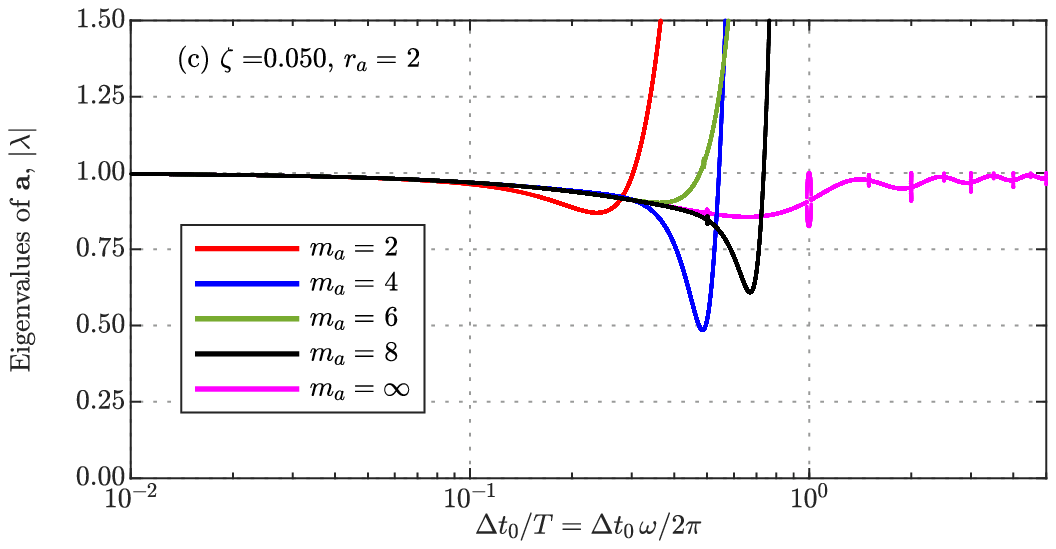} & 
			\includegraphics[width=7.6cm]{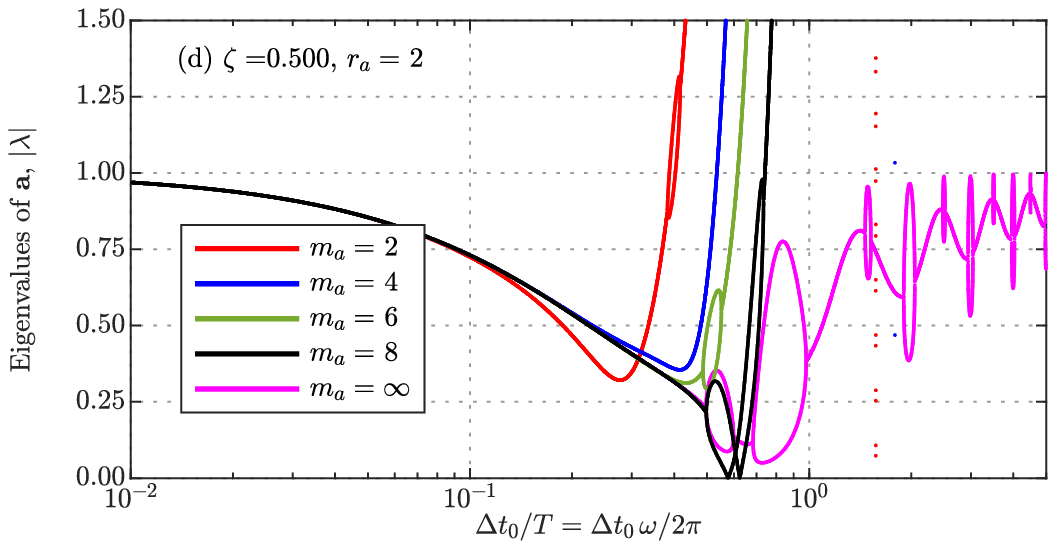}  \\ 
			\includegraphics[width=7.6cm]{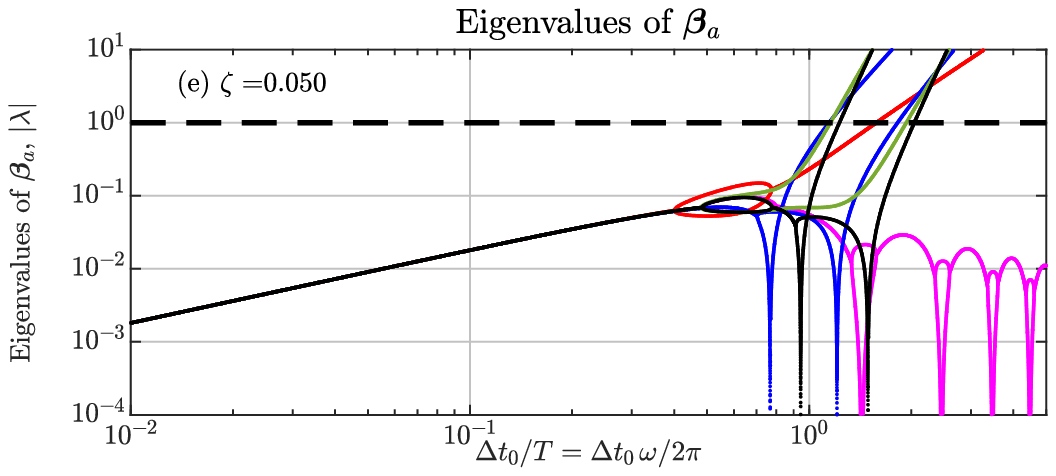} & 
			\includegraphics[width=7.6cm]{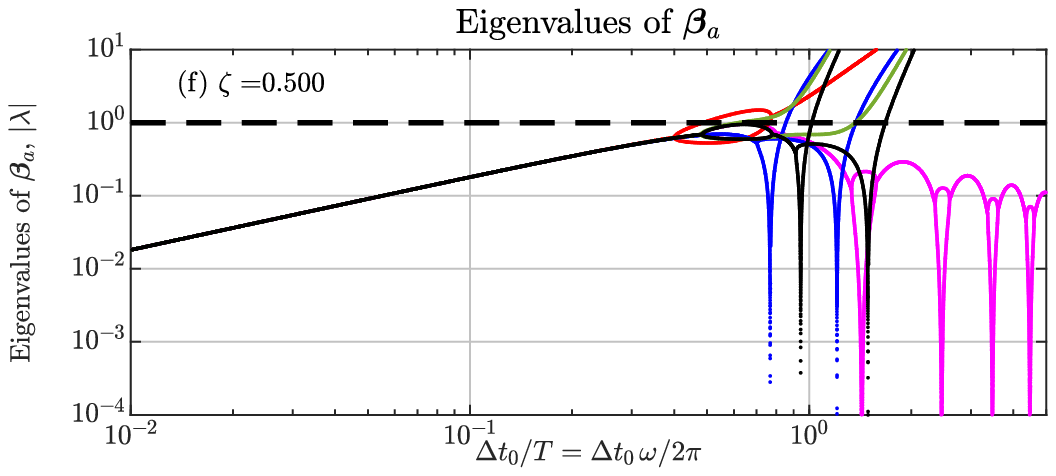} 
		\end{tabular}
		\caption{Absolute value of eigenvalues of matrix $\mathbf{a}(\dt_0) = (\mbf{I}_{2} - \bm{\beta} )^{-1} \p{\mbf{T} \,  + \bm{\alpha} } $ in a single dof system with natural frequency $\omega= \sqrt{K/M} = 2 \pi / T$ and damping ratio $\zeta = c/2M\omega$, for different cases of damping level: (a) $\zeta=0$, (b) $\zeta = 0.005$, (c) $\zeta = 0.05$ and (d) $\zeta = 0.50$. Abscissas represent the reduced time step $\dt_0/T$. Curves for different truncation orders $m_a$ are plotted in colors: $m_a = 2$ (red), $m_a = 4$ (blue), $m_a = 6$ (green), $m_a = 8$ (black), $m_a = \infty$ (magenta). The inverse matrix is approximated by $(\mathbf{I} - \bm{\beta})^{-1} \approx \mathbf{I} + \bm{\beta}_a + \bm{\beta}_a^2$ ($r_a=2$). Plots (e) and (f): absolute value of eigenvalues of matrix $\bm{\beta}_a$ for two damping cases: (e) $\zeta = 0.05$, (f) $\zeta = 0.50$}
		\label{fig02}%
	\end{center}
\end{figure}

The proposed scheme
\begin{equation}
\mbf{U}_{k+1}   = \mbf{a}(\dt)  \, \mbf{U}_k+ \mbf{b}_k(\dt) \ , \quad k = 0,1,2,\ldots
\label{eq118}
\end{equation}
is stable provided that
\begin{equation}
\rho \cor{\mathbf{a}(\dt)} \leq 1
\label{eq119}
\end{equation}
To discuss the stability of the method it is sufficient to study the single dof oscillator of mass $M$, spring rigidity $K$ and damping coefficient $c = 2 M  \omega \, \zeta$. Where $\omega = \sqrt{K/M} = 2\pi / T$ stands for the undamped natural frequency, $T$ for the natural period and $\zeta$ is the damping ratio.  The matrices of the system are then the numbers $\mathbf{A} = \cor{K/M} = \cor{\omega^2}$ and $\mathbf{M}^{-1}\mathbf{C} = \cor{2 \omega \zeta}$. Since the matrix $\mathbf{a}(\dt)$ is computed using the $2^p$ algorithm, the stability should be studied when it is evaluated at the reduced time step $\dt_0 = \dt / 2^p$. As described above, the order of truncation of the series involved is controlled by both the parameters $m_a$ and $r_a$. The parameter $m_a$ denotes the truncation order  of matrices $\mathbf{T}$, $\bm{\alpha}$ and $\bm{\beta}$ while $r_a$ stands for that of the Neumann series for $(\mathbf{I} - \bm{\beta})^{-1}$. We can then write
\begin{equation}
\mathbf{a}(\dt_0) = (\mathbf{I}_2 - \bm{\beta})^{-1} (\mathbf{T} + \bm{\alpha}) 
\approx \p{\sum_{n=0}^{r_a}\bm{\beta}_a^n} (\mathbf{T}_a + \bm{\alpha}_a)
\label{eq121}
\end{equation}
where the truncated matrices, denoted with the subscript $(\bullet)_a$, results
\begin{eqnarray}
\mathbf{T}_a &=& \sum_{j=0}^{m_a/2} \frac{(-1)^j \tau^{2j}}{(2j+1)!} 
\left[ 
\begin{array}{cc}
2j + 1		&		\dt_0		\\
- \tau^2 / \dt_0	&	2j +1 
\end{array}
\right]  \nonumber \\
\bm{\alpha}_a &=& \sum_{j=0}^{m_a/2} \frac{(-1)^j \tau^{2j+1}\zeta}{(2j+4)!} 
\left[ 
\begin{array}{cc}
24(j + 1)		&		-4(2j+1)(j+1) \dt_0		\\
24(2j+1)(j+2) / \dt_0	&	-8 (2j +1) (j+2) 
\end{array}
\right]  \nonumber \\
\bm{\beta}_a &=& \sum_{j=0}^{m_a/2} \frac{(-1)^j \tau^{2j+1}\zeta}{(2j+4)!} 
\left[ 
\begin{array}{cc}
-24(j + 1)		&		4(2j+1) \dt_0		\\
-24(2j+1)(j+2) / \dt_0	&	16 j  (j+2) 
\end{array}
\right] 
\label{eq120}
\end{eqnarray}
The  variable $\tau = \dt_0 \omega = 2\pi \dt_0 / T$ represents a dimensionless measure of the reduced time step. Eigenvalues $\lambda$ of matrix $\mathbf{a}(\dt_0)$ can be determined as function of the nondimensional time step $\dt_0 / T = \tau / 2 \pi$. As shown in Eqs.~\eqref{eq120}, $\abs{\lambda}$ depends also on the parameters $m_a$, $r_a$ and on the damping ratio $\zeta$. Since the reduced time step necessary to compute $\mathbf{a}(\dt_0)$ is a very small number, i.e. $\dt_0 = \dt / 2^p$, with $p=20$, then in general the inequality $\rho(\bm{\beta}_a) \ll 1$ holds and only a few terms in the Neumann series are necessary. In fact, considering $r_a = 2$ leads to very accurate results in the majority of cases, something that is supported by the numerical experiments. Figs.~\ref{fig02}(a),  \ref{fig02}(b),  \ref{fig02}(c) and ~\ref{fig02}(d) show the absolute value of eigenvalues of matrix $\mathbf{a}(\dt_0)$ corresponding to four cases of damping intensity, from $\zeta = 0$ (undamped case) to $\zeta = 0.5$, and for several truncation orders from $m_a = 2$ to $m_a = 8$. Intervals of stability, found as intersections of curves with $\abs{\lambda}=1$, are listed in Table \ref{tab05} in terms of $\dt_0/T$. For the undamped case ($\zeta=0$), the orders of truncation for $m_a = 4$ and $m_a=8$ are unstable for small values of $\dt_0$. However the case $m_a=2$ results to be stable for the undamped case in the range $0 < \dt_0 / T < 0.2757$. Therefore, taking $p=20$, the stability extends up to the limits $0 < \dt / T < 2.89 \times 10^5$,  much wider than the usual range of $\dt$ used to guarantee a minimum of accuracy in real problems. As Figs.~\ref{fig02} show, dissipative forces are in general favorable for stability in the proposed scheme. Most conflictive cases in terms of stability are those with null or very light damping, as for example $m_a = 4$ and $m_a = 8$ (see Table~\ref{tab05}). The case $m_a=2$ behaves very favorably, showing conditional stability in an interval that increases slightly as the damping increases.  \\

In figs.~\ref{fig02}(e) and~\ref{fig02}(f), the absolute value of eigenvalues of matrix $\bm{\beta}_a$ are plotted for two cases of damping $\zeta = 0.05$ and $\zeta = 0.5$. Since $\bm{\beta}_a$ is directly proportional to $\zeta$, the spectral radius decreases also proportionally to the damping level. For our approach to be convergent, we need the spectral radius of the $\bm{\beta}_a$ matrix to be less than unity. The convergence interval, with respect to the computation of $\mathbf{a}(\dt_0)$, can be specified as the range of values of $\dt_0/T$ for which $\rho(\bm{\beta}_a)<1$ is verified. Hence, highly damped structures could have an interval of convergence narrower than the interval of stability, in terms of $\dt_0 / T$. However, due to the very small size of the reduced time step, this point is not relevant and convergence will still be guaranteed for most practical cases.

\begin{table}
	\begin{center}
		{\footnotesize \begin{tabular}{llrr}
Damping ratio 			&			Truncation order				&			Stability boundaries			\\
\hline
$\zeta = 0.000$				  & 		  $m_a = 2$										 &           $0.0000 < \dt_0/T  <  0.2757$  \\
					    				 & 		  $m_a = 4$										 &           $0.2964 < \dt_0/T  <  0.5405$  \\
					    				 & 		  $m_a = 6$										 &           $0.0000 < \dt_0/T  <  0.2808$  \\
					    				 & 		  $m_a = 8$										 &           $0.2749 < \dt_0/T  <  0.7279$  \\
\hline
$\zeta = 0.005$				  & 		  $m_a = 2$										 &         $0.0000 < \dt_0/T  <  0.2791$  \\
										& 		  $m_a = 4$										 &           $0.0000 < \dt_0/T  <  0.5406$  \\
										& 		  $m_a = 6$										 &           $0.0000 < \dt_0/T  <  0.3741$  \\
										& 		  $m_a = 8$										 &           $0.0000 < \dt_0/T  <  0.7287$  \\
\hline
$\zeta = 0.050$				  & 		  $m_a = 2$										 &         $0.0000 < \dt_0/T  <  0.3024$  \\
										 & 		  $m_a = 4$										 &           $0.0000 < \dt_0/T  <  0.5421$  \\
										 & 		  $m_a = 6$										 &           $0.0000 < \dt_0/T  <  0.4766$  \\
										 & 		  $m_a = 8$										 &           $0.0000 < \dt_0/T  <  0.7343$  \\
\hline
$\zeta = 0.500$				  & 		  $m_a = 2$										 &         $0.0000 < \dt_0/T  <  0.3871$  \\
										& 		  $m_a = 4$										 &           $0.0000 < \dt_0/T  <  0.5342$  \\
										& 		  $m_a = 6$										 &           $0.0000 < \dt_0/T  <  0.6156$  \\
										& 		  $m_a = 8$										 &           $0.0000 < \dt_0/T  <  0.7407$  \\
\hline
\end{tabular}}
	\end{center}
	\caption{Intervals of stability in terms of $\dt_0/T$, found as intersections of curves with $\abs{\lambda}=1$. The reduced time step is $\dt_0 = \dt / 2^p$, with $p=20$. $\zeta$ stands for the damping ratio of the single dof oscillator of frequency $\omega = 2\pi / T$.}
	\label{tab05}
\end{table}

\section{Order of accuracy}
\label{Accuracy}

The order of accuracy of an iterative time-step numerical method is the minimum order of $\dt$ in the error term. Denoting by $\bm{\mathcal{U}}(t_k)$ to the $2N$-size array with the exact displacements and velocities at $t=t_k$ and by $\mbf{U}_k$ to the proposed numerical approach, then the order of accuracy is said to be $h$ if
\begin{equation}
\bm{\mathcal{U}}(t_k) = \mbf{U}_k + \bm{\mathcal{O}}(\dt^h)
\label{eq105}
\end{equation}
where $\bm{\mathcal{O}}(\dt^h)$ represents a $2N$ vector expressed as an infinite series in terms of $\dt$ whose minimum exponent is precisely $h$. Throughout this section the notation $\bm{\mathcal{O}}(\bullet)$  will be used interchangeably to refer the order in the error terms. It may represent either vectors or matrices. The context of each expression will allow us to deduce what it refers to in each case. Along the developed methodology several assumptions have been considered, leading to different sources of error. The objective of this section is to address the order of the error associated to each approximation. The analysis will be carried out considering separately each assumption, deriving the corresponding order of magnitude of the error. 
\subsection{Approximation of $\mathbf{x}^{(n-1)}(t)$ in $[t_k,t_{k+1}]$ by interpolation}
\label{Appr_xn}

In Eqs. \eqref{eq016}, the derivation of terms $\mathbf{x}_{k+1}^{(n)}, \ \dot{\mathbf{x}}_{k+1}^{(n)}$ requires to solve the integrals with the transient solution of the $(n-1)$th iteration, $\mathbf{x}^{(n-1)}(t)$. To this end, the function $\mathbf{x}^{(n-1)}(t)$ is approximated by 3rd order cubic splines using displacements and velocities at both limits of the interval, i.e. $t=t_k$ and $t=t_{k+1}$. Under this assumption, Eq.~\eqref{eq018} can be written in terms of the error order as
\begin{equation}
\mbf{x}^{(n-1)} \p{t} =  \mathcal{N}_1(\xi) \, \mbf{x}_k^{(n-1)}		+ 		
\mathcal{D}_1(\xi) \, \dt \, \dot{\mbf{x}}_k^{(n-1)} + 
\mathcal{N}_2(\xi) \, \mbf{x}_{k+1}^{(n-1)}		+ 		
\mathcal{D}_2(\xi) \, \dt \, \dot{\mbf{x}}_{k+1}^{(n-1)} + \bm{\mathcal{O}}(\dt^4) 
\label{eq122}
\end{equation}
Plugging Eq.~\eqref{eq122} into the integrals of Eq.~\eqref{eq016} and following the same steps, we can rewrite Eq.~\eqref{eq021} in simplified form in terms of the error as 
\begin{equation}
\bm{\mathcal{X}}_{k+1}^{(n)} = \mathbf{X}_{k+1}^{(n)}  + 
\left\| \mathbf{M}^{-1}\mathbf{C}\right\| \, \bm{\mathcal{O}}(\dt^5) 
\ , \quad n \geq 1 \ , \quad k \geq 0
\label{eq123}
\end{equation}
where $\bm{\mathcal{X}}_{k+1}^{(n)}$ denotes the exact solution of the $n$th iteration at instant $t_{k+1}$ and 
$$
\mathbf{X}_{k+1}^{(n)} = \mbf{T} \, \mbf{X}_k^{(n)}  + \bm{\alpha} \, \mbf{X}_k^{(n-1)}  +\bm{\beta} \, \mbf{X}_{k+1}^{(n-1)}
$$
stands for the numerical approximate value. The exact solution is by definition
\begin{equation}
\bm{\mathcal{U}}(t_k) = \sum_{n=0}^\infty \bm{\mathcal{X}}_k^{(n)} 
\label{eq124}
\end{equation}
Therefore, from Eqs.~\eqref{eq123} and \eqref{eq124}, it can be established that the order of the error due to the interpolation of $\mathbf{x}^{(n-1)}(t)$ is proportional, something that can be written as
\begin{equation}
\left\| \bm{\mathcal{U}}(t_k) - \mathbf{U}_k \right\| \leq  \left\| \mathbf{M}^{-1}\mathbf{C}\right\| \, \bm{\mathcal{O}}(\dt^5) 
\label{eq125}
\end{equation}
for some norm $\left\| \bullet \right\| $. The above order of accuracy only measures the error of the approximation given by Eq.~\eqref{eq122}.  Truncation of series involved in matrices $\mathbf{a}(\dt)$ and $\mathbf{b}_k(\dt)$ also give rise to errors, whose order of magnitude will be estimated in the following points.

\subsection{Approximation of $\delta \mathbf{a}(\Delta t_0)$ by  truncation}
\label{Appr_a}

As stated in Sec. \ref{PreciseComp_a},  truncation of series is needed for computation of the incremental part, $\delta \mbf{a}(\dt_0) = \mathbf{a}(\dt_0) - \mbf{I}_{2N}$ in Eq.~\eqref{eq098}. Let us denote by $\delta \bm{a}(\dt_0)$ to the exact incremental matrix after retaining all terms of the series and by $\delta \mathbf{a}(\dt_0)$ to the matrix numerically determined. Remind that $\dt_0 = \dt / 2^p$ represents the reduced time step, with $p=20$. Then, Eq.~\eqref{eq082} can be written as
\begin{equation}
\delta \bm{a}(\dt_0) =  \p{\mathbf{I}_{2N} - \bm{\beta} }^{-1} \p{ \mathbf{T} + \bm{\alpha}} - \mathbf{I}_{2N}
= \delta \mbf{T} + \bm{\alpha} + \delta \bm{\beta}  
+ \delta \bm{\beta} \,  \delta \mbf{T} + \delta \bm{\beta} \, \bm{\alpha}
\label{eq106}
\end{equation}
where $\delta \mbf{T} , \ \bm{\alpha} , \ \delta \bm{\beta} $ stand for the exact version of the matrices, which in turn can be expanded as function of the approximate  matrices using the parameter $m_a$, yielding
\begin{eqnarray}
\delta \mbf{T}  & =  &  \delta \mbf{T}_a + 
\left[ 
\begin{array}{rr}
\bm{\mathcal{O}}(\dt_0^{m_a + 2})						   & \bm{\mathcal{O}}(\dt_0^{m_a + 3}) \\
\bm{\mathcal{O}}(\dt_0^{m_a + 3})  & 	\bm{\mathcal{O}}(\dt_0^{m_a + 2})
\end{array}
\right] 
\ , \quad 
\bm{\alpha}    =   \bm{\alpha}_a + 
\left[ 
\begin{array}{rr}
\bm{\mathcal{O}}(\dt_0^{m_a + 3})  & \bm{\mathcal{O}}(\dt_0^{m_a + 4}) \\
\bm{\mathcal{O}}(\dt_0^{m_a + 2})  & 	\bm{\mathcal{O}}(\dt_0^{m_a + 3})
\end{array}
\right] 
\nonumber \\
\delta \bm{\beta}  &=& \sum_{n=1}^\infty \bm{\beta}^n 
=  \sum_{n=1}^{r_a} \bm{\beta}^{n} + \bm{\mathcal{O}}\p{\left\|\bm{\beta} \right\|^{r_a+1} } 
\ , \quad 
\bm{\beta}   =   \bm{\beta}_a + 
\left[ 
\begin{array}{rr}
\bm{\mathcal{O}}(\dt_0^{m_a + 3})  & \bm{\mathcal{O}}(\dt_0^{m_a + 4}) \\
\bm{\mathcal{O}}(\dt_0^{m_a + 2})  & 	\bm{\mathcal{O}}(\dt_0^{m_a + 3})
\end{array}
\right] 
\label{eq104} 
\end{eqnarray}
The minimum order of the error in the three matrices is $\bm{\mathcal{O}}(\dt_0^{m_a + 2})$. Hence,  plugging Eq. \eqref{eq104} into Eq.~\eqref{eq106} we have
\begin{multline}
\delta \bm{a}(\dt_0)  
= \delta \mbf{a}(\dt_0) + \bm{\mathcal{O}}(\dt_0^{m_a + 2}) + 				\left\|\mbf{M}^{-1} \mbf{C} \right\|^{r_a+1} \bm{\mathcal{O}}\p{\dt_0^{r_a+1}  } 
=  \\
\delta \mbf{a}(\dt_0) + \bm{\mathcal{O}} \p{\frac{\dt^{m_a + 2}}{2^{p(m_a + 2)}}} + 
\left\|\mbf{M}^{-1} \mbf{C} \right\|^{r_a+1} \, \bm{\mathcal{O}} \p{\frac{\dt^{r_a + 1}}{2^{p(r_a + 1)}}}
\label{eq103}
\end{multline}
Considering the common values $m_a=2$, $p=20$ and $r_a=2$, the error in the computation of $\mbf{a}(\dt)$ results to be of the order
%
\begin{equation}
\norm{\delta \bm{a}(\dt_0)  -  \delta \mbf{a}(\dt_0) } \leq 
8.27 \times 10^{-25} \bm{\mathcal{O}} \p{  \dt^{4}} + 
8.67 \times 10^{-19} \left\|\mbf{M}^{-1} \mbf{C} \right\|^{3} \,  \bm{\mathcal{O}} \p{  \dt^{3}}
\end{equation}
which represents a very small error in terms of computer accuracy.

\subsection{Approximation of $\mathbf{b}_k(\dt)$ by  truncation}
\label{Appr_b}

We will now study the order of the error committed in the computation of the vector $\mbf{b}_k(\dt)$, associated to the nonhomogeneous solution. In the context of this section, we will denote by $\bm{b}_k(\dt) $ the exact nonhomogeneous term and reserve the notation $\mbf{b}_k(\dt)$ for the approximate numerical result. In this vector, we find several sources of error: 
\begin{itemize}
	\item [(i)] In Eq.~\eqref{eq009b} the vector $\mathbf{f}(t)$ is approximated in $\cor{t_k,t_{k+1}}$ by 3rd order Laplace polynomials.
	\item [(ii)] in Eq.~\eqref{eq088} The inverse matrix $ (\mbf{I}_{2N} - \bm{\beta})^{-1} $ is approximated by a Neumann series up to the $r_b$--th order. 
	\item [(iii)] The truncation of both matrices $ \bm{\beta}$ and  $\mathbf{L}$ also introduces an error, controlled by the parameter $m_b$. 
\end{itemize}

First, let us see the order of the error (i), as consequence of the approximation of the external forcing vector $\mbf{f}(t)$ shown in Eq.~\eqref{eq009b} by polynomial interpolation  in the interval $[t_k,t_{k+1}]$. Since the type of interpolation functions are 3rd order Lagrange polynomials, Eq.~\eqref{eq009b} can now be written as
\begin{equation}
\mbf{f}(t) = \sum_{i=1}^4 \mathcal{L}_i \p{\frac{t-t_k}{t_{k+1} - t_k}} \, \mbf{f} \p{t_k + \frac{(i-1)\dt}{3}} + \bm{\mathcal{O}}(\dt^4) \quad , \quad t_k \leq t \leq t_{k+1} 
\label{eq107}
\end{equation}
where $\mathcal{L}_i(\xi), \ 1 \leq i \leq 4$ stand for the  3rd order Lagrange polynomials. The aim of this interpolation is the computation of the integrals within Eqs.~\eqref{eq011}, whose order of approximation results
\begin{equation}
\left\lbrace 
\begin{array}{c}
\ds \int_{t=t_k}^{t_{k+1}} \mbf{H}(t_{k+1}-t) \mbf{M}^{-1} \mbf{f}(t) \, dt  \\
\ds \int_{t=t_k}^{t_{k+1}} \mbf{G}(t_{k+1}-t) \mbf{M}^{-1} \mbf{f}(t) \, dt 
\end{array}
\right\rbrace =
\mathbf{L} \, \mathbf{g}_k +  \bm{\mathcal{O}}(\dt^5) 
\label{eq011c}
\end{equation}
where $\mathbf{g}_k$, defined in Eq.~\eqref{eq083}, contains the forcing terms at the data points $t_k + (i-1)\dt/3, \ i=1,2,3,4$.  Therefore the exact form of the vector $\bm{b}_k(\dt)$ is
\begin{equation}
\bm{b}_k(\dt) = (\mbf{I}_{2N} - \bm{\beta})^{-1} \, \cor{\mathbf{L} \, \mathbf{g}_k +  \bm{\mathcal{O}}(\dt^5) } 
\label{eq126}
\end{equation}
Now, let us see how the truncation of the series in matrices $\bm{\beta}$ and $\mbf{L}$ affect to the order of the error. Considering the parameters $r_b$ and $m_b$, 
\begin{equation}
(\mbf{I}_{2N} - \bm{\beta})^{-1}  =  \sum_{n=0}^{r_b} \bm{\beta}^{n} + \bm{\mathcal{O}}\p{\left\|\bm{\beta} \right\|^{r_b+1} } 
\quad , \quad
\bm{\beta}   =   \bm{\beta}_b + \bm{\mathcal{O}}(\dt^{m_b+2}) \ , \quad 
\mathbf{L} =\mathbf{L}_b + \bm{\mathcal{O}}(\dt^{m_b + 3})
\label{126}
\end{equation}
where $\bm{\beta}_b $ and $\mathbf{L}_b$ are defined in Eq.~\eqref{eq089}. It has been shown in Eq.~\eqref{eq060} that, if the ratio $\dt / T$ lies within certain range,  $\bm{\beta}$ can be considered  proportional to the damping matrix and to the time step. 
\begin{eqnarray}
\bm{b}_k(\dt) &=& (\mbf{I}_{2N} - \bm{\beta})^{-1} \, \cor{\mathbf{L} \, \mathbf{g}_k +  \bm{\mathcal{O}}(\dt^5) }   \nonumber \\
&=& \cor{ \sum_{n=0}^{r_b} \bm{\beta}_b^n  +  \bm{\mathcal{O}}(\dt^{m_b+2}) + \left\|\mbf{M}^{-1} \mbf{C} \right\|^{r_b+1} \bm{\mathcal{O}}\p{\dt^{r_b+1}  } }  \, \cor{\mathbf{L} \, \mathbf{g}_k + \bm{\mathcal{O}}(\dt^5)   }
\nonumber \\
&=& \cor{\sum_{n=0}^{r_b} \bm{\beta}_b^n  +  \bm{\mathcal{O}}(\dt^{m_b+2}) + \left\|\mbf{M}^{-1} \mbf{C} \right\|^{r_b+1} \bm{\mathcal{O}}\p{\dt^{r_b+1}  }}  \, 
\cor{\mathbf{L}_b \, \mathbf{g}_k   +  \bm{\mathcal{O}}(\dt^{m_b+2}) + \bm{\mathcal{O}}(\dt^5)  } \nonumber \\
&=& \p{ \sum_{n=0}^{r_b} \bm{\beta}_b^n }  \, 
\mathbf{L}_b \, \mathbf{g}_k  
+ \bm{\mathcal{O}}(\dt^{m_b+2}) +  
\left\|\mbf{M}^{-1} \mbf{C} \right\|^{r_b+1} \bm{\mathcal{O}}\p{\dt^{r_b+1}  } + \bm{\mathcal{O}}(\dt^5) \nonumber \\
& \equiv &
\mbf{b}_k(\dt) 
+ \bm{\mathcal{O}}(\dt^{m_b+2}) +  
\left\|\mbf{M}^{-1} \mbf{C} \right\|^{r_b+1} \bm{\mathcal{O}}\p{\dt^{r_b+1}  } + \bm{\mathcal{O}}(\dt^5) 
\label{eq110}
\end{eqnarray}

\begin{table}
	\begin{center}
		\scriptsize{\begin{tabular}{lllc}
Source of error			&			Order of error, $\left\| \bm{\mathcal{U}}(t_k) - \mathbf{U}_k\right\| $	 & Description\\
\hline
Sec. \ref{Appr_xn}    &           $\left\| \mathbf{M}^{-1}\mathbf{C}\right\| \, \bm{\mathcal{O}}(\dt^5)  $
																				& Approximation of $\mathbf{x}^{(n-1)}(t)$ by cubic splines \\
Sec. \ref{Appr_a}   & $\bm{\mathcal{O}} \p{\frac{\dt^{m_a + 2}}{2^{p(m_a + 2)}}} + 
\left\|\mbf{M}^{-1} \mbf{C} \right\|^{r_a+1} \, \bm{\mathcal{O}} \p{\frac{\dt^{r_a + 1}}{2^{p(r_a + 1)}}}$				& Approximation of $\mathbf{a}(\dt)$ \\	
Sec. \ref{Appr_b} & $\bm{\mathcal{O}}(\dt^{m_b+2}) + 
								\left\|\mbf{M}^{-1} \mbf{C} \right\|^{r_b+1} \bm{\mathcal{O}}\p{\dt^{r_b+1}  } + \bm{\mathcal{O}}(\dt^5) $
							&  Approximation of vector $\mathbf{b}_k(\dt)$ \\																			
\hline																				
\end{tabular}}
	\end{center}
	\caption{Order of error of the proposed method from the different sources of approximation}
	\label{tab06}
\end{table}

A summary of the order of accuracy obtained from each error source is presented in Table \ref{tab06}. For practical purposes, we will consider fixed the parameters $m_a=2$ and $r_a=2$ for the determination of $\mathbf{a}(\dt)$. In view of the results in Table \ref{tab06}, the order of approximation would be theoretically $\mathcal{O}(\dt^3)$ for these particular values, although in practice these terms are affected by such a small number that the computation  $\mathbf{a}(\dt)$ with the precise method of Sec.~\ref{PreciseComp_a} can be considered exact for computer precision. From the accuracy point of view, the other orders shown in Table \ref{tab06} are more interesting since they show a dependence on the damping level, further affected by the parameter $r_b$, which will explain the behavior of the response error depicted later in the numerical simulations.

\section{Computation effort}
\label{ComputationalEffort}

In this section an estimation of the computational cost by counting the number of required operations will be addressed. In particular the number of matrix multiplications of the proposed method will be counted and compared with other explicit methods: the Modified Precise Integration Method (MPIM) the 4th order Runge--Kutta method (RK4). The computational cost of the calculation of $\mathbf{M}^{-1}$ will not be considered, because it is a common  requirement for the three methods.

\subsection{Proposed method (PER)}

The proposed method is based on the perturbation of the damping terms, hence it will be called the PER method in this section. We will distinguish between the two clearly differentiated stages: (1) a previous set-up of the system matrices, focused on the computation of main matrices $\mathbf{a}$ and $\mathbf{b}_k$ and (2) the computation of the solution at the different time instants by the iterative procedure, $\mathbf{U}_1,\ldots, \mathbf{U}_k$, where $k = t_{\max{}}/\dt$. In order to derive a closed form of the operations count, let us establish some hypotheses
\begin{itemize}
	\item The parameters $m_a$ and $m_b$, which measure the maximum order of $\dt$ in the computation of $\mathbf{a}$ and $\mathbf{b}_k$, are both even numbers and they verify that $m_b \geq m_a \geq 2$.
	\item The series $	\mathbf{I}_{2N} + \bm{\beta} + \bm{\beta}^2 + \cdots + \bm{\beta}^r$ can be evaluated using $r/2$ matrix products by the algorithm 
	$$
	\mathbf{I}_{2N} + \bm{\beta} + \bm{\beta}^2 + \cdots + \bm{\beta}^r + \cdots= 	\mathbf{I}_{2N} + \bm{\beta} +
	\bm{\beta}^2 ( \mathbf{I}_{2N} + \bm{\beta} +  \bm{\beta}^2 ( \mathbf{I}_{2N} + \bm{\beta} +  \bm{\beta}^2 (\mathbf{I}_{2N} + \bm{\beta} +   \bm{\beta}^2( \mathbf{I}_{2N}+ \cdots))))
	$$	
	consequently the parameters $r_a$ and $r_b$ will be even numbers
\end{itemize} 
Table \ref{tab01} shows the operation counting of each step of the proposed algorithm, resulting in the total number of operations
\begin{equation}
\mathcal{C}_{\text{PER}} = \left[33 + 4(r_a + r_b) + 8p + m_b\right] \,N^3 + 
\left[ 11m_a/2 + 8 m_b + 24 + 4 t_{\max{}} / \dt \right]\, N^2 + 
(8 t_{\max{}}/\dt) \, N 
\label{eq127}
\end{equation}
Some values that can be fixed are $m_a = r_a = 2$ and $p = 20$

\begin{table}
	\begin{center}
		\begin{tabular}{llr}
	Type of solution 							&		Algorithm step		&		Operation count						\\
	\hline 
	Common matrices 		&		
	Matrices $\mathbf{M}^{-1}\mathbf{C}$, $\mathbf{M}^{-1}\mathbf{K}$  &				$2N^3$ 									\\
	& Matrices $\mathbf{A}^2,\mathbf{A}^3,\ldots,\mathbf{A}^{m_b/2}$		&				$(m_b/2 - 1) N^3$				  \\	
	& Matrices $\mathbf{A}(\mathbf{M}^{-1}\mathbf{C}),\ldots,\mathbf{A}^{m_b/2}(\mathbf{M}^{-1}\mathbf{C})$		
																													&				$(m_b/2) N^3$				  \\	
	\hline 																													
	Matrix $\mathbf{a}(\dt)$ 	&	
	Matrix $\bm{\alpha}$, products $\bm{\alpha}_j(\dt_0) \otimes \mathbf{A}^{j}(\mathbf{M}^{-1}\mathbf{C})$			& $(m_a/2 + 1) (4 N^2)$ \\
			 &
	Matrix $\bm{\beta}_a(\dt_0)$, products $\bm{\beta}_j(\dt_0) \otimes \mathbf{A}^{j}(\mathbf{M}^{-1}\mathbf{C})$			& $(m_a/2 + 1) (4 N^2)$				\\
	& Matrix $\mathbf{T}$, products $\mathbf{t}_j(\dt_0) \otimes \mathbf{A}^{j}$			& $(m_a/2) (3N^2)$ \\																			
	&	Matrix $\bm{\beta}_a^2 + \cdots + \bm{\beta}_a^{r_a}$																& $(r_a/2) \, (2N)^3$  \\
	&  Matrix $\delta \mathbf{a}(\dt_0)$																		& $ 2 \, (2N)^3 $			\\
	& Matrix $\delta \mathbf{a}(\dt)$																		     & $p \, (2N)^3 $			\\
	\hline 	
	Vector $\mathbf{b}_k(\dt)$  &
	Matrix $\mathbf{L}_b$, products $\bm{l}_j(\dt) \otimes \mathbf{A}^{j}$			& $(m_b/2) \, (8 N^2)$   \\
			&		
	Matrix $\bm{\beta}_b(\dt)$, products $\bm{\beta}_j(\dt) \otimes \mathbf{A}^{j}(\mathbf{M}^{-1}\mathbf{C})$			& $(m_b/2 + 1) (4 N^2)$				\\
	& 	Matrix $\bm{\beta}_b^2 + \cdots + \bm{\beta}_b^{r_b}$																& $(r_b/2) \, (2N)^3$  \\
	&  Vector $\mathbf{g}_k$ (assumed $k$--independent)															& $4N^2$ 			\\
	&  Product $(\sum_{j=0}^{r_b/2} \bm{\beta}_b^j ) \, \mathbf{L}_b$												 & $16N^3$				\\
	&  Product $(\sum_{j=0}^{r_b/2} \bm{\beta}_b^j ) \, \mathbf{L}_b \, \mathbf{g}_k$						& $8 N^2$				\\
	\hline
	Time-step iteration     &   Homogeneous terms															&    $(t_{\max{}} / \dt) \, (4N^2) $ \\
										&  Nonhomogeneous terms														& 	  $(t_{\max{}} / \dt) \, (8N) $ \\				
	\hline 	
	 Total       &      \multicolumn{2}{r}{$\left[33 + 4(r_a + r_b) + 8p + m_b\right] \,N^3$} \\
	       			 &    	 \multicolumn{2}{r}{$+ \left[ 11m_a/2 + 8 m_b + 24 + 4 t_{\max{}} / \dt \right]\, N^2$} \\
	        	    &       \multicolumn{2}{r}{$ + (8 t_{\max{}}/\dt) \, N $}
\end{tabular}
	\end{center}
	\caption{Operation count of the proposed method. Hypotheses: $m_a$, $m_b$, $r_a$ and $r_b$ even number, $\dt_0 = \dt / 2^p$, $m_b >  m_a\geq 2$ }
	\label{tab01}
\end{table}

\subsection{Modified Precise Integration Method (MPIM)}%
\label{PIM}

As shown in Sec.~\ref{PreciseComp_a}, a precise computation of the matrix $\mathbf{a}(\dt)$ must be carried out using the $2^p$ algorithm. These feature is shared with the well-known Precise Integration Method~\cite{Zhong-1994,Zhong-2004,WangMF-2005}, which make it a direct competitor from the computational point of view. The transient problem can be written under the state-space formalism
\begin{equation}
\der{\bm{\mathcal{U}}}{t}  = \mbf{W} \, \bm{\mathcal{U}}(t) + \mbf{h}(t) 
\label{eq090}
\end{equation}
where
\begin{equation}
\bm{\mathcal{U}}(t) = 
\left\lbrace 
\begin{array}{c}
\mbf{u}(t) \\
\dot{\mbf{u}}(t) 
\end{array}
\right\rbrace 
\quad ,
\quad
\mbf{W} =
\left[ 
\begin{array}{rr}
\mbf{O}_N & \mbf{I}_N \\
-\mbf{M}^{-1} \mbf{K} & -\mbf{M}^{-1} \mbf{C}
\end{array}
\right] \ , \quad
\mbf{h}(t) = 
\left\lbrace 
\begin{array}{c}
\mbf{0} \\
\mbf{M}^{-1} \mbf{f}(t)
\end{array}
\right\rbrace 
\label{eq092}
\end{equation}
Considering the time step $\dt$, the iterative scheme yields~\cite{Zhong-1994}
%
\begin{equation}
\mbf{U}_{k+1} = e^{\mathbf{W} \dt} \, \mbf{U}_k + 
\int_{t_k}^{t_{k+1}}   e^{\mathbf{W} (t_{k+1} - \tau)}\, \mbf{h}(\tau) \, d\tau
\label{eq091b}
\end{equation}
The main contribution of the Precise Integration Method is the proposal of a highly efficient computation of the exponential matrix. Denoting by $\bm{\mathcal{H}}(t) = e^{\mathbf{W}t}$, then 
\begin{equation}
\bm{\mathcal{H}}(\dt) =  \cor{\bm{\mathcal{H}}(\dt / 2^p)}^{2^p}
\label{eq129}
\end{equation}
where, as in Sec.~\ref{PreciseComp_a}, the reduced time step is $\dt_0 =  \dt / 2^p$. As known, in the $2^p$ algotithm, the products of Eq.~\eqref{eq129} are carried out from the incremental matrix $\delta \bm{\mathcal{H}}(\dt_0) = \bm{\mathcal{H}}(\dt_0) - \mathbf{I}_{2N} $, which in the first step is evaluated from the Taylor series expansion of the exponential matrix up to an order which guaranties the stability (the fourth order is usual~\cite{Zhong-1994,WangMF-2005}), yielding
\begin{equation}
\delta \bm{\mathcal{H}}(\dt_0) = e^{\mathbf{W}\dt_0} - \mathbf{I}_{2N} 
\approx  \mathbf{W} \dt_0 + \frac{\dt_0^2}{2!} \, \mathbf{W}^2 +  \frac{\dt_0^3}{3!} \, \mathbf{W}^3 + 
\frac{\dt_0^4}{4!} \, \mathbf{W}^4 
= \mathbf{W} \dt_0 + \frac{\dt_0^2 }{2!} \, \mathbf{W}^2 
\p{\frac{\dt_0}{3} \, \mathbf{W}   + \frac{\dt_0^2}{12} \, \mathbf{W}^2   }    
\label{eq130}
\end{equation}
From the above equation, two matrix products of size $2N$ are needed, leading to $2(2N)^3$ operations. The PIM  can currently be considered more as a family of methods, some of which propose different approaches for the determination of the non-homogeneous term. In the original paper of Zhong~\cite{Zhong-1994} a linear interpolation of the applied force was proposed, requiring the computation of the inverse matrix $\mathbf{W}^{-1}$. Since then, other solutions avoiding the computation of the inverse matrix $\mathbf{W}^{-1}$ and keeping a high accuracy in the solution have been proposed~\cite{Caprani-2013,Huang-2018,Leung-1986,Lin-1995,WangMF-2005}. Among them, the so-called Modified PIM (MPIM~\cite{WangMF-2005}) evaluates the non-homogeneous term using Gauss integration points, resulting in the following expression
%
\begin{equation}
\mbf{w}_k = \int_{t_k}^{t_{k+1}} e^{\mathbf{W}(t_{k+1} - \tau)}  \, \mbf{h}(\tau) \, d\tau \approx 
\frac{\dt}{2} \, \sum_{i=1}^g \, \mathcal{W}_i \, \exp \left[ \frac{\dt}{2}(1-\eta_i) \, \mathbf{W}\right] \,  \mbf{h} \left[ t_k + \frac{\dt}{2}(1-\eta_i)\right] 
\label{eq093}
\end{equation}
where $\eta_i,  \ 1 \leq i \leq g$ denote the Gauss integration points in the interval $-1 \leq \eta \leq 1$ and $ \mathcal{W}_i$ the corresponding weights. The MPIM presents several advantages that lead us to choose it in this article for comparison with the proposed method: (i) it reduces the non-homogeneous term to a simple expression easy to implement, (ii) it does not need the computation of the inverse matrix, (iii) it keeps a high accuracy in the calculation of the system response, (iv) it works independently of the form of the forcing vector, which does not need a preprocessing. On the other side, it needs the computation of the exponential matrix $\exp \left[ \frac{\dt}{2}(1-\eta_i) \mathbf{W}\right]  $ as many times as the number of Gauss points, i.e. $g$ times \cite{WangMF-2005}. 
\begin{table}
	\begin{center}
		\begin{tabular}{llr}
	Type of solution 							&		Algorithm step		&		Operation count						\\
	\hline 
	Matrix $\bm{\mathcal{H}}(\dt) = e^{\mathbf{W}\dt}$ 		&		
	Matrices $\mathbf{M}^{-1}\mathbf{C}$, $\mathbf{M}^{-1}\mathbf{K}$  &				$2N^3$ 									\\
	&  Matrix $\delta \bm{\mathcal{H}}(\dt_0)$	(up to 4th order)																	& $ 2 (2N)^3 $			\\
	& Matrix $\delta \bm{\mathcal{H}}(\dt )$																		     & $p \, (2N)^3 $			\\
	\hline 	
	Nonhomogeneous term  &
	Matrices  $\delta \bm{\mathcal{H}}(\dt (1 - \eta_i)/2)$, $1 \leq i \leq g$			&    $ g \, p \, (2N)^3 $   \\
	& Vectors, $\bm{\mathcal{H}}(\dt (1 - \eta_i)/2) \mathbf{h}(t_k + \dt (1 - \eta_i)/2)$	&     	$g  \, (2N)^2 $ \\
\hline 	
	Time-step iteration     &   Matrix--vector products															&    $(t_{\max{}} / \dt) \, (4N^2) $ \\
\hline 	
	Total & 	\multicolumn{2}{r}{$\left[  18 +  8p(1 + g) \right] \,N^3$} \\	
	        &  		\multicolumn{2}{r}{$\left[   4g + 4 t_{\max{}} / \dt \right] \, N^2 $}	\\
\end{tabular}
	\end{center}
	\caption{Operation count of the Modified Precise Integration Method (MPIM). Parameter $g$ stands for the number of integration points in the evaluation of the exponential matrix  (nonhomogeneous term).}
	\label{tab02}
\end{table}
In Table \ref{tab02} the number of operations involved in the matrix products of the MPIM are counted, resulting
\begin{equation}
\mathcal{C}_{\text{MPIM}} = \left[18 + 8p(1+g)\right] \,N^3 + 
\left[ 4g+ 4 t_{\max{}} / \dt \right]\, N^2 
\label{eq128}
\end{equation}

\subsection{4th order Runge-Kutta (RK4)}
\label{RK4_method}

The Runge-Kutta method approximates the value of the function $\mathbf{U}_{k+1} = \bm{\mathcal{U}}(t_k+ \dt)$ from the previous step $\bm{\mathcal{U}}(t_k)$ in the ordinary differential equation $\dot{\bm{\mathcal{U}}}(t) = \bm{\mathcal{F}}(\bm{\mathcal{U}},t)$, where the vector field 
$$
\bm{\mathcal{F}}(\bm{\mathcal{U}},t) = \mathbf{W} \, \bm{\mathcal{U}} + \mathbf{h}(t)
$$
has already been introduced in Eq.~\eqref{eq090} for oscillatory linear damped systems. In particular, the well--known fourth--order Runge--Kutta method (RK4,  \cite{Butcher-2008,Burden-2001}) achieves a truncation error at each time-step of the order $\mathcal{O}(\dt^5)$ and an overall error of order $\mathcal{O}(\dt^4)$ after four evaluations of the function $\bm{\mathcal{F}}(\mathbf{U},t)$. Since only two matrix products are necessary, indeed $\mathbf{M}^{-1}\mathbf{K}$ and $\mathbf{M}^{-1}\mathbf{C}$, the operations count of the RK4 yields
\begin{equation}
\mathcal{C}_{\text{RK4}} = 2 N^3 + \left(16 t_{\max{}} / \dt \right) N^2 + + 8 t_{\max{}} \, N / \dt 
\label{eq097}
\end{equation}

\subsection{Comparison}

In view of the computational cost analysis, the following conclusions can be drawn:
\begin{itemize}
	\item As the size of the system increases, the terms which really matters are those proportional to $N^3$. Comparing the proposed method with the MPIM, the former will be numerically more efficient than the latter if 
	\begin{equation}
	\lim_{N\to \infty} \frac{\mathcal{C}_{\text{PER}}}{\mathcal{C}_{\text{MPIM}} }= \frac{33 + 4(r_a + r_b) + 8p + m_b}{18 + 8 (1 + g) p} < 1
	\label{eq111}
	\end{equation}
	For instance, considering $m_a = r_a = 2$ and $p = 20$, Eq.~\eqref{eq111} leads to 
	\begin{equation}
	23 + m_b + 4r_b > 160 g
	\label{eq113}
	\end{equation}
	where $g$ stands for the number of Gauss integration points, which in general can be taken between 3 and 5~\cite{WangMF-2005}. 
	\item When the number of iterations is small compared to the system size, the Runge-Kutta method is the most effective. It is very favorable the fact that no previous operations are needed to build the algorithm matrices (to construct the iterative scheme only costs $2N^3$ multiplications). However, it is known that the RK4 method has stability boundaries and requires relatively small time steps. This is detrimental to the computational efficiency since the computational cost per step is 4 times slower than the other two methods considered in this section (PER and MPIM). Thus, a high number of steps can tip the balance in favor of these latter methods. 
\end{itemize}

\section{Numerical examples}

\subsection{Example 1. Discrete system}

\begin{figure}[ht]%
	\begin{center}
		\begin{tabular}{c}
			\includegraphics[width=15.0cm]{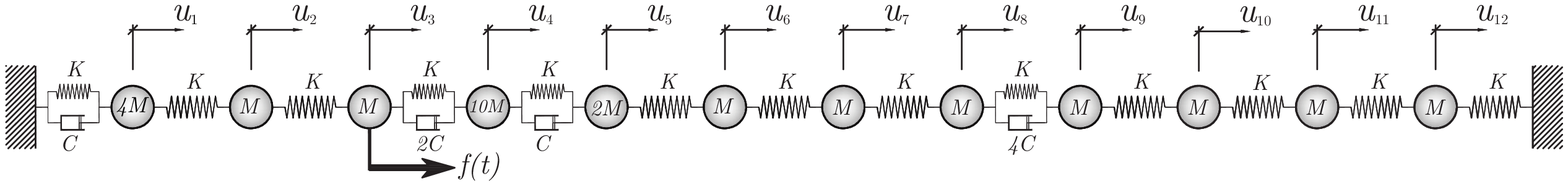} \\
			\includegraphics[width=13cm]{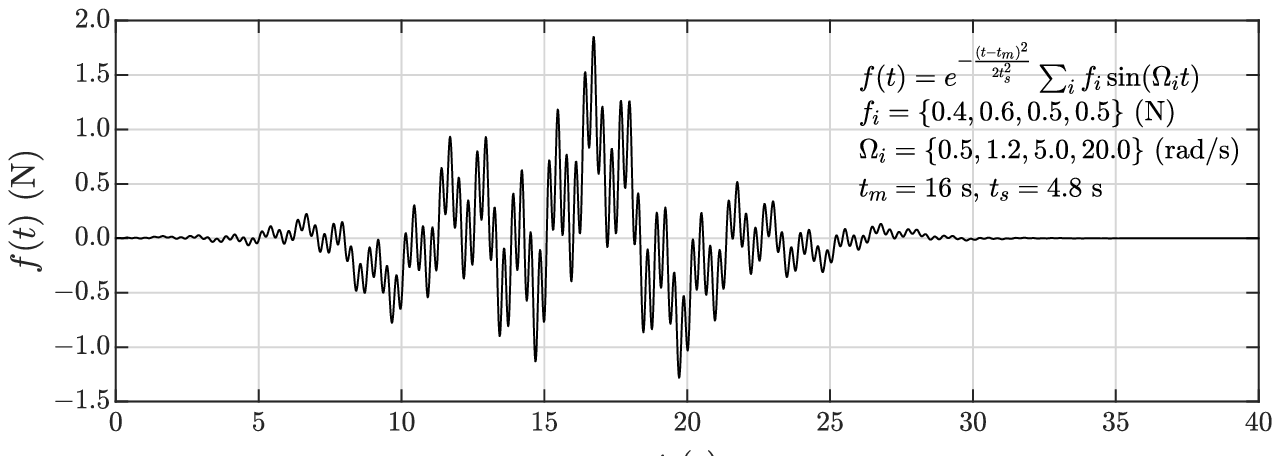} \\		
		\end{tabular}
		\caption{Discrete system of Example 1. (Top) 12-dof sketch of the lumped-mass structure. External applied force is located at dof \# 3. (Bottom) Time function of external force $f(t)$, formed by a superposition of harmonic functions weighted by a Gaussian function.}
		\label{fig04}%
	\end{center}
\end{figure}

We will first consider a numerical example consisting of a 12--dof discrete lumped spring--mass system shown in Fig.~\ref{fig04}(top).  Mass, rigidity and damping coefficient are respectively $M = 1$ kg, $K = 100$ N/m and $C = 2 \sqrt{KM} \zeta$, where $\zeta$ is a dimensionless parameter used to quantify the level of damping in the model. The value $\zeta=0$ leads to the undamped state with natural frequencies given in Table \ref{tab03}. As $\zeta$ increases, the coefficients of the damping matrix also increase linearly. Table~\ref{tab03} also lists the  relationships between the modal damping ratios and the parameter $\zeta$, i.e. $\xi_j / \zeta$, where $\xi_j = C'_{jj} / 2 \omega_j$ and $C'_{jj}$ denotes the $j$th main diagonal entree of the modal damping matrix, $\mathbf{C}' = \bm{\Phi}^T\mathbf{C}\bm{\Phi}$. The minimum period of the structure is $T = 0.32$ s. Although the damping model is non-proportional, i.e. $C'_{jk} \neq 0, \ j \neq k$, the modal damping ratios represent a good indicator to measure the dissipative level. 
\begin{table}
	\begin{center}
		\begin{tabular}{lrrrrrrrrrrrrr}
Mode, $j$			& 1	& 2 	& 3	& 4	& 5	& 6	\\ 
\hline 
Natural frequency, $\omega_j$ (rad/s) &
1.5482  &   3.7604 &    5.4508 &    6.7428 &    9.2233  &  11.0287  \\
Modal damping ratio, $C'_{jj} / 2 \omega_j \zeta $ &
0.0972  &    0.0979 &    0.5608 &    0.7106 &    0.3382 &    1.0004  \\ 
\hline \\
Mode $j$	&		7	& 8	& 9	& 10	& 11 & 12  \\
\hline 
$\omega_j$ (rad/s) &    11.8809  &   14.5269  &   16.8189  &   17.5123  &   18.5560  &   19.6346 \\
$C'_{jj} / 2 \omega_j \zeta $  &     0.8391  &   0.3872 &    1.3389  &   0.5817 &    0.1603 &    1.8695 \\
\hline
\end{tabular}
	\end{center}
	\caption{Example 1. Natural frequencies $\omega_j$ and modal damping ratios $\xi_j = C'_{jj}/2\omega_j$ (relative to the damping parameter $\zeta$) }
	\label{tab03}
\end{table}
The external force $f(t)$ is located at the 3rd dof and it reproduces a multiharmonic input affected by a Gaussian weight. The mathematical model and the corresponding graph are shown in Fig.~\ref{fig04}(bottom) and it is formed by 4 harmonic functions with different amplitudes. The proposed method will be compared with some existing iterative numerical methods (listed below). The time step $\dt$ is a common parameter for all the methods.

\begin{description}
	\item[Runge--Kutta (RK4).] The fourth order Runge-Kutta method (RK4, \cite{Butcher-2008,Burden-2001}), already presented in Sec. \ref{RK4_method},  produces a truncation error at each time-step of the order $\mathcal{O}(\dt^5)$ and an overall error of order $\mathcal{O}(\dt^4)$. The main computational cost is caused by the evaluation of the vector field $\bm{\mathcal{F}}(\mathbf{U},t) = \mathbf{W} \, \mathbf{U} + \mathbf{h}(t)$ four times in each iteration. 
	\item[Newmark method (NEW).] The Newmark integration scheme \cite{Newmark-1959,Bathe-2014}, also known as Newmark--$\beta$ method, is based on the application of the extended mean value theorem in the estimation of displacements and velocities, $\mathbf{u}_{k+1}, \dot{\mathbf{u}}_{k+1}$, from accelerations. For that, two parameters $\gamma$ and $\beta$ are used in the method,  which will be taken in our particular examples as $\gamma = 1/2$ and $\beta = 1/4$, guaranteeing the stability.  
	\item[Wilson method (WIL).] The Wilson-$\theta$ method \cite{Wilson-1973,Bathe-2014}	 assumes a linear change of acceleration within the time range $[t,t + \theta \dt]$. The method is implemented in the current numerical examples with $\theta = 1.40$ which ensures the stability of the system \cite{Wilson-1973}. Both Newmark's and Wilson's method are examples of implicit methods, which requires updating the system matrices in each iteration. 
	\item[Bathe method (BATHE).] The Bathe's method~\cite{Bathe-2005,Bathe-2007,Bathe-2012b} is an implicit scheme which consists of splitting up the time-step into two sub-steps of length $\gamma \dt $ and $(1-\gamma)\dt$ respectively where $\gamma$ is a parameter. Two different approaches are applied within each sub-step: the trapezoidal rule is used for the first sub-step and the 3-point Euler backward method for the second one. For the current numerical examples, the case $\gamma=1/2$ will be consider, which makes both time sub-steps of equal length, exhibiting in addition unconditional stability. 
	\item[Modified Precise Integration Method (MPIM).] \cite{Zhong-1994,Zhong-2004,WangMF-2005} The Modified Precise Integration Method combines two techniques: on one hand, it uses the original approach of Zhong~\cite{Zhong-1994} for calculating the exponential matrix in the homogeneous solution based on the $2^p$ algorithm. On the other hand, the integral of the  non-homogeneous term are accomplished using the Gauss quadrature  method with $g$ integration points~\cite{WangMF-2005}. In reference~\cite{WangMF-2005}, simulations are carried out for a number of integration points equal to $g=\{3,4,5\}$. So, in the current paper the same number will be used for comparison. 
	\item[Proposed method.] As described in the theoretical developments, several parameters can be tunned in the proposed method. Along the numerical examples the main matrix $\mathbf{a}(\dt)$ is computed using $m_a = r_a = 2$ since highly accurate results are achieved for these values, guaranteeing stability for a very wide range of time steps. Parameters $m_b$ and $r_b$ need to be fixed in order to obtain the non-homogeneous term $\mathbf{b}_k(\dt)$. As will be seen below, the appropriate values of these parameters depend on the level of damping and the choice of time step. Therefore, they will be left free for the time being. 
\end{description}
\begin{figure}[H]%
	\begin{center}
		\includegraphics[width=13.5cm]{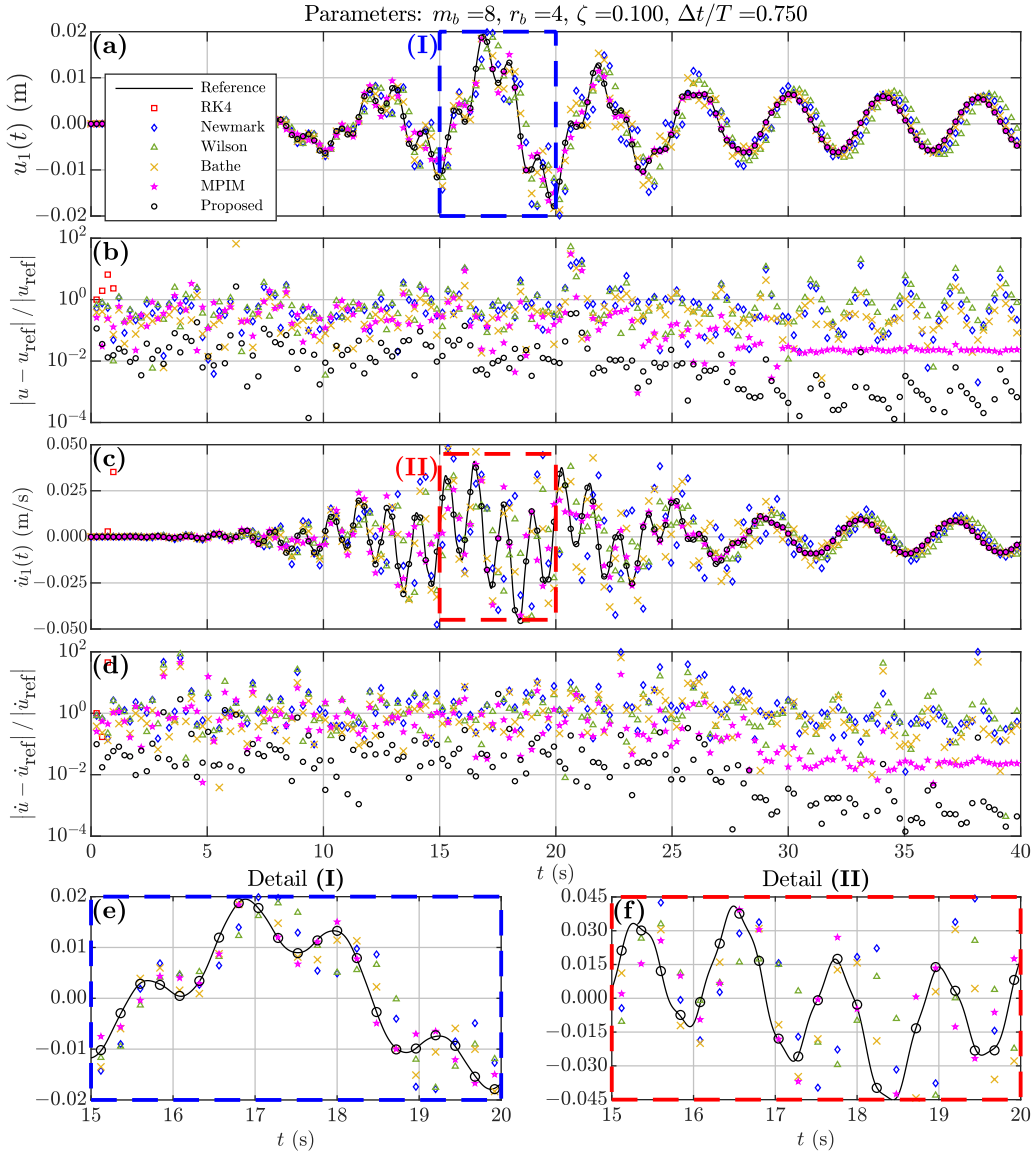} \\
		\caption{(Example 1) Time-domain numerical solutions for the different methods at dof $u_1(t)$. Dimensionless damping parameter $\zeta = 0.10$. Time-step $\dt = 0.24$ s = $0.75T$. Newmark--$\beta$ method with parameters $\gamma = 1/2, \ \beta = 1/4$. Wilson--$\theta$ method is calculated for $\theta = 1.4$. Bathe's method is implemented with $\gamma=1/2$.  Modified PIM (MPIM) uses $g=4$ Gauss quadrature points in integration of nonhomogeneous term. Order of truncation in proposed method in the non-homogeneous term: $m_b = 8$ and $r_b=4$. The reference solution has been evaluated using RK4 with a time-step $\dt/500$. (a) Displacement $u_1(t)$, (b) Local relative errors of displacements. (c) Velocity of dof--1, $\dot{u}_1(t)$, (d) Local relative error of velocity, (e) and (f) Detail plots of displacement and velocity in the range $15 \leq t \leq 20$ s. }
		\label{fig05}%
	\end{center}
\end{figure}
In fig.~\ref{fig05} the transient problem for the current example has been solved using the six different approaches presented above. Figs.~\ref{fig05}(a) and~\ref{fig05}(c) show respectively the simulated response of the first degree of freedom, i.e. displacement $u_1(t)$ and the velocity $\dot{u}_1(t)$. The interval of simulation is $0 \leq t \leq t_{\max{}}=40$ s, with a time-step of $\dt = 0.24 = 0.75T$ s. Damping parameter has been taken $\zeta = 0.10$, something that leads to a set of modal damping ratios in the range $0.97 \% \leq \xi_j \leq 18.69 \%$.
The relationship between the spectral radius of matrices $\mathbf{M}^{-1} \mathbf{C} $ and $\sqrt{\mathbf{M}^{-1} \mathbf{K} }$, 
$
\rho (\mathbf{M}^{-1} \mathbf{C} )  / \rho (\sqrt{\mathbf{M}^{-1} \mathbf{K} }) = 0.815
$
can be considered as a global measure of the dissipative forces. In order to calculate the relative errors of each numerical approach, a reference solution has been obtained using $\dt/500 = 4.8 \times 10^{-4}$ s. The relative errors along the time range have been plotted in figs.~\ref{fig05}(b) and \ref{fig05}(d). Since the time step chosen is relatively large, in fact $\dt/T \approx 0.75$, the methods based on the discretization of the derivatives (Newmark, Wilson and Bathe) are those ones that carry more numerical errors. The value of $\dt$ used in fig.~\ref{fig05} lies out of the stability region of the RK4 method. This can be observed at the first instants of the simulation, where the marker points of the RK4 (red squares) diverge. The MPIM and the proposed method have the best agreement with the reference solution. The parameter $m_b=8$ has been used for this simulation, leading to a spectral radius $\rho(\bm{\beta}_b)=0.8617<1$, guaranteeing the convergence although with a value very close to 1. We take $r_b=4$ terms  in the Neumann series noticing that from this value the errors do not become lower. The number of Gauss integration points for the MPIM are $g=4$, revealing also that for $g > 4$ the accuracy barely change. Under this choice of the parameters, better results are perceived with the proposed method, reaching relative errors between one and two orders of magnitude lower, both in displacements and velocities. In figs.~\ref{fig05}(e) and~\ref{fig05}(f), detail windows of the solution within $ 15\leq t \leq 20$ s are shown.  In them, a very satisfactory agreement of the proposed approach is perceived compared with the rest of the simulated methods. The differences are especially evident in the evaluation of velocities, where for this relatively high time step, the other approaches show errors significantly higher than those recorded with the proposed method. \\

\begin{figure}[ht]%
	\begin{center}
		\includegraphics[width=16cm]{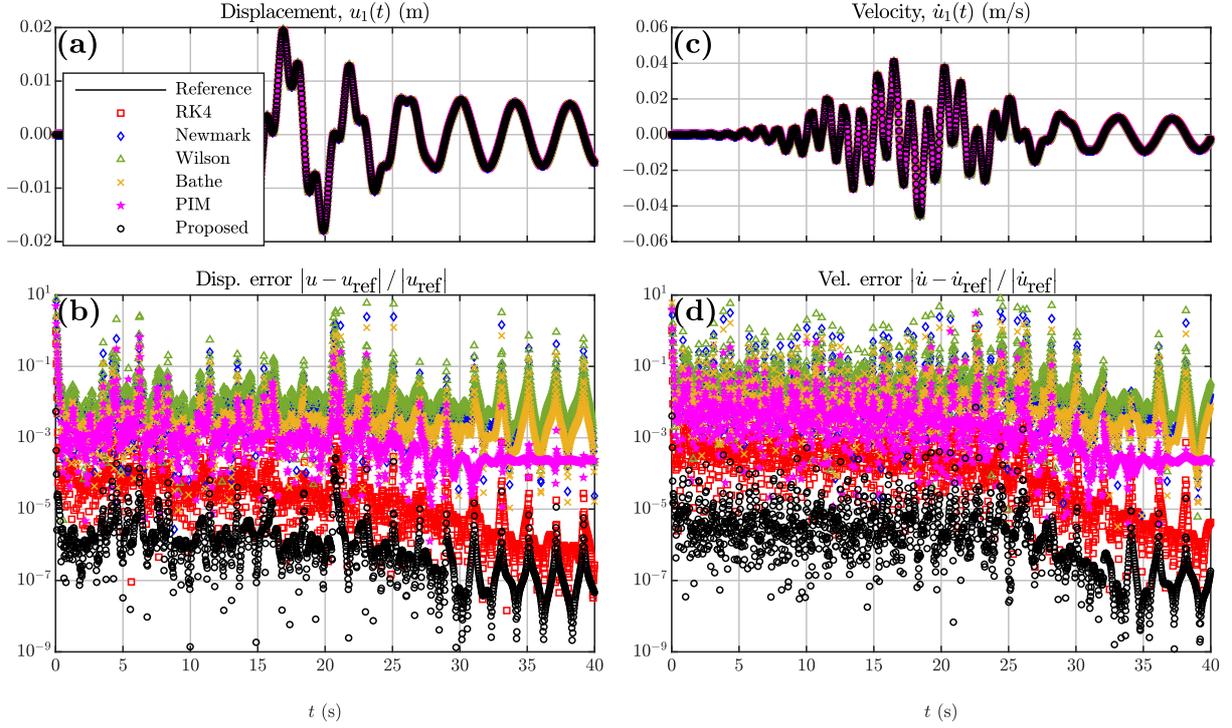} \\
		\caption{Example 1. Simulation of time-domain response of dof--1. Dimensionless damping parameter $\zeta = 0.10$. Time-step $\dt = 0.024 = 0.075T$ s. Newmark--$\beta$ method with parameters $\gamma = 1/2, \ \beta = 1/4$. Wilson--$\theta$ method is calculated for $\theta = 1.4$.  Bathe's method is implemented with $\gamma=1/2$. Modified PIM (MPIM) uses $g=4$ Gauss quadrature points in integration of nonhomogeneous term. Order of truncation in proposed method in the non-homogeneous term: $r_b=4$, $m_b = 4$. (a) Displacement, $u_1(t)$. (b) Relative error of displacement response. (c) Velocity, $\dot{u}_1(t)$. (d) Relative error of velocity response.}
		\label{fig03}%
	\end{center}
\end{figure}

It is expected that a reduction of the time step will induce an improvement in the numerical solution for all methods. This effect is shown in fig~\ref{fig03} where the displacement and velocity response together with the corresponding relative errors have been plotted for a 10--times smaller time step respect to that used for fig.~\ref{fig05} simulations, i.e. $\dt = 0.024 = 0.075T$. The damping level is left as the same value as before, keeping the value of the dissipative parameter $\zeta = 0.10$. With a time step 10 times smaller, and no change in damping, a reduction in the truncation error emerges as expected. Therefore, in order to computationally optimize the proposed algorithm, a smaller value of the truncation parameter $m_b$ can be taken. In the results of fig~\ref{fig03} the values $m_b=4$ and $r_b=4$ are chosen. In general, the relative errors drop by 2 to 4 orders of magnitude with respect to those shown in fig.~\ref{fig05}. Moreover, the RK4 method is now stable and shows very accurate results, even better than those of the MPIM. In the proposed method, the reduction of the time step leads to a substantial improvement of the predictions. Along the entire range of time simulation, the proposed approach is visibly the one that, among all the simulated approaches, yields the highest accurate results. The error in both displacement and velocity improve appreciably in all methodologies within the range $30<t<40$ s, where the non-homogeneous terms are reduced to practically zero. \\

\begin{figure}[H]%
	\begin{center}
		\begin{tabular}{cc}
			\textbf{LIGHT DAMPING} & 			\textbf{HIGH DAMPING}  \\ \\
			\includegraphics[width=7.5cm]{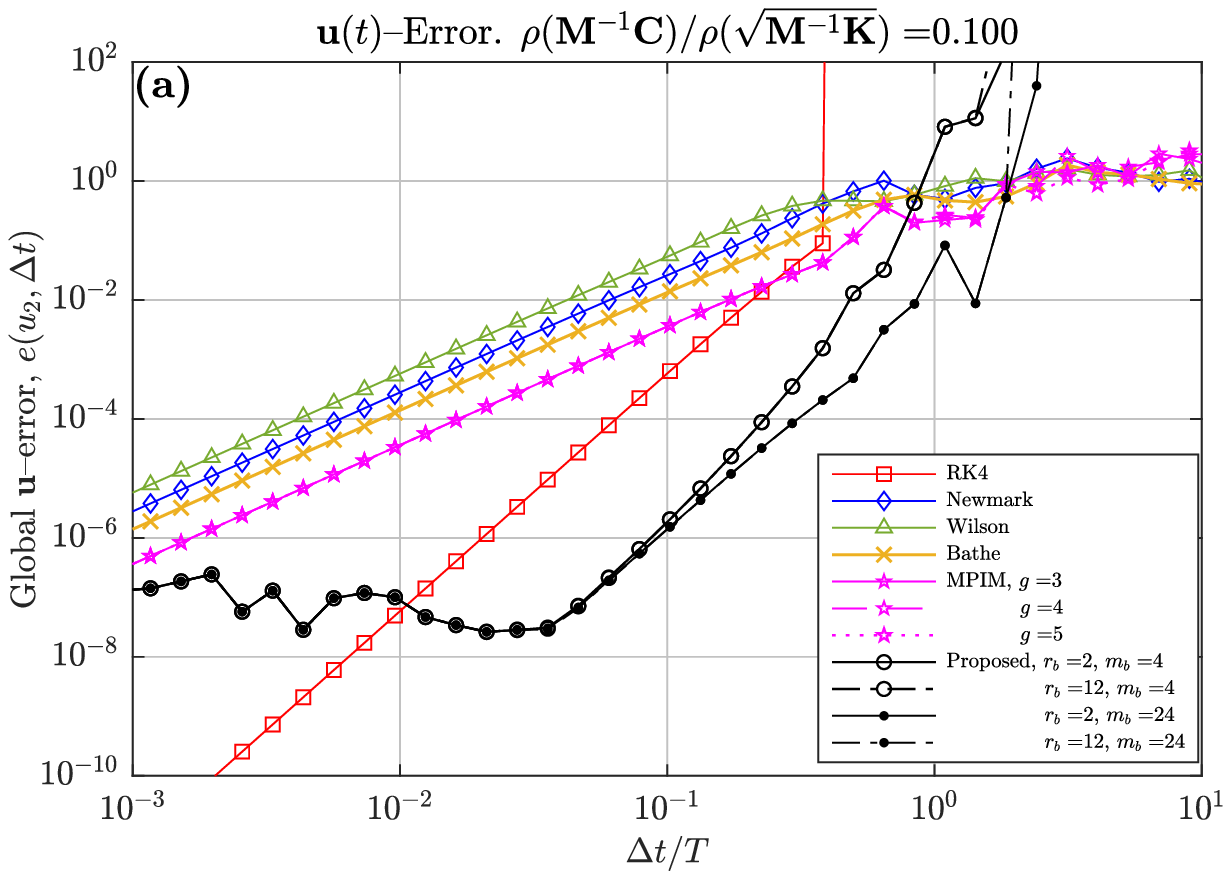} &
			\includegraphics[width=7.5cm]{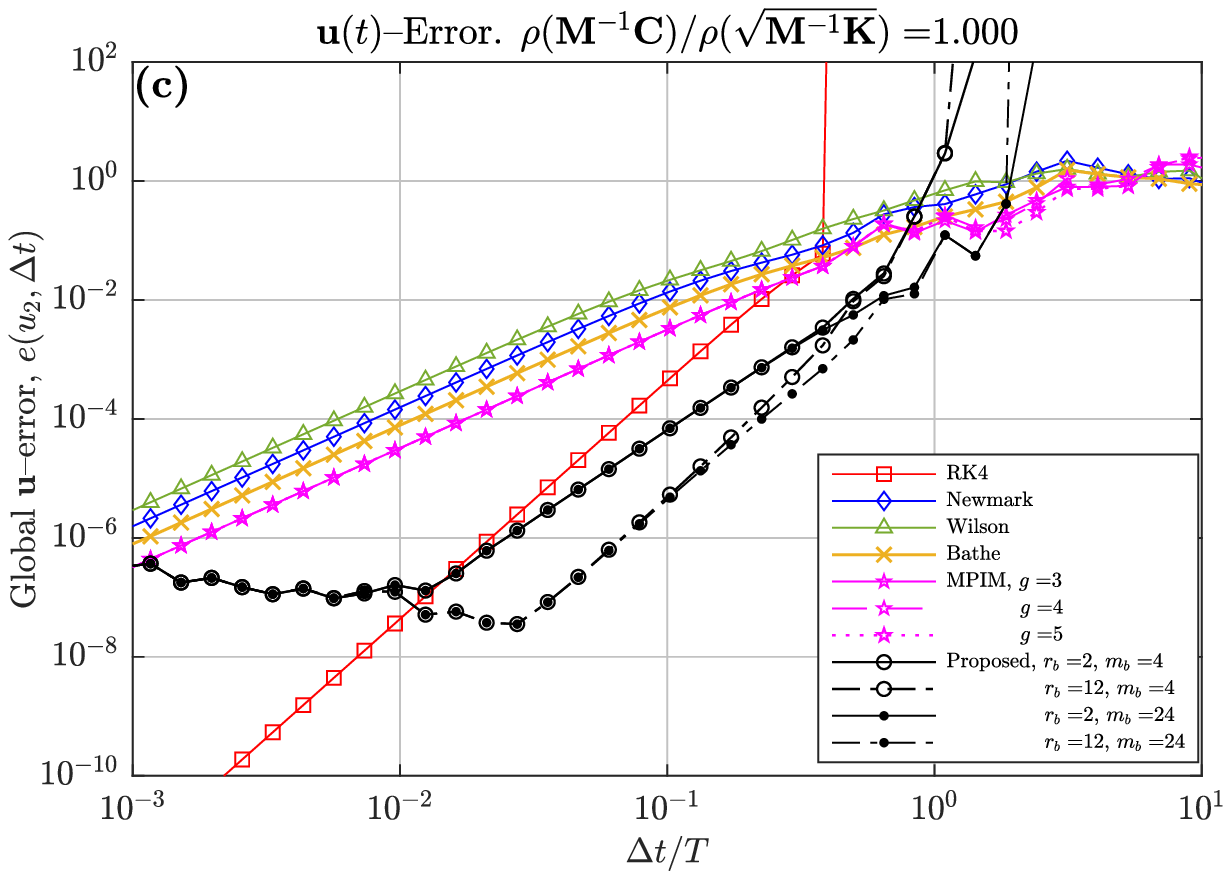} \\ 
			& \\
			\includegraphics[width=7.5cm]{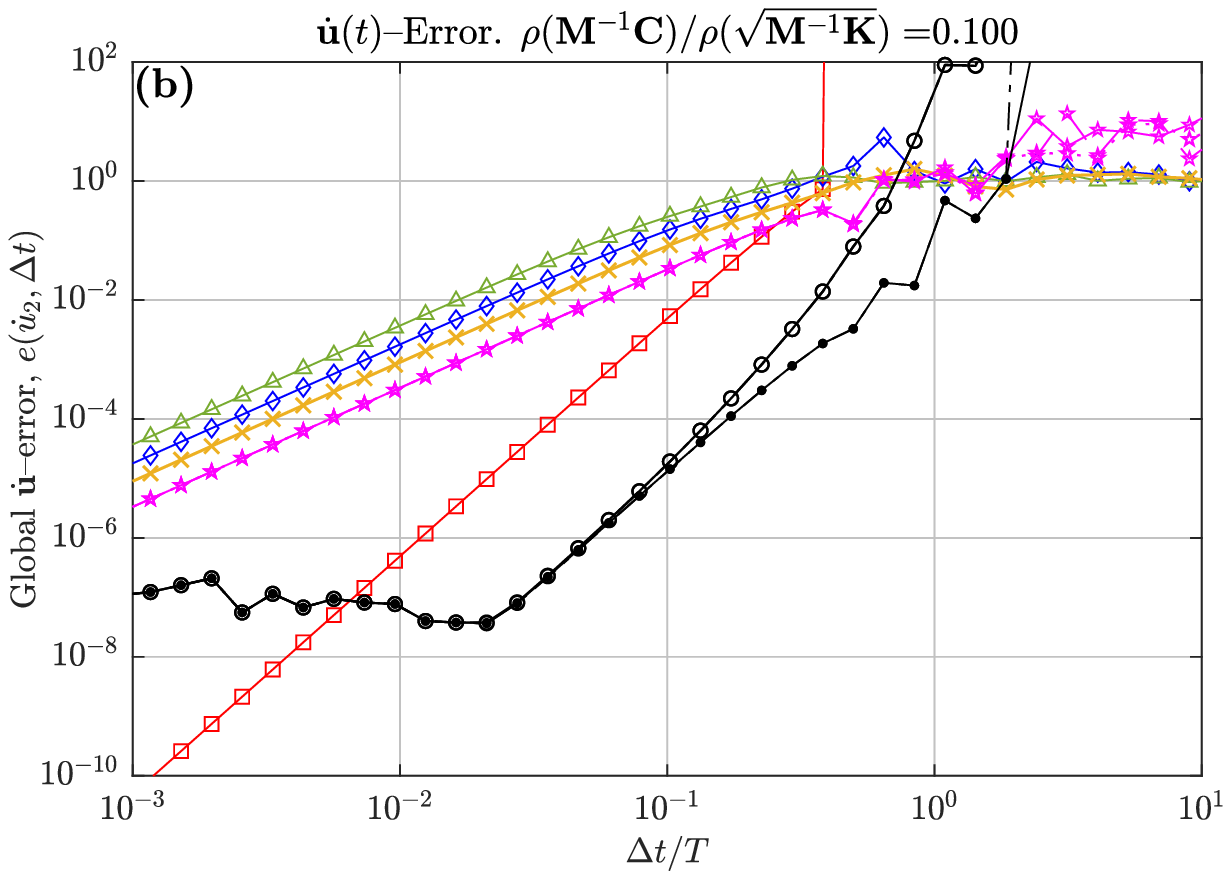} &					
			\includegraphics[width=7.5cm]{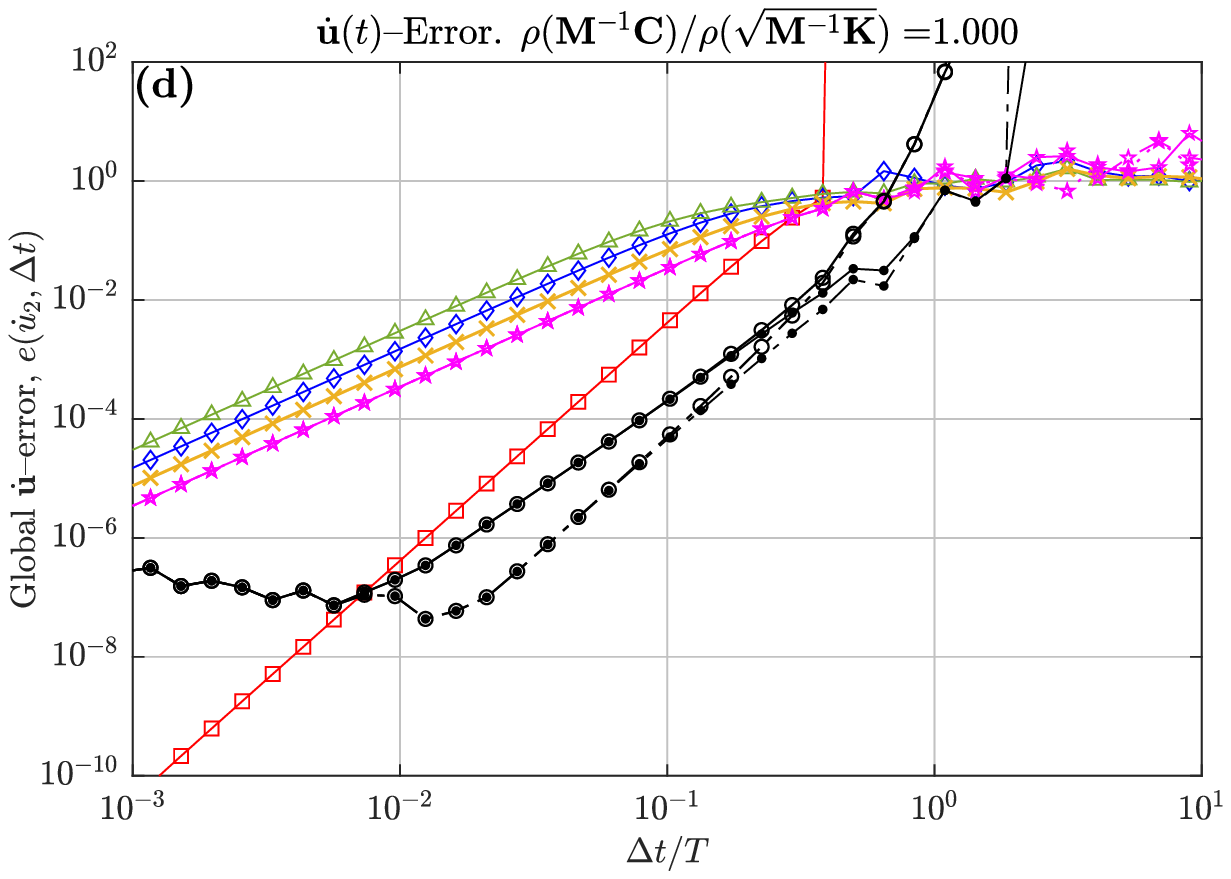} \\					
		\end{tabular}
		\caption{Relationship between global error of dof \# 2, $e(u_2,\dt)$,$e(\dot{u}_2,\dt)$ and time step $\dt / T$ (example 1), for two cases of damping level:  (a) and (b) global error plots for $\zeta = 0.012$, $\rho (\mathbf{M}^{-1} \mathbf{C} )  / \rho (\sqrt{\mathbf{M}^{-1} \mathbf{K} }) = 0.10$. (c) and (d) global error plots for $\zeta = 0.123$, $\rho (\mathbf{M}^{-1} \mathbf{C} )  / \rho (\sqrt{\mathbf{M}^{-1} \mathbf{K} }) = 1.00$.  Newmark--$\beta$ method with parameters $\gamma = 1/2, \ \beta = 1/4$. Wilson--$\theta$ method is calculated for $\theta = 1.4$.  Bathe's method is implemented with $\gamma=1/2$. Orders of truncation for the proposed method in the non-homogeneous term: $m_b \in \{4,24\}$ and $r_b \in \{2,12\}$. }
		\label{fig06}%
	\end{center}
\end{figure}

We will now focus on studying the influence of different parameters in the global error. The proposed method depends mainly on three relevant parameters, 
\begin{itemize}
	\item [(i)] The time step, $\dt$, which as usual will be measured in dimensionless form respect to the minimum period of the system. For the current example $T = 2\pi / \omega_{\max{}} = 0.32$ s.
	\item [(ii)] The damping level, quantified throughout the entries of the matrix $\mathbf{C}$, becomes of special importance in the proposed method since this latter is constructed on the basis of its perturbation. There are several forms of measuring how strong or light are the dissipative forces. In this example the coefficients of the dampers are proportional to the dimensionless parameter $\zeta$. In the results the dimensionless ratio 
	$ \rho (\mathbf{M}^{-1} \mathbf{C} )  / \rho (\sqrt{\mathbf{M}^{-1} \mathbf{K} })$ will be used as a relative measure of the damping level.
	\item [(iii)] The truncation orders,  $m_b$ and $r_b$, introduced to construct the non--homogeneous term $\mathbf{b}_k(\dt)$ in Eqs.~\eqref{eq084} and~\eqref{eq088}. 
\end{itemize}
Let us define a measure of the total error (global error) over the entire time range. Given an interval $[0,t_{\max{}}]$, a time step $\dt$ and a function $y(t)$ defined at all points $y(t_k)$, $0\leq k \leq k_{\max{}}$, $t_k = k \dt$, with $k_{\max{}} = t_{\max{}}/\dt$. Then we define the following function-norm $\left\| \bullet \right\| $ as
\begin{equation}
\left\| y(t) \right\|  = \sqrt{ \sum_{k=0}^{k_{\max{}}} y^2(t_k)    }
\label{eq114}
\end{equation}
and the global error of the $j$th dof displacement and velocity can then be defined as
\begin{equation}
e(u_j,\dt) = \frac{\left\| u_j(t) - u_{j,\textnormal{ref}}(t) \right\|}{\left\| u_{j,\textnormal{ref}} (t)\right\|}
\ , \quad
e(\dot{u}_j,\dt) = \frac{\left\| \dot{u}_j(t) - \dot{u}_{j,\textnormal{ref}}(t) \right\|}{\left\| \dot{u}_{j,\textnormal{ref}} (t)\right\|}
\label{eq115}
\end{equation}
The norm defined above will be employed as global measure of the error for each simulation and for each numerical method. The relationship between the global error and the time step is shown in figs.~\ref{fig06}(a) and~\ref{fig06}(b) show respectively the errors $e(u_2,\dt)$ (displacement) and $e(\dot{u}_2,\dt)$ (velocity)  at degree of freedom \#2, for a relatively low damping, given by a $\rho (\mathbf{M}^{-1} \mathbf{C} )  / \rho (\sqrt{\mathbf{M}^{-1} \mathbf{K} }) = 0.10$, corresponding to $\zeta = 0.0123$. The left-side plots (low damping) can be compared with those of the right side (high damping), fig.~\ref{fig06}(c) and~\ref{fig06}(d), which show the same variables but for a higher level of damping, i.e. $\rho (\mathbf{M}^{-1} \mathbf{C} )  / \rho (\sqrt{\mathbf{M}^{-1} \mathbf{K} }) = 1.00$, calculated from $\zeta = 0.123$. In view of the results, the proposed method shows evidence of high accuracy compared to the other methods studied. Although it seems to reach an error of the order of $10^{-7}$ below which it fails to decrease, unlike the rest of the methods which do undergo a monotonic decreasing as $\dt / T$ approaches to zero. Only the proposed method seems to be sensitive to the changes of the damping level. The plots shown for the rest of the approaches are hardly changed. The RK4 method shows high accuracy in the simulations, specially for low time step. However, this method exhibits a limitation by stability, something that can be observed  in figs.~\ref{fig06}, where the RK4 becomes unstable approximately for $\dt/T > 0.2$. The loss of accuracy shown by the proposed method approximately at $\dt/T \approx 1$ is not due to instability but to divergence of the sum of the infinite series \eqref{eq059}. As shown in Sec.~\ref{Convergence}, the spectral radius $\rho(\beta)$ is strongly dependent on the time-step and on the damping level.  In the algorithms developed in Sec.~\ref{Computation} for the evaluation of the scheme matrices, the $\bm{\beta}$-matrix needs to be computed twice: the first one, (matrix $\bm{\beta}_a$m, for the computation of the matrix $\mathbf{a}(\dt)$ and the second one (matrix $\bm{\beta}_b$) for the evaluation of the non-homogeneous terms within $\mathbf{b}_k(\dt)$. It is indeed this latter which somewhat limits the range of validity and may diverge for relatively high values of $\dt/T$ as seen in the example. The MPIM has been simulated for different integration points $g = \{3,4,5\}$. However, it is appreciated that the global error is barely sensitive to this parameter. In fact, curves for $g=3$, $g=4$ and $g=5$ show slight differences between each other only for high time steps. \\

As expected, the error increases with the time step also for the proposed method. However, different configurations between the $r_b$ and $m_b$ parameters result in different error distributions. In order to identify the influence of these parameters, different simulations have been carried out for $r_b = \{2, 12\}$ and $m_b = \{4, 24\}$. The observed pattern in the error evolution can be partly explained with the theoretically derived order of accuracy in Table \ref{tab06}, where the error order associated to the non-homogeneous terms is
\begin{equation}
\underbrace{
	\mathcal{O}(\dt^{m_b+2})
}_{(E1)} + 
\underbrace{ 
	\left\| \mathbf{M}^{-1} \mathbf{C} \right\| ^{r_b+1} \, \mathcal{O}(\dt^{r_b+1})
}_{(E2)}
\label{eq131}
\end{equation}
For lightly damped systems, figs \ref{fig06}(a) and \ref{fig06}(b), it is verified that 
$\rho (\mathbf{M}^{-1} \mathbf{C} )  / \rho (\sqrt{\mathbf{M}^{-1} \mathbf{K} })   \ll 1$. Thus,  the term (E2) in Eq.~\eqref{eq131} is very small compared to term (E1) and the error is actually governed this latter, i.e. $\mathcal{O}(\dt^{m_b+2}) $. This explains why the benefit of taking a high value such as $m_b=24$, respect to $m_b = 4$, is only appreciated for high values of $\dt/T$. Moreover, figs. \ref{fig06}(a) and \ref{fig06}(b) are hardly insensitive to the parameter $r_b$.  However, the result of increasing the damping is a proportional increase in the spectral radius $\rho(\bm{\beta}_b)$. If the coefficient $r_b$ is not changed, then the computation of $(\mathbf{I} - \bm{\beta}_b)^{-1}$ by Neumann series is poorer and  the error will increase. This fact can be seen in figs.~\ref{fig06}(c) and~\ref{fig06}(d)  where the error for the $r_b=2$ curves increases significantly with respect to those of the lightly damped case. If we seek to avoid this loss of accuracy, the number of $r_b$ terms in the Neumann series must be increased. Thus, increasing $r_b$ from $r_b=2$ to $r_b=12$, leads to a substantial decreasing in the error. Furthermore, the slope of the curves (close related to the order of precision) is also slightly affected by increasing $r_b$. This validates the theoretical result shown in Eq.~\eqref{eq131}, which predicts a precision order dependent on the parameter $r_b$. \\

\begin{figure}[ht]%
	\begin{center}
		\begin{tabular}{cc}
			$\dt/T = \mathbf{0.05}$  & 			$\dt/T =\mathbf{ 0.50}$ \\ \\
			\includegraphics[width=8cm]{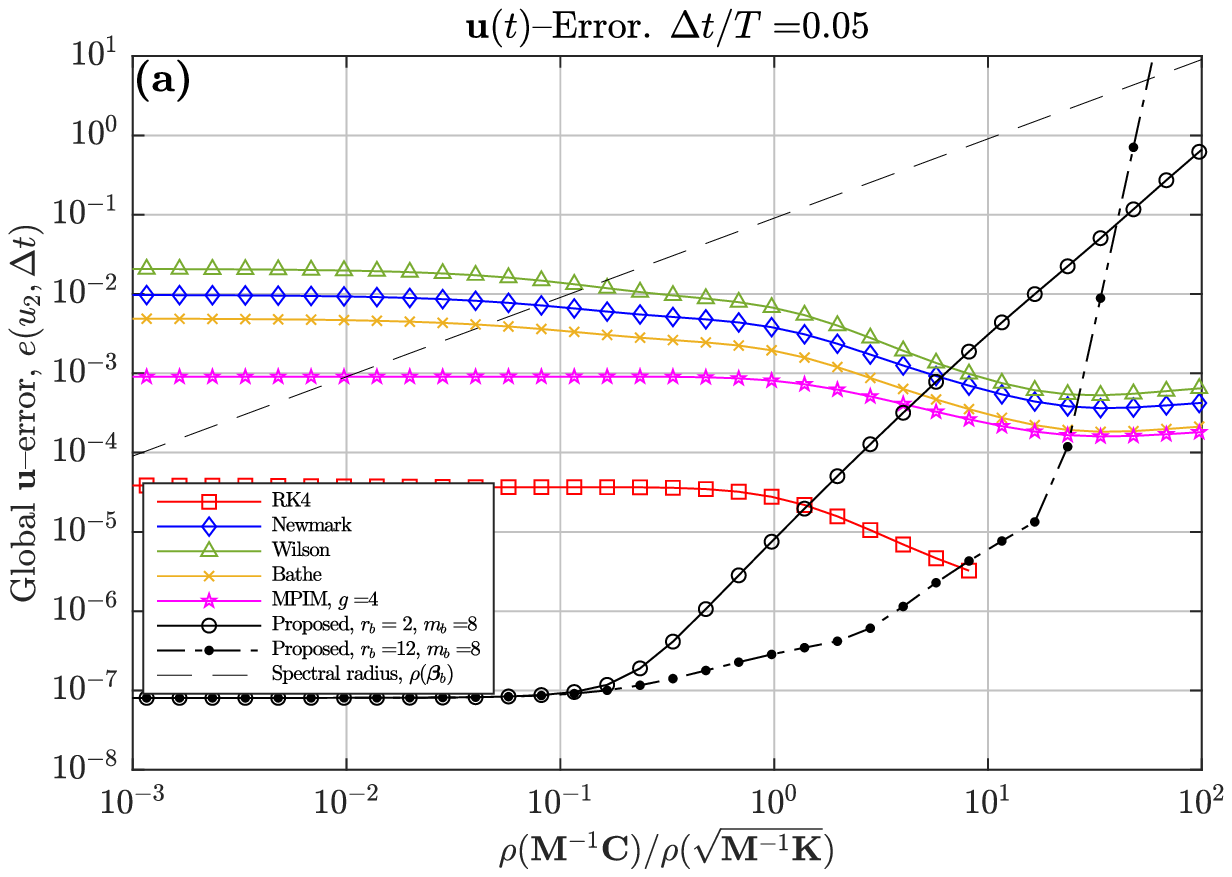} &
			\includegraphics[width=8cm]{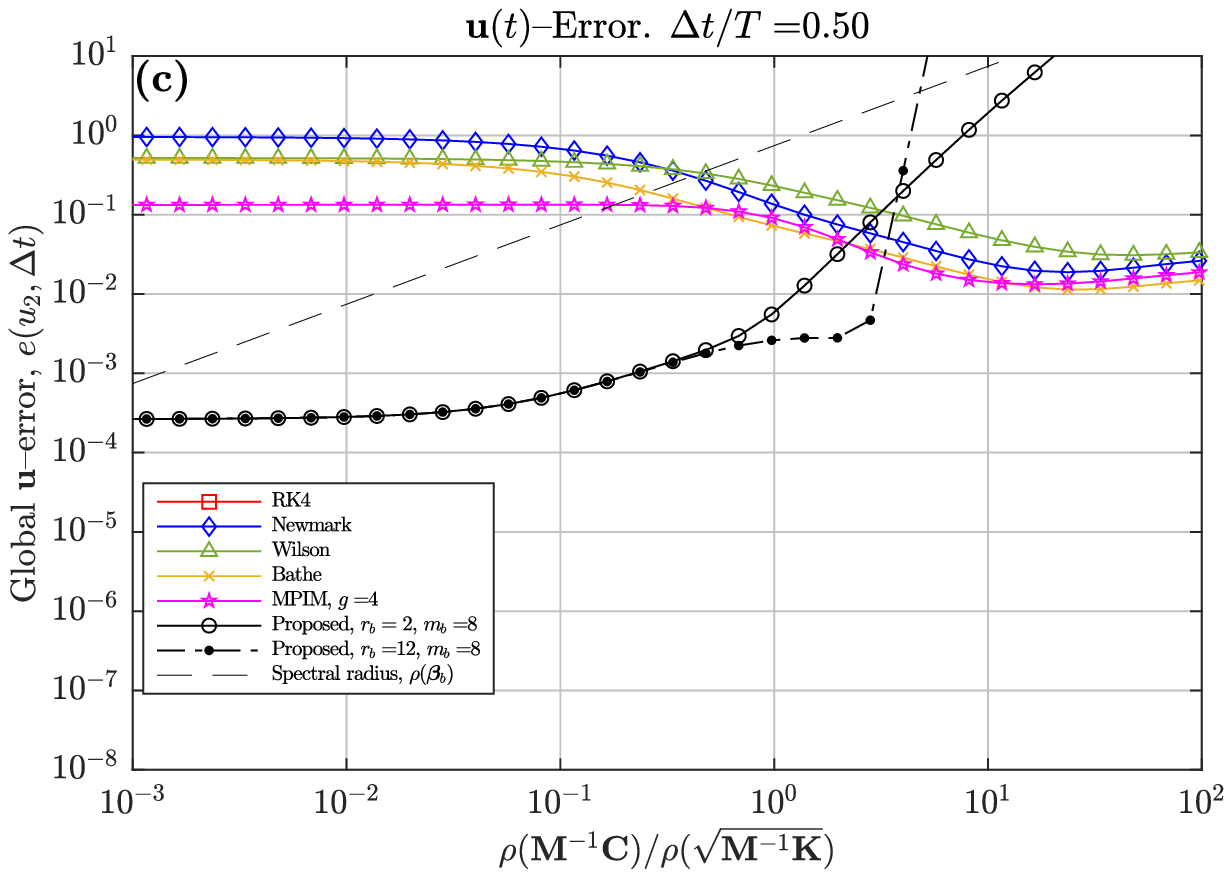} \\ 
			& \\
			\includegraphics[width=8cm]{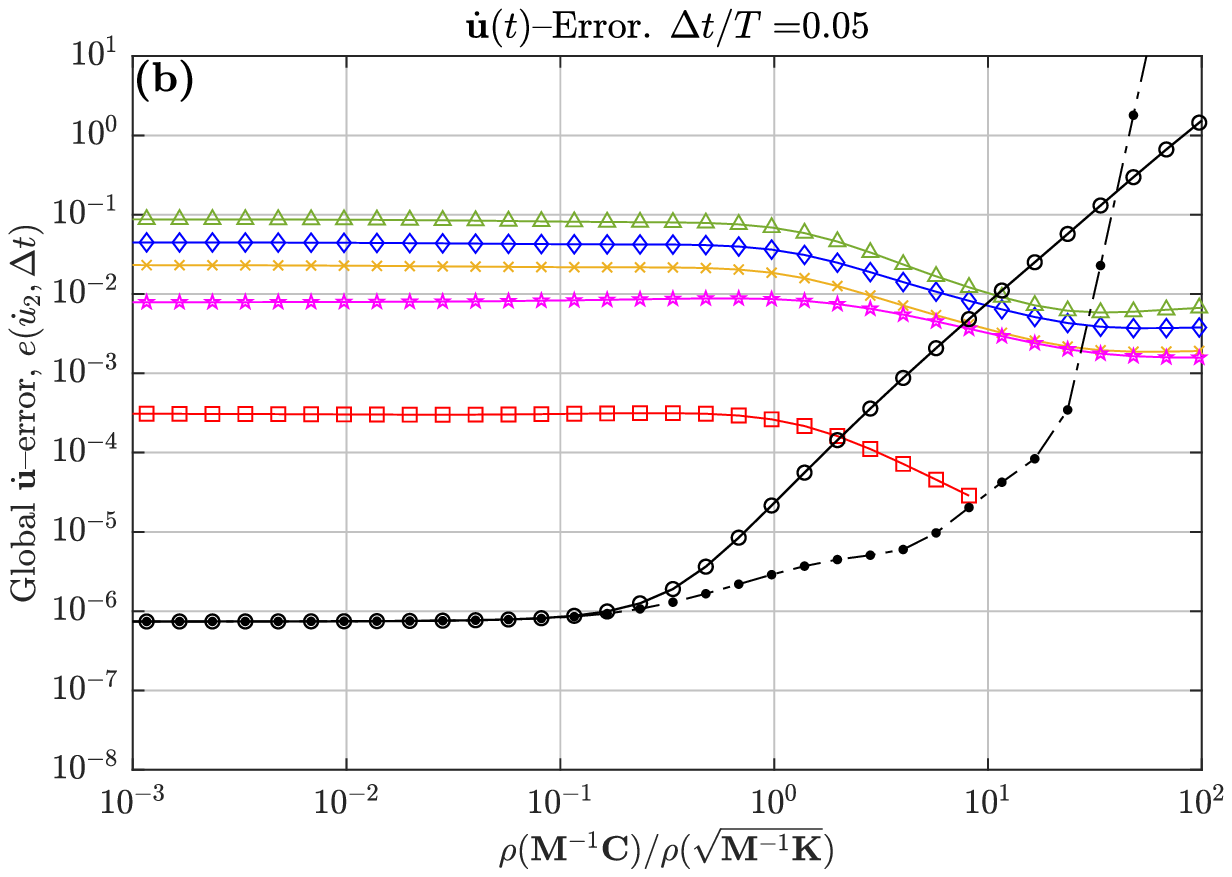} &					
			\includegraphics[width=8cm]{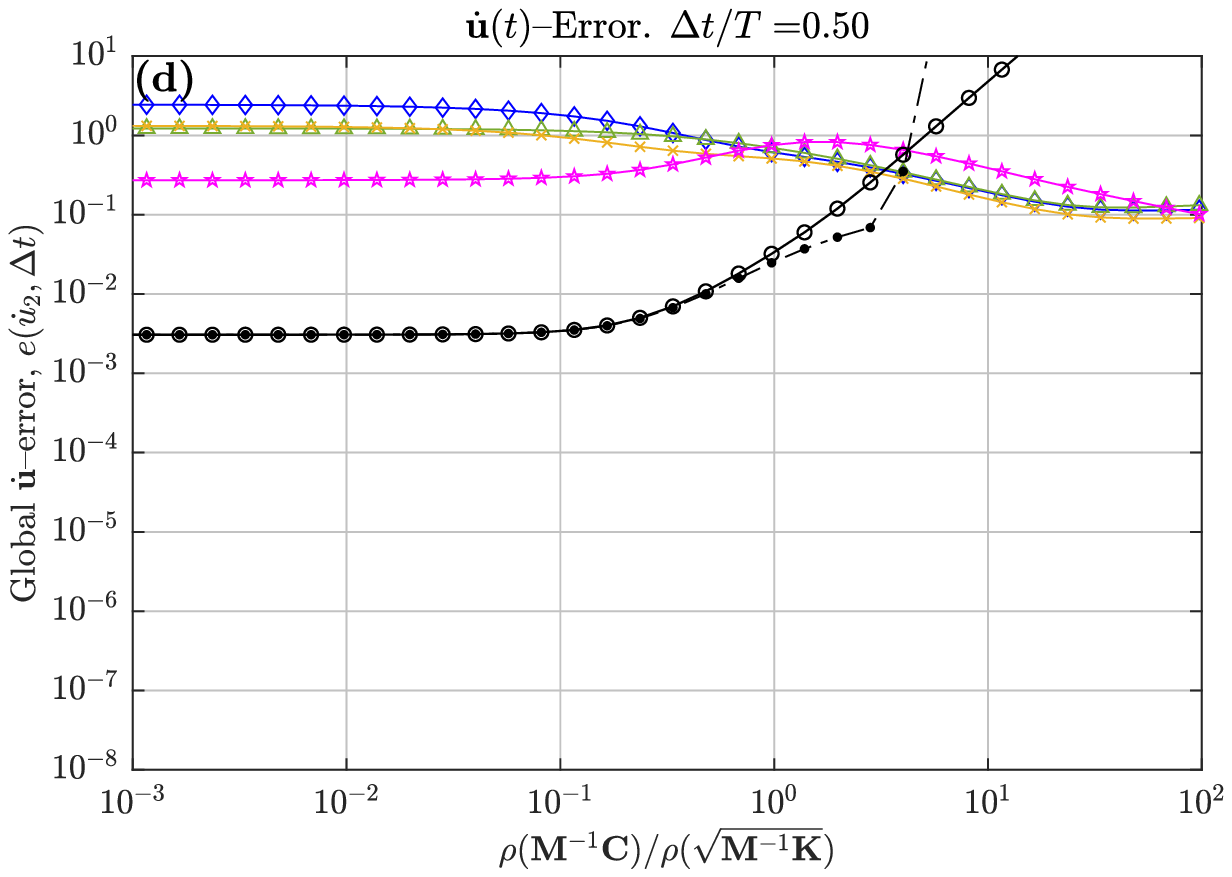} \\					
		\end{tabular}
		\caption{Relationship between global error of dof \# 2, $e(u_2,\dt)$,$e(\dot{u}_2,\dt)$ and damping level, measured as $\rho (\mathbf{M}^{-1} \mathbf{C} )  / \rho (\sqrt{\mathbf{M}^{-1} \mathbf{K} }) $ (example 1), for two time steps:  (a) and (b) global error plots for $\dt/T = 0.05$. (c) and (d) global error plots for $\dt / T = 0.5$.  Newmark--$\beta$ method with parameters $\gamma = 1/2, \ \beta = 1/4$. Wilson--$\theta$ method is calculated for $\theta = 1.4$.  Bathe's method is implemented with $\gamma=1/2$. Orders of truncation for the proposed method in the non-homogeneous term: $m_b =8$ and $r_b \in \{2,12\}$. Dashed black line  (no markers) in global error of displacements denotes the spectral radius of $\bm{\beta}_b$}
		\label{fig07}%
	\end{center}
\end{figure}
The underlying nature of the proposed method, based on the perturbation of the dissipative model, suggests that the results will depend strongly on damping. The trend observed in figs.~\ref{fig06} makes us wonder how the quality of our proposal will evolve as we increase the damping level. To this end, fig.~\ref{fig07} shows the dependency of the global error with the damping intensity, this latter measured by the magnitude $\rho (\mathbf{M}^{-1} \mathbf{C} )  / \rho (\sqrt{\mathbf{M}^{-1} \mathbf{K} })$. Left and right side plots of fig.~\ref{fig07} represent  the relative global error for two time steps, $\dt/T = \{0.05, 0.50\}$, respectively. In all performed simulations, the proposed method exhibits the lowest error and consequently the most accurate results for lightly damped systems. It is also observed that the choice of $r_b$ does not affect the error when the damping is small, something that, as above, can be explained by the fact that the error term (E2) in Eq. \eqref{eq131} is negligible with respect to (E1). In Figs \ref{fig07}(a) and \ref{fig07}(c) the spectral radius $\rho(\bm{\beta}_b)$ has also been plotted, showing its linear dependency with damping. Thus, in order to keep the errors bounded, the Neumann series needs to be updated with more terms as the damping increases. From the limit $\rho(\bm{\beta}_b)=1$ it is useless to increase $r_b$ since the Neumann series is divergent. It can be seen in the four graphs of fig. \ref{fig07} the divergence of the method when this threshold is exceeded. The other methods under study are also insensitive to damping at low values of $\rho (\mathbf{M}^{-1} \mathbf{C} ) $. However, above certain limiting values, damping improves the accuracy of the result slightly, contrary to what is observed with the current approach. Despite this improvement, the damping affects the stability conditions for the RK4 method and the accuracy for the MPIM. For the time step $\dt/T = 0.05$, the RK4 method becomes unstable for $\rho (\mathbf{M}^{-1} \mathbf{C} )  / \rho (\sqrt{\mathbf{M}^{-1} \mathbf{K} })\approx 8$. On the right side, for $\dt/T = 0.5$, the RK4 is unstable in the whole range of damping level and MPIM exhibits a pronounced loss of accuracy for high damping levels. Newmark, Wilson and Bathe methods show decreasing errors as the damping level is increased. We have purposely left the same scale on the ordinate (error) axes for all graphs. This allows to observe at a simple glance the effect of multiplying the time step by 10, from $\dt/T = 0.05$ in figs.~\ref{fig07}(a) and ~\ref{fig07}(b) to $\dt/T = 0.50$ in figs.~\ref{fig07}(c) and ~\ref{fig07}(d). As expected, all methods exhibit an increasing error, although the proposed method for these particular values still yields the best numerical results for light and moderate damping problems.

\subsection{Example 2. Continuous system}

\begin{figure}[ht]%
	\begin{center}
		\includegraphics[width=16.5cm]{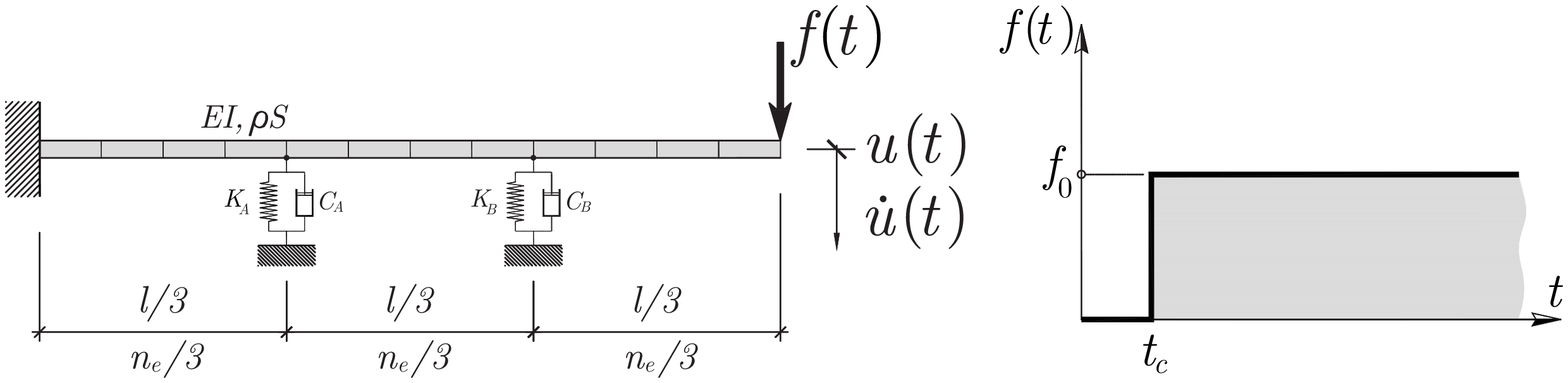} 
		\caption{Example 2. Continuous system}
		\label{fig10}%
	\end{center}
\end{figure}

\begin{figure}[H]%
	\begin{center}
		\begin{tabular}{rr}
			\includegraphics[width=7.20cm]{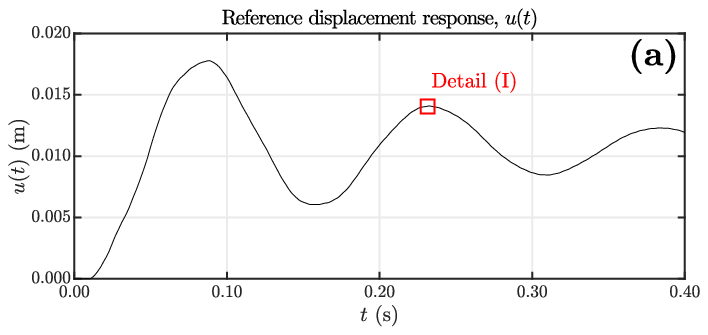} & 
			\includegraphics[width=7.20cm]{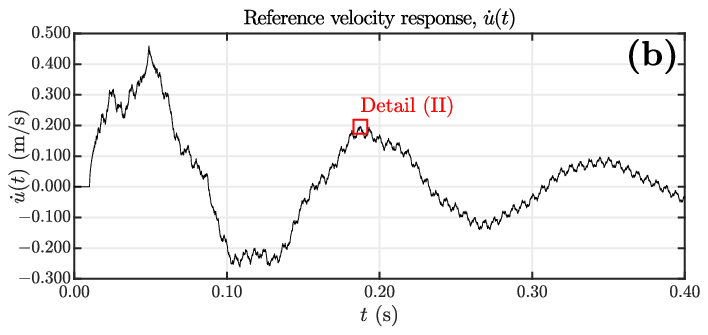} \\			
			\includegraphics[width=7.00cm]{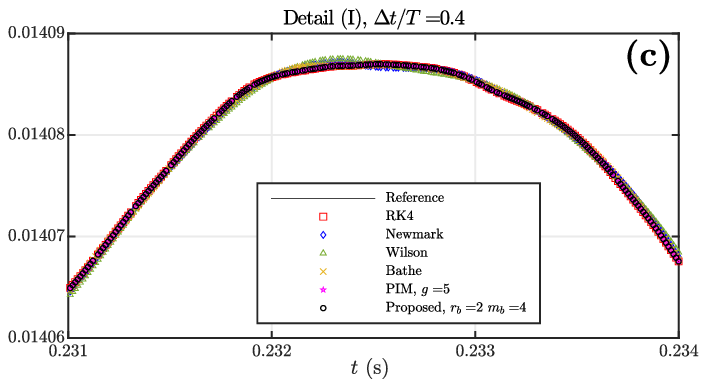} & 
			\includegraphics[width=7.00cm]{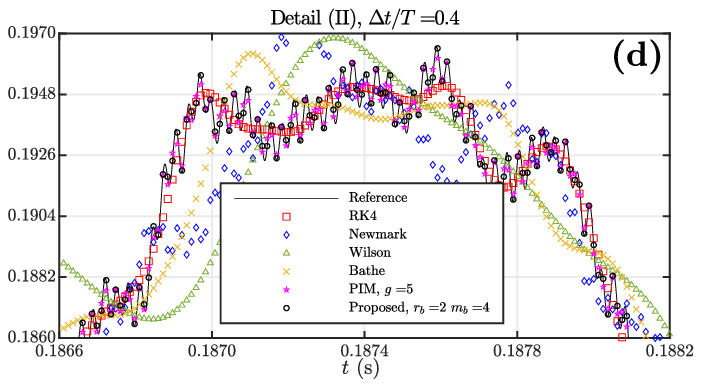} \\ 			
			\includegraphics[width=7.0cm]{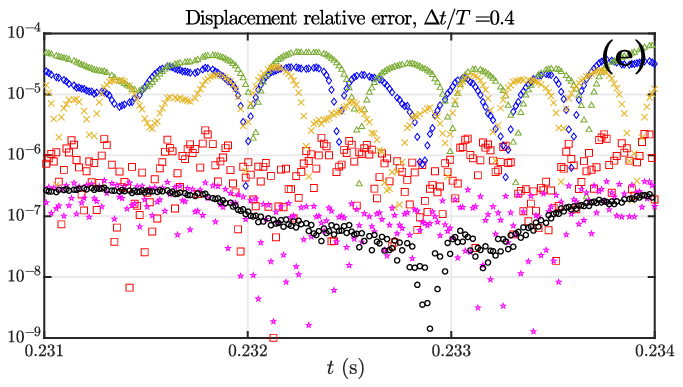} & 
			\includegraphics[width=7.0cm]{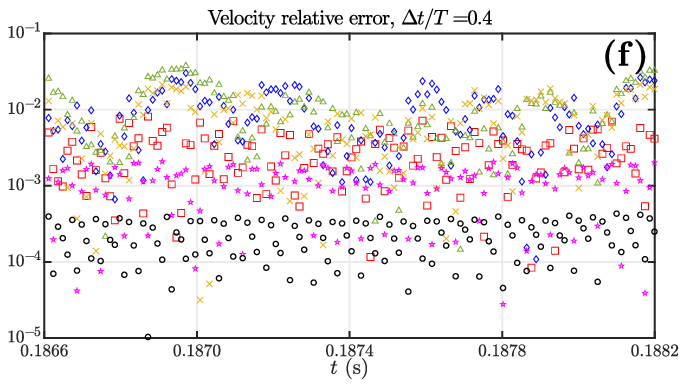} \\ \\			
			\includegraphics[width=7.00cm]{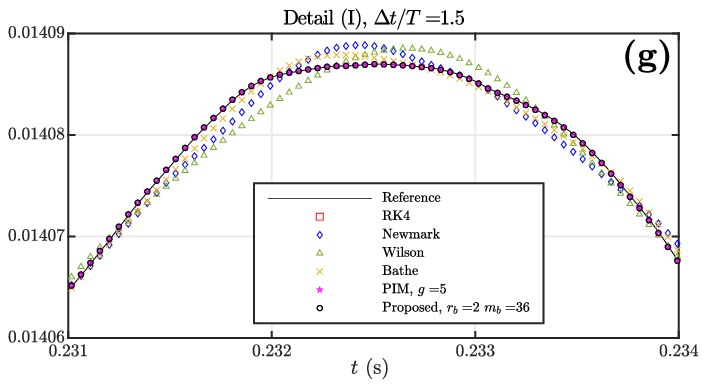} & 
			\includegraphics[width=7.00cm]{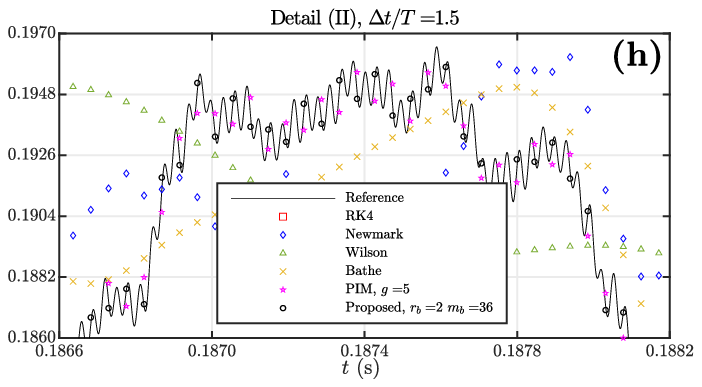} \\ 			
			\includegraphics[width=7.0cm]{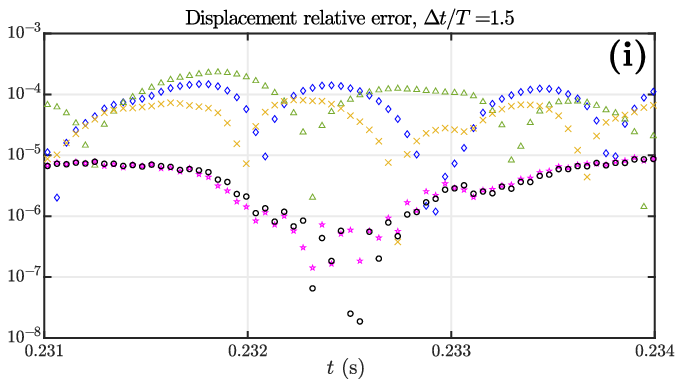} & 
			\includegraphics[width=7.0cm]{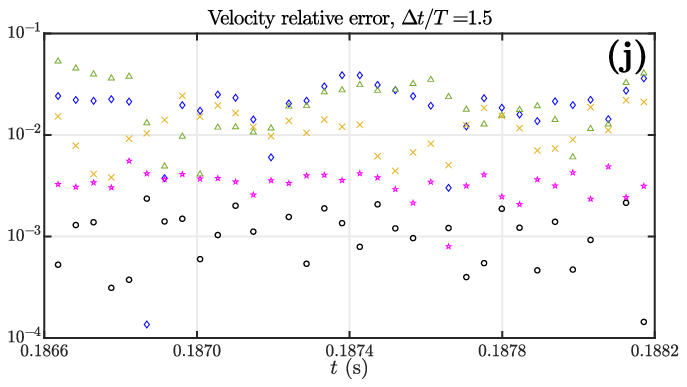} \\ 			
		\end{tabular}
		\caption{Numerical simulations of the continuous beam of Example 2 with $n_e = 24$ elements and $N = 48$ dof with a minimum period of $T=3.1013 \times 10^{-5}$ s. The damping parameter of supports is $\zeta_A = \zeta_B=0.50$. (a) and (b): Reference response, displacement $u(t)$ and velocity $\dot{u}(t)$ at the right-edge of the beam respectively. (c) and (d): Details (I) and (II) of the displacement and velocity for $\dt = 0.4T$. (e) and (f): relative error of displacement and velocity respectively, for $\dt = 0.4T$. (g) and (h): Details (I) and (II) of the displacement and velocity for $\dt = 1.5T$. (i) and (j): relative error of displacement and velocity respectively, for $\dt = 1.5T$}
		\label{fig08}%
	\end{center}
\end{figure}

In this example, the application of the proposed numerical method to continuous dynamic problems is explored. A  cantilever beam of fig. \ref{fig10} supported by two elastic constraints with viscous dampers is considered. The total length is $l=3$ m and the cross-section is steel-made with bending stiffness $EI = 437.5$ kNm$^2$ and total mass $m_0 = 235.5$ kg. A finite element model based on Euler-Bernoulli beam theory with $n_e$ elements is used. The stiffnesses of the supports are $K_A = 20 EI/l^3 = 324$ kN/m and $K_B = 10 EI/l^3 = 162$ kN/m . The damping coefficients are respectively $c_A = 2 m_0 \omega_r \zeta_A$ and $c_B = 2 m_0 \omega_r \zeta_B$, where $\omega_r = \sqrt{ EI/m_0 l^3}$ = 8.29 rad/s is a reference frequency. The damping parameters $\zeta_A, \ \zeta_B$ will be fixed later for each simulation. A point force $f(t) $ is applied downwards at the free right-end of the beam, as shown in fig. \ref{fig10}(left). The force vs time distribution is shown in the graph on in fig. \ref{fig10}(right), where
\begin{equation}
f(t) =
\begin{cases}
0  &  t < t_c \\
f_0  &  t \geq t_c
\end{cases}
\end{equation}
The following parameters are taken: $t_c = 0.01$ s and $f_0= 1$ kN. In fig.~\ref{fig08} some simulations of the beam have been carried out using $n_e = 24$ elements and $N=48$ degrees of freedom. The structure has a minimum period of $T = 3.1013 \times 10^{-5}$ s. The damping parameter for both elastic supports has been taken as $\zeta_A = \zeta_B=0.50$ resulting modal damping ratios in the range $0.6 \% \leq \xi_j \leq 11.2 \%$ in the six first modes, while  higher modes present  values lower than 0.4 \%. The time domain response has been computed for both displacement and velocity $u(t)$ and $\dot{u}(t)$ at the right-end of the beam, downwards, as indicated in fig. \ref{fig10}. Figs.~\ref{fig08}(a) and~\ref{fig08}(b) represent the reference solution  (displacement and velocity, respectively) calculated using a time step of $\dt/500$. Two different cases corresponding to two time steps have been considered, namely
\begin{description}
	\item [Case $\dt/T = 0.4$] [figs.~\ref{fig08}(c), \ref{fig08}(d), \ref{fig08}(e) and \ref{fig08}(f)]: The truncation orders in the non-homogeneous term of the proposed method are chosen as $r_b = 2$ and $m_b=4$ and the spectral radius of the $\bm{\beta}_b$ matrix is $\rho(\bm{\beta}_b)=0.0011$. In general, something common to all considered methods is the fact that errors in displacements are significantly smaller than those of the velocities. It is observed that the methods based on the finite difference approximation are those ones with the highest errors. The displacement error curves show similar accuracy in the calculation of $u(t)$ for the MPIM (using $g=5$ integration points) and the proposed method, see Fig. \ref{fig08}(e). However, velocities calculated with the proposed method is reported to be considerably more accurate than the others, as shown in Fig. \ref{fig08}(f).
	\item [Case $\dt/T = 1.5$] [figs.~\ref{fig08}(g),  \ref{fig08}(h), \ref{fig08}(i) and \ref{fig08}(j)]: Increasing the time step obviously leads to increase the relative error in all methods. In particular, the RK4 method becomes unstable so that its results are not displayed in the plots. Again, results observed in displacements are more accurate than those of velocities. The precision of MPIM and that of the  proposed method are very similar in the calculation of displacement, but the evaluation of velocities turns out to be more accurate with the proposed method. 
	To achieve this accuracy, it has been necessary to increase the parameter $m_b$ up to the value $m_b=36$. In general, increasing $m_b$ is the most effective when faced with a significant increase in time step. While, on the other hand, including more terms in the Neumann series (increasing $r_b$) is only efficient for cases where $\rho(\bm{\beta}_b)$ is moderately high (but not reaching one).  In case of study, since the spectral radius is $\rho(\bm{\beta}_b)=0.0033\ll 1$, it is not necessary to increase $r_b$, and the value $r_b=2$ results sufficient. 
\end{description}
%

\begin{figure}[ht]%
	\begin{center}
		\begin{tabular}{cc}
			\includegraphics[width=8cm]{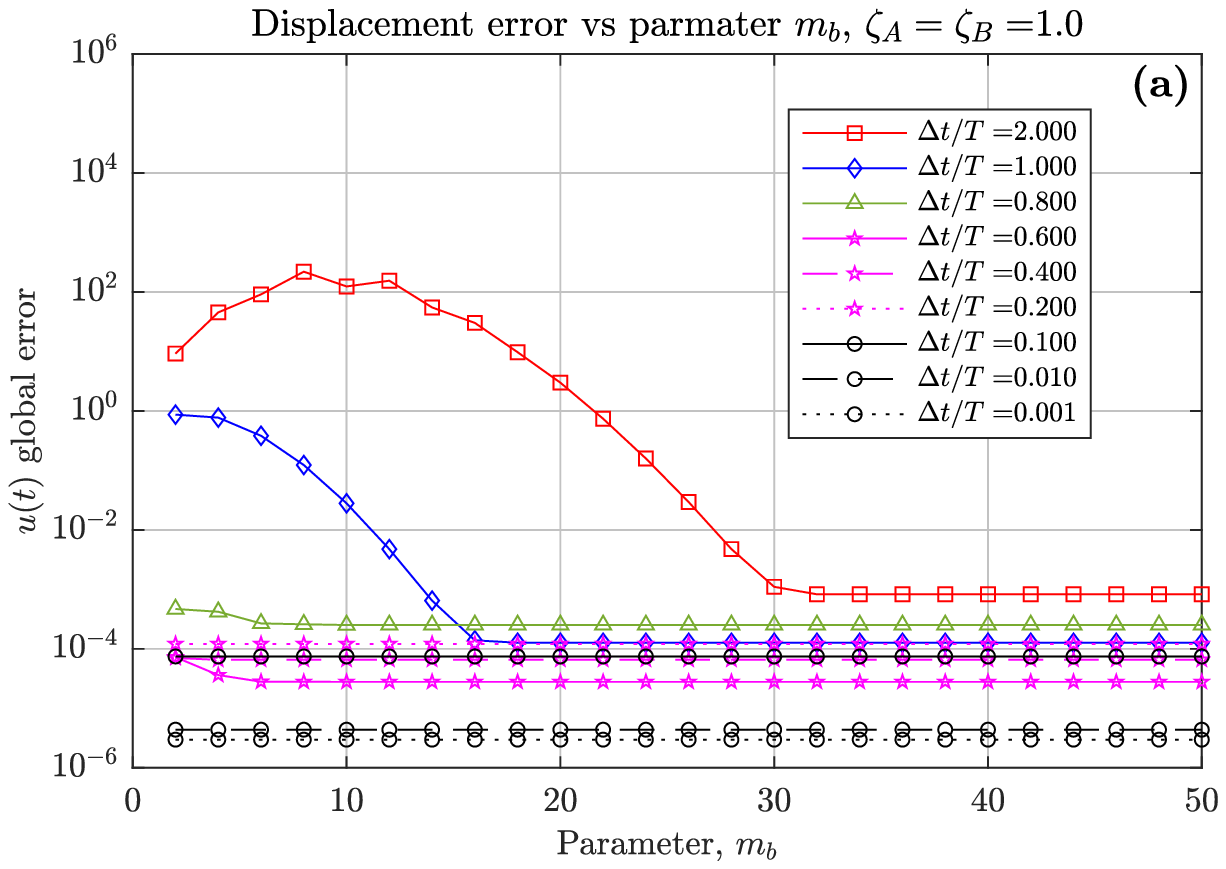} & 
			\includegraphics[width=8cm]{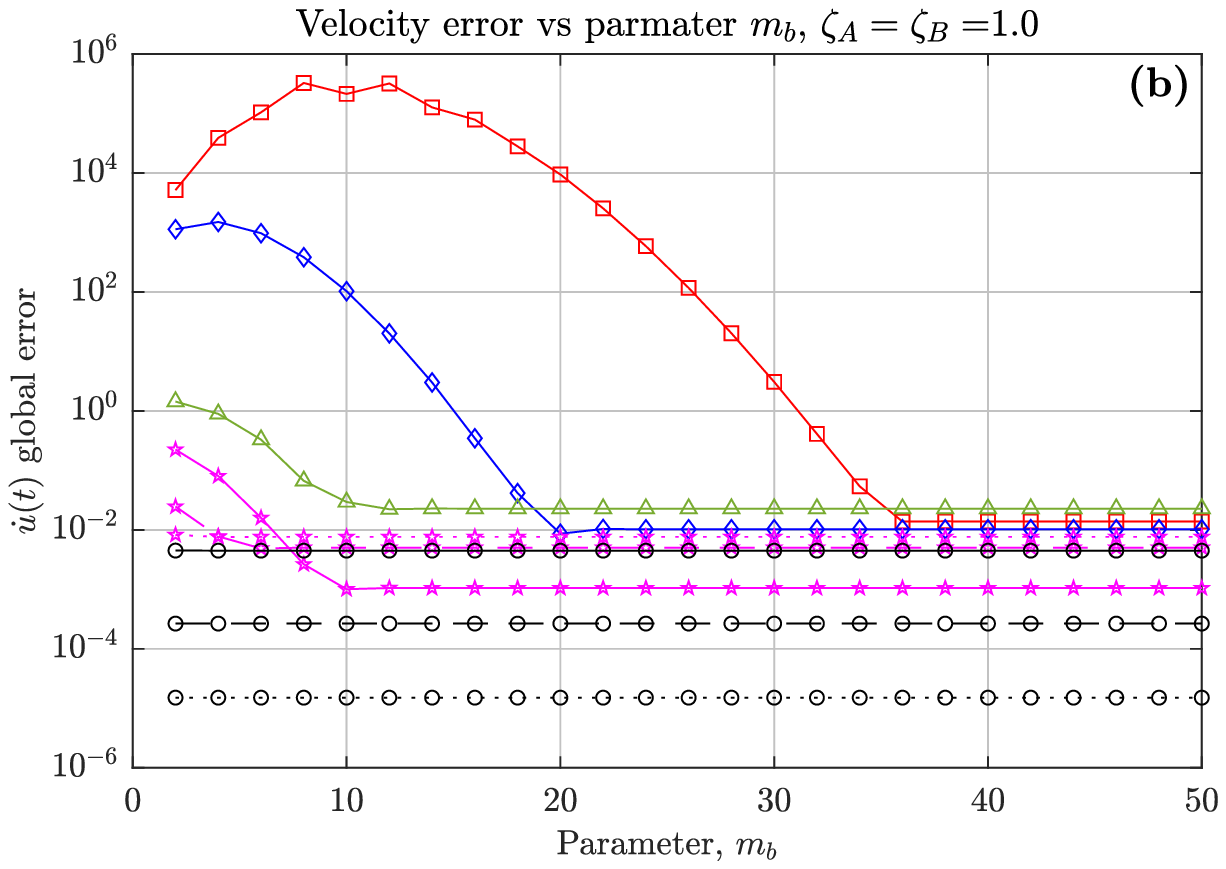} \\			
		\end{tabular}
		\caption{Influence of the order of truncation $m_b$ in the global errro. (a) Global error of the displacement $u(t)$. (b) Global error of the velocity $\dot{u}(t)$. Damping parameters of the elastic supports: $\zeta_A=\zeta_B=1$. The different curves depict the effect under several time steps from $\dt / T = 0.001$ up to $\dt / T = 2$. }
		\label{fig09}%
	\end{center}
\end{figure}

After several numerical examples simulated, the truncation order $m_b$ of the $\bm{\beta}_b$ series contributes to a significant reduction of the error for relatively high values of the time step. However, for low or moderate values, no visible effect is perceived. In order  to better understand the influence that the parameter $m_b$ has on the solution, some numerical simulations have been performed by changing this parameter for different time steps. The results in Fig. \ref{fig09} show that the errors of displacements seem to be practically insensitive to the parameter $m_b$ for $\dt / T <1$. Only from this value (in Fig. \ref{fig09}(a), the curves $\dt/T=1$ and $\dt/T=2$), it is really effective to increase $m_b$ to bound the error, although this latter stabilizes from a limit value of $m_b$. Thus, for instance, when $\dt/T=1$ the error does not decrease for $m_b \geq 16$, while in the case $\dt/T=2$, this limit is at $m_b\geq 32$. The velocity errors follows the same pattern as the displacements. Although the parameter $m_b$ starts to be sensitive already for values for $\dt/T=0.5$ onwards. Thus, for example, the error in the velocities for $\dt/T=2$ stabilizes for $m_b \geq 36$. As can be seen in Fig. \ref{fig09}, although the error stabilizes, the stabilization value grows considerably with increasing time step, especially in the computation of the velocities. \\

The behavior observed in Fig. \ref{fig09} makes useless and unnecessary to take very high values of $m_b$. This is advantageous since instabilities have been identified in the evaluation of the series for very high values of $m_b$, in the current example for the order $m_b>60$. Although they are under study, it seems likely that they are due to the numerical computation of terms of the type $\mathbf{A}^j / (2j+4)! $ with $j$ very high. \\

\begin{figure}[ht]%
	\begin{center}
		\begin{tabular}{cc}
			\includegraphics[width=8cm]{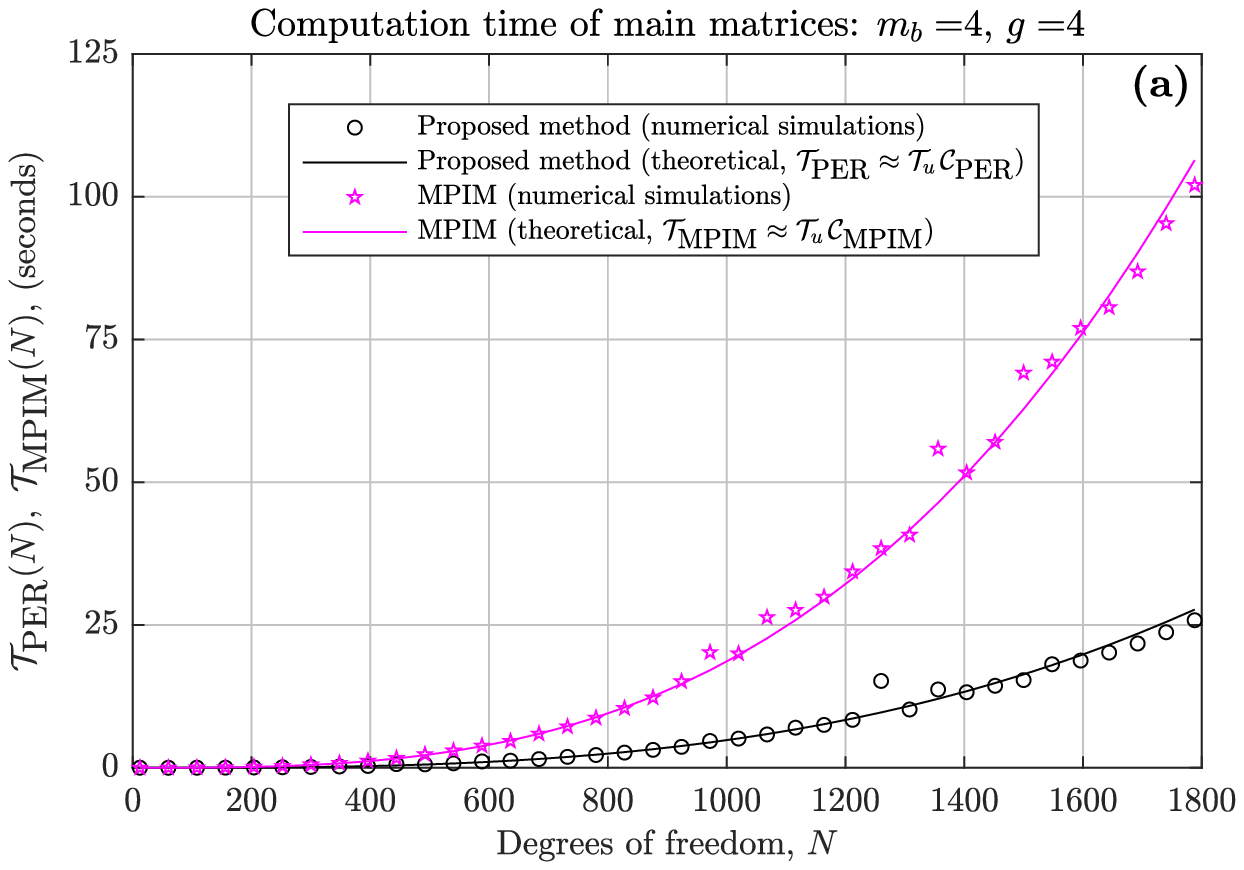} &
			\includegraphics[width=8cm]{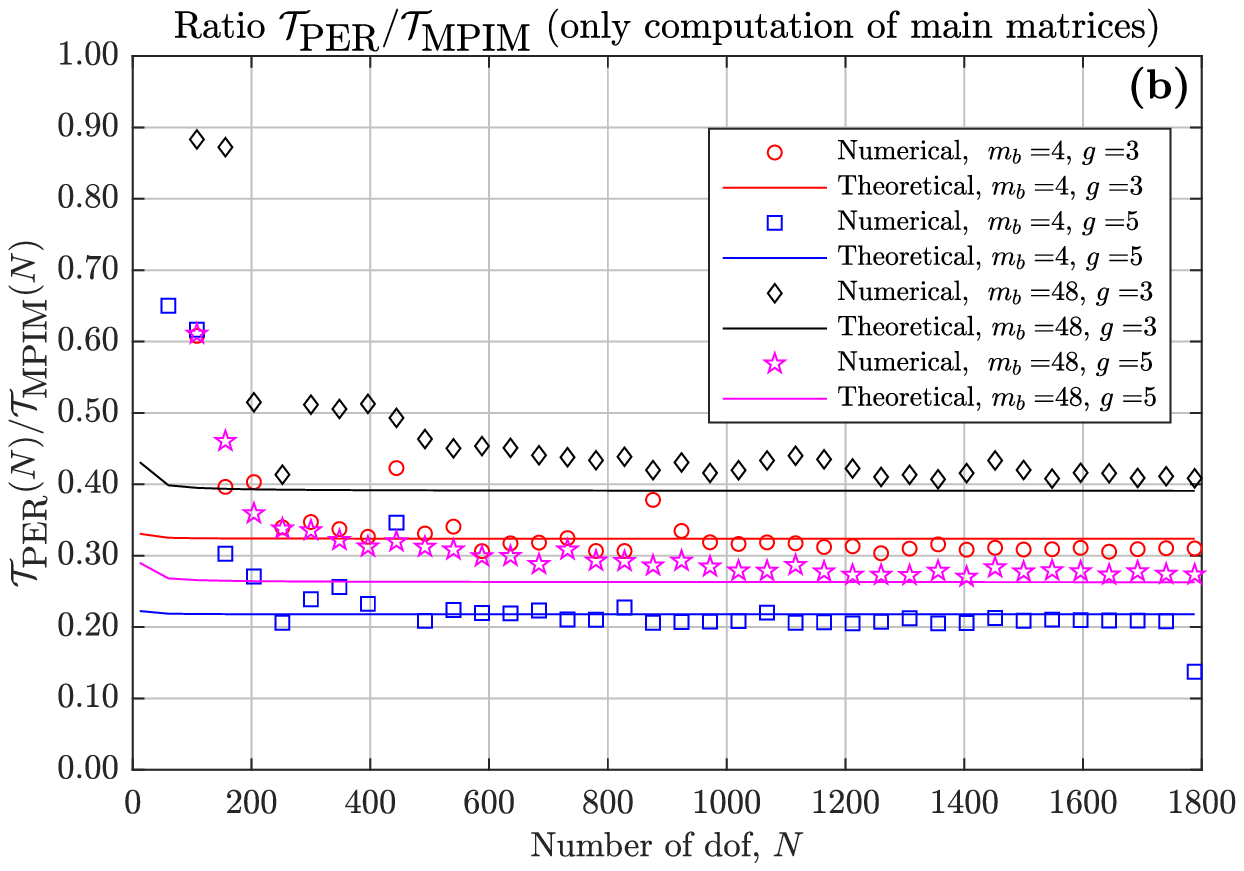} \\					
		\end{tabular}
		\caption{(a) Computation time of the main matrices for the proposed method ($\mathcal{T}_{\text{PER}}$) and for the MPIM ($\mathcal{T}_{\text{MPIM}}$) as function of the number of dof. The coefficient $\mathcal{T}_u$ represents the average time per operation, calculated by fitting curves it leads to $\mathcal{T}_u = 2.216 \times 10^{-9}$ s. (b) Ratio of computation times $\mathcal{T}_{\text{PER}} / \mathcal{T}_{\text{MPIM}}$ with number of dof $N$, for different parameters $m_b$ and $g$. Other parameters: $\zeta_A=\zeta_B=0.50$, $r_a=r_b=2$, $m_a = 2$.}
		\label{fig11}%
	\end{center}
\end{figure}

Next, the theoretical results obtained about the  computational effort of the proposed method will be validated. First, a comparison of the computational time of the main matrices of the iterative scheme in both the MPIM and the proposed method will be carried out. Let us denote by $\mathcal{T}_\text{PER}(N)$ to the computation time required to obtain both matrices  $\mathbf{a}(\dt)$ and $\mathbf{b}(\dt)$ using $N$ degrees of freedom and by $\mathcal{T}_\text{MPIM}(N)$ the same for the computation of $\bm{\mathcal{H}}(\dt)$ and $\mathbf{w}_k(\dt)$, Eqs.~\eqref{eq129} and~\eqref{eq093}, respectively. In the theoretical developments of Sec. \ref{ComputationalEffort}, closed form expressions for the number of multiplications needed to obtain these matrices, depending on the number of degrees of freedom, have been determined for both methods. The corresponding computation time arises  then from multiplying the counting expression by the unitary time per operation, $\mathcal{T}_u$ (measured in seconds). Considering only terms proportional to $N^3$, it yields
\begin{eqnarray}
\mathcal{T}_{\text{PER}} &=& \mathcal{T}_u \, \mathcal{C}_{\text{PER}} 
\approx \mathcal{T}_u \, (33 + 4(r_a + r_b) + 8p + m_b) \, N^3  \nonumber \\
\mathcal{T}_{\text{MPIM}} &=& \mathcal{T}_u \, \mathcal{C}_{\text{MPIM}} \approx 
\mathcal{T}_u \, (18 + 8 (1 + g) p) \, N^3  \label{eq133}
\end{eqnarray}
The proposed model and the MPIM have been run varying the number of dof $N$, keeping fixed the next parameters: for the proposed method $m_a = 2$, $r_a = 2$, $p=20$, $r_b = 2$, $\zeta_A=\zeta_B = 0.50$. The computation times have been monitored and plotted in fig.~\ref{fig11} as numerical simulations. From these numerical experiments, the parameter $\mathcal{T}_u$ has been obtained by curve-fitting resulting an average of $\mathcal{T}_u = 2.216 \times 10^{-9}$ seconds per operation. As it can be observed in fig.~\ref{fig11}(a), the model of Eq.~\eqref{eq133} results suitable to describe the computation time needed for the preparation of the main matrices in both methods, for the particular values of $m_b = 4$ and $g=4$. In order to check the goodness of these expressions for other cases, we have plotted in fig.~\ref{fig11}(b) the relation between computation times $\mathcal{T}_{\text{PER}}/ \mathcal{T}_{\text{PER}}$, which for a high number of degrees of freedom should tend to 
\begin{equation}
\lim_{N\to \infty} \frac{\mathcal{T}_{\text{PER}}}{\mathcal{T}_{\text{MPIM}}} = 
\frac{33 + 4(r_a + r_b) + 8p + m_b}{18 + 8 (1 + g) p}
\label{eq132}
\end{equation}
The graphical results shown in fig.~\ref{fig11}(b) show that indeed the relationship between the computation times obtained from the numerical experiments fit the model of Eq. \eqref{eq132} quite accurately, reflecting that within the values considered in practice as adequate, the preparation of the matrices is more effective with the proposed method, which in general takes less than half the time long. 
~\\

\begin{figure}[ht]%
	\begin{center}
		\begin{tabular}{cc}
			\includegraphics[width=7.5cm]{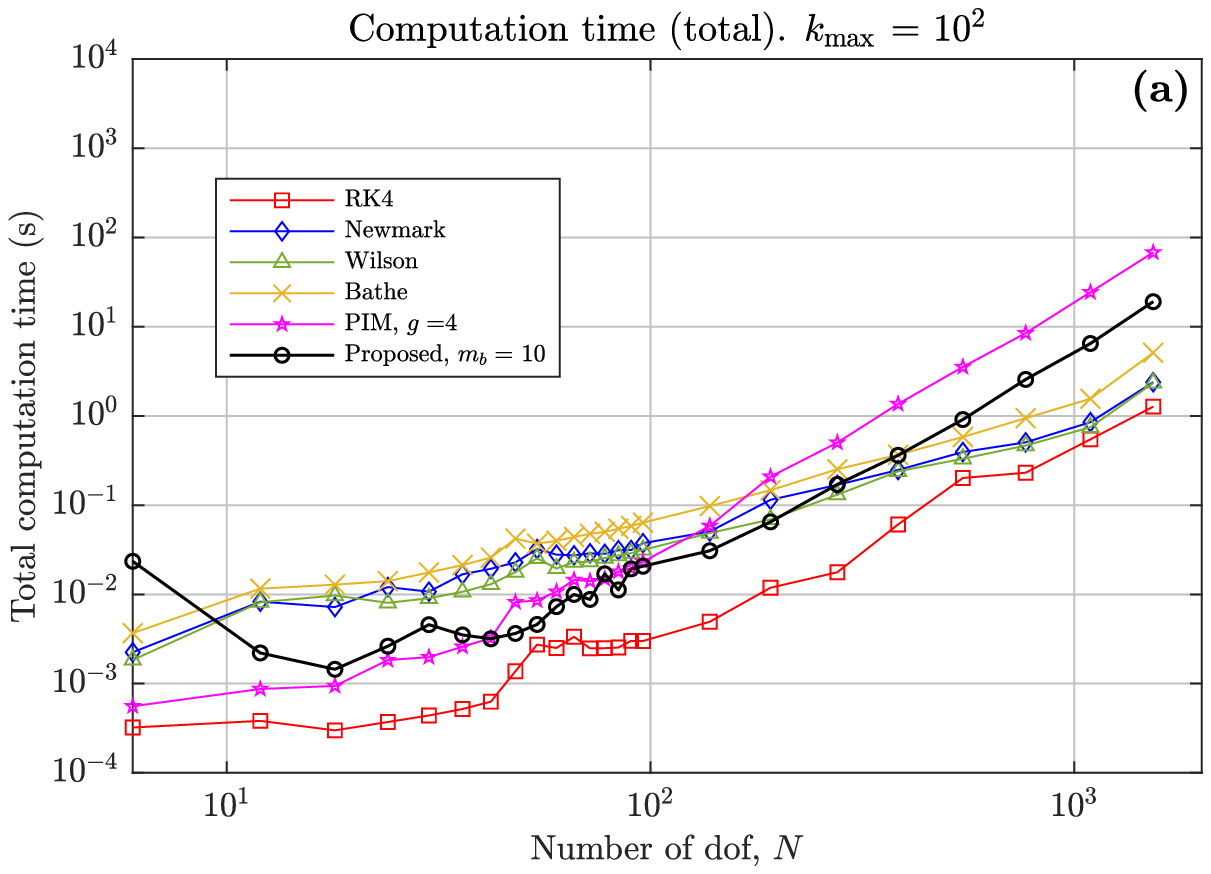} & 
			\includegraphics[width=7.5cm]{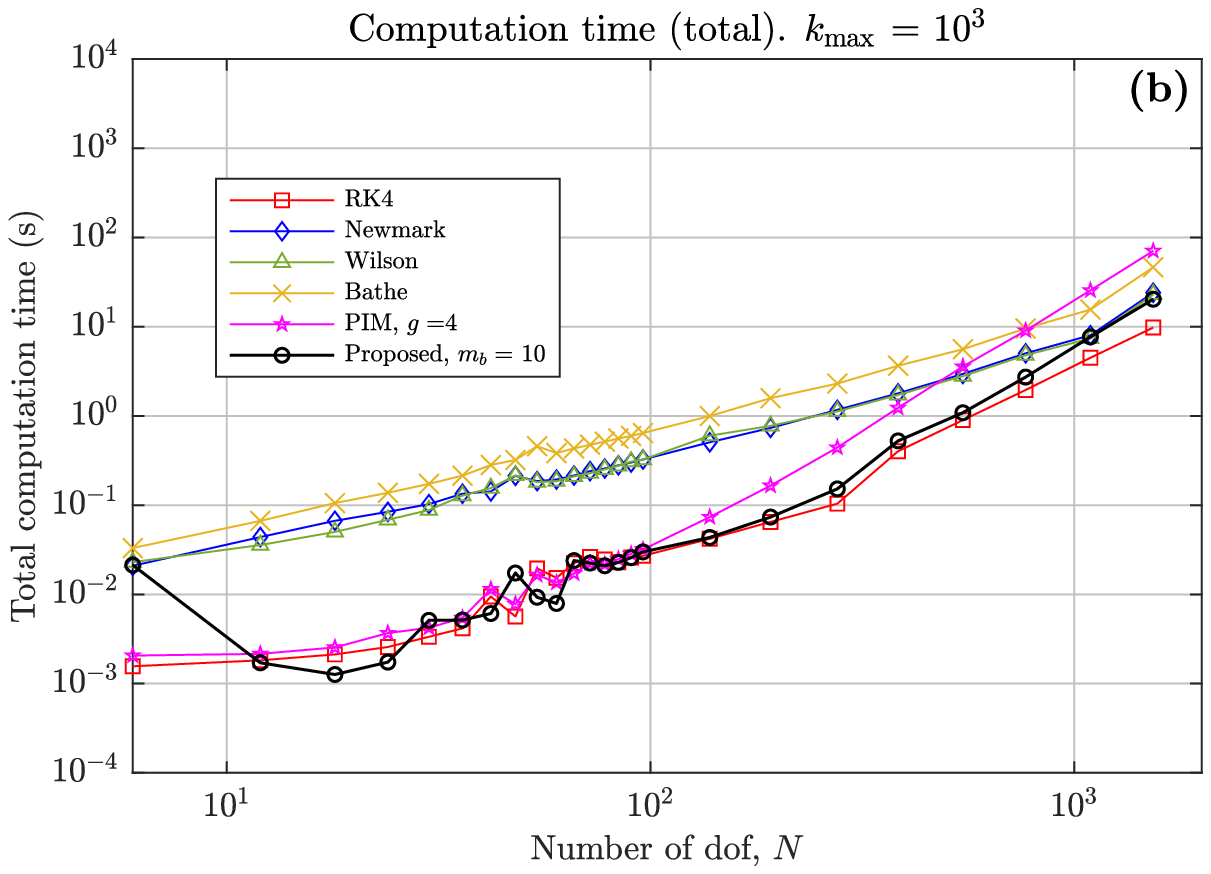}  \\ \\
			\includegraphics[width=7.5cm]{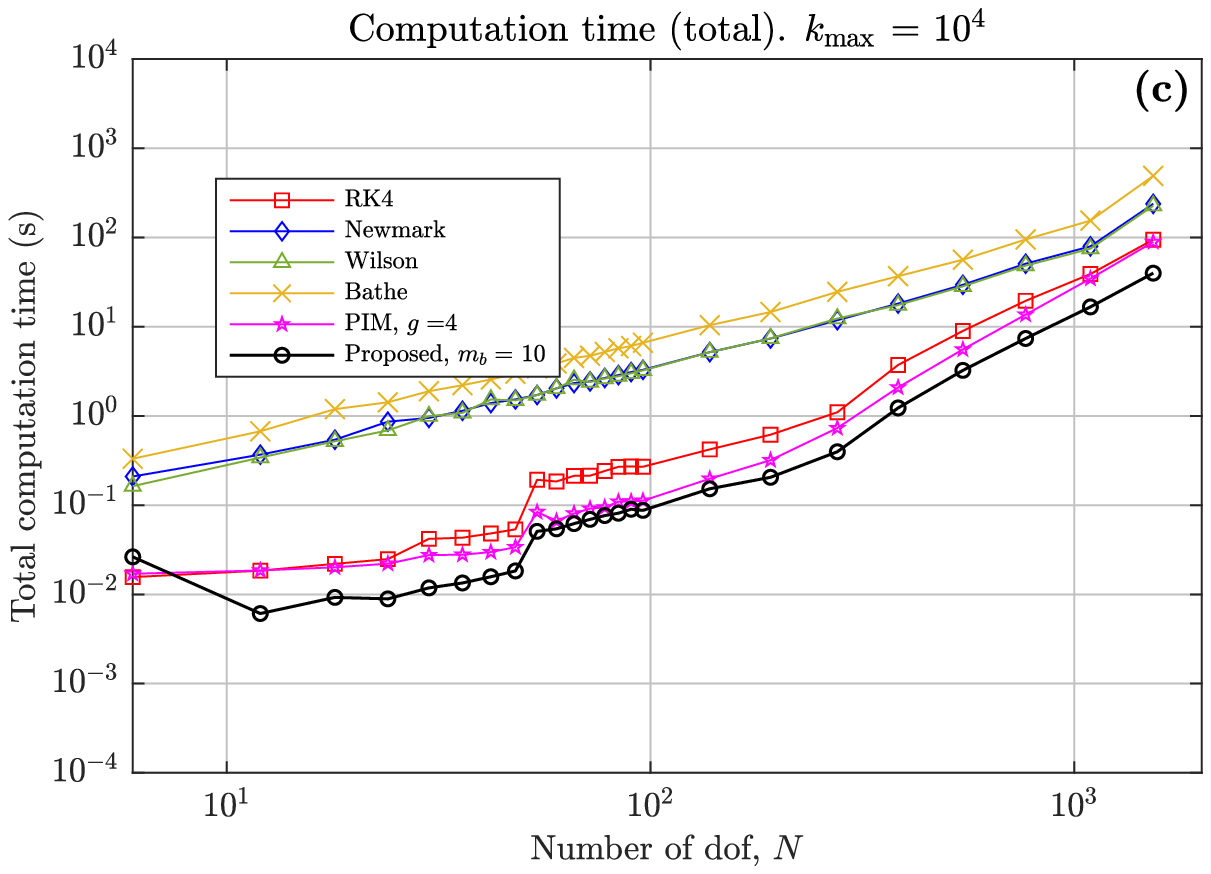} & 
			\includegraphics[width=7.5cm]{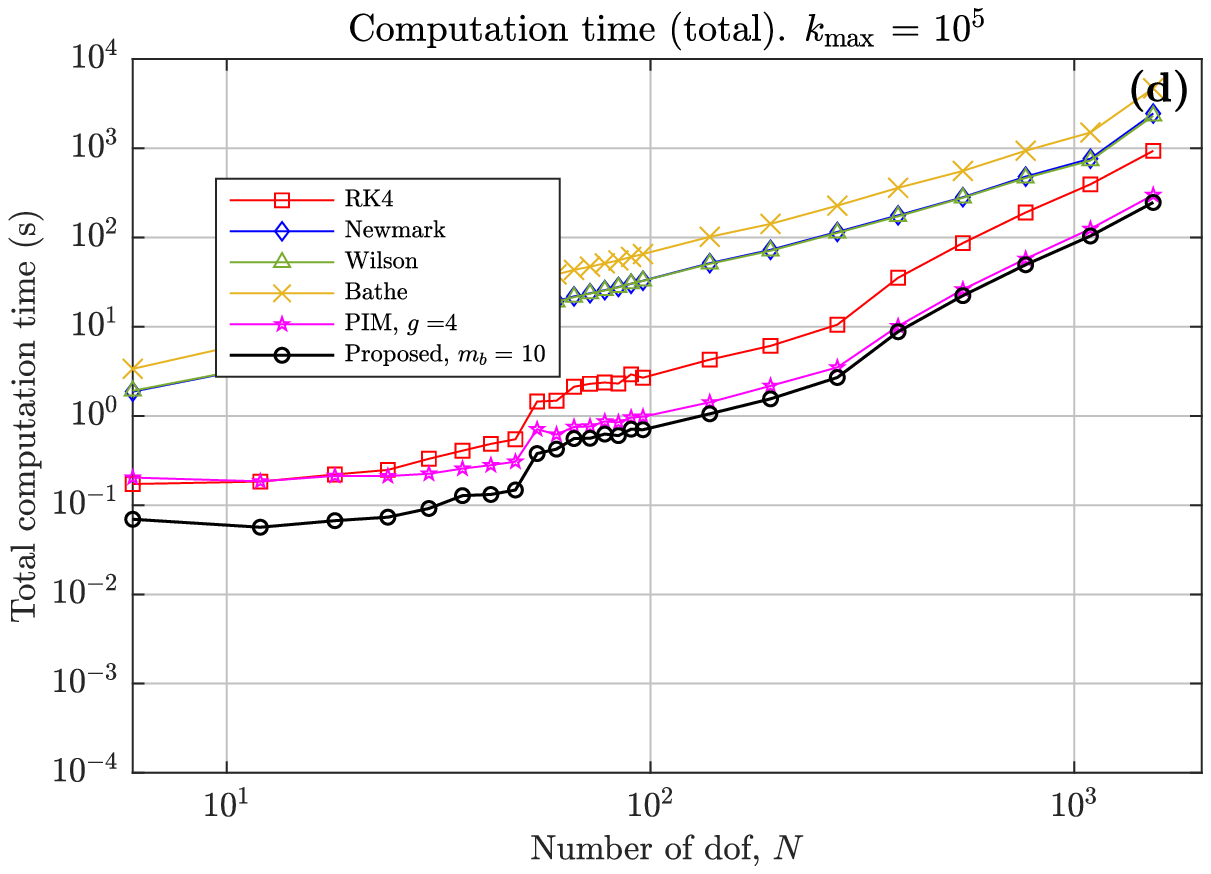} 
		\end{tabular}
		\caption{Total computation time (counting preparation of matrices and iterative process)  vs. size of the system $N$. Comparison for the different numerical methods and different sizes of time iterations, $k_{\max{}}$. (a) $k_{\max{}}=10^2$;  (b) $k_{\max{}}=10^3$; (c) $k_{\max{}}=10^4$; (d) $k_{\max{}}=10^5$.}
		\label{fig12}%
	\end{center}
\end{figure}

It is also interesting to compare the computation time considering the whole process of the algorithm, including the calculation of the model matrices and the iterative process that includes $k_{\max{}} = t_{\max{}}/ \dt$ steps. As derived from the theory, the RK4 procedure relatively consumes hardly any computational time in the model preparation, although within the iterative stage, it takes 4 times longer in each step than the MPIM and the proposed methods. This difference is expected to be noticed in simulations with a large number of iterations in relation to the system size. In fig~\ref{fig12} the total computation time for the six methods considered has been plotted as a function of the number of degrees of freedom, $N$. In addition, 4 different simulations with $k_{\max{}} = \{10^2, 10^3, 10^4, 10^5\}$ iterations have been considered. 
With only $k_{\max{}} = 10^2$ iterations, the computation of the main matrices of the algorithm takes most of the computation time, making both the MPIM and the proposed method the least efficient. Roughly, for $k_{\max{}} > 10^3$  it can be seen that the RK4 method is slower than both the proposed method and the MPIM, something that is especially noticeable in fig. \ref{fig12}(d) with a large number of iterations. In this latter plot, it can be observed that the time used for the preparation of the algorithm matrices in the MPIM and the proposed method is no longer relevant with respect to that of iterative process. It is observed that in this two methods, the computation time is similar for systems with a large number of degrees of freedom.

\section{Conclusions}

This paper considers the numerical resolution of transient problems in structural dynamics. Taking advantage that dissipative forces are in general much smaller than inertia and elastic forces, we propose to artificially perturb the damping terms to derive an iterative asymptotic procedure. The theoretical developments lead to the summation of the resulting infinite series culminating in an explicit iterative process. Convergence of the series, stability of the method, order of accuracy and computational effort are studied in depth. Detailed algorithms for the computation of the main matrices are also presented. The proposed approach has been validated by means of two numerical examples. The first one consists of a lumped-mass discrete structure. The influence of the error with different parameters such as time step, damping level and numerical parameters intrinsic to the method is studied. The results are also compared with other existing methods in the literature, showing great accuracy both in displacements and velocities for a wide range of time steps. In the second example, the transient response in a continuous beam problem has been determined. The goodness of the proposed method has been validated with respect to the other methods and times of computation have been compared, verifying that the proposed method has important computational advantages.

\section{Acknowledgments}

This work was supported by the Ministerio de Ciencia e Innovaci\'on (Spain), grant number PID2020-112759GB-I00.



\section*{References}
\bibliographystyle{elsarticle-num} 
\bibliography{bibliography}





\end{document}